\documentclass[mathpazo]{fldauth}

\usepackage{array}
\usepackage{color}
\usepackage{tabularx}
\usepackage{graphicx}
\usepackage{amsmath}
\usepackage{amssymb}
\usepackage{amsfonts}	
\usepackage{moreverb}
\usepackage{dsfont}
\usepackage{grffile}
\usepackage{bm}
\usepackage{multirow}
\usepackage{soul}

\newcommand\BibTeX{{\rmfamily B\kern-.05em \textsc{i\kern-.025em b}\kern-.08em
T\kern-.1667em\lower.7ex\hbox{E}\kern-.125emX}}

\def\volumeyear{2016}

\newcommand{\Q}{\mathbf{Q}}

\newcommand{\q}{\mathbf{q}}
\newcommand{\F}{\mathbf{F}}
\newcommand{\f}{\mathbf{f}}
\newcommand{\g}{\mathbf{g}}
\newcommand{\h}{\mathbf{h}}
\renewcommand{\v}{\mathbf{v}}
\newcommand{\B}{\mathbf{B}}

\newcommand{\p}{\textnormal{P}}

\newcommand{\Path}{\boldsymbol{\Psi}}

\newcommand{\tens}[1]{{\mathbf{#1}}}
\newcommand{\K}{\tens{K}}

\newcommand{\D}{\mathbf{D}}

\newcommand{\vv}{\vec v^{\,2}}

\renewcommand{\S}{\mathbf{S}} 
\newcommand{\A}{\mathbf{A}} 
\newcommand{\G}{\mathbf{G}} 
\newcommand{\x}{\mathbf{x}} 
\newcommand{\xxi}{\boldsymbol{\xi}} 
 
\renewcommand{\P}{\mathbf{P}} 
\newcommand{\w}{\mathbf{w}}

\newcommand{\halb}{\frac{1}{2}}

\newcommand{\be}{\begin{equation}}
\newcommand{\ee}{\end{equation}}
\newcommand{\bdm}{\begin{displaymath}}
\newcommand{\edm}{\end{displaymath}}
\renewcommand{\v}{\mathbf{v}}

\newcommand{\quotew}[1]{``#1''}

\newcommand{\dev}{\textnormal{dev}} 

\newcommand{\AAA}{{\boldsymbol{A}}}

\newcommand{\GG}{{\mathbf{G}}}

\renewcommand{\vv}{{\mathbf{v}}}

\newcommand{\II}{{\mathbf{I}}}
\newcommand{\JJ}{{\mathbf{J}}}

\newcommand{\BS}{{\boldsymbol{\sigma}}}

\newfont{\numerikEleven}{ecrm1000}
\newfont{\numerikTen}{cmss10}
\newfont{\numerikNine}{cmss9}
\newfont{\numerikEight}{cmss8}

\begin{document} 

\runningheads{W.~Boscheri et al.}{ADER-WENO-ALE schemes for nonlinear hyperelasticity}

\title{On direct Arbitrary-Lagrangian-Eulerian ADER-WENO finite volume schemes for the HPR model of nonlinear hyperelasticity}


\author{
Walter Boscheri$^1$ , 
Michael Dumbser$^1$, 
Rapha{\"e}l Loub{\`e}re$^2$} 

\address{ \centering
$^1$\ Laboratory of Applied Mathematics, Department of Civil, Environmental and \\ Mechanical Engineering,  University of Trento, I-38123 Trento, Italy. \\
$^2$ CNRS and Institut de Math\'{e}matiques de Toulouse (IMT) \\ Universit{\'e} Paul-Sabatier, Toulouse, France
          				     }	
										
\corraddr{
{\tt walter.boscheri@unitn.it} (W.~Boscheri), 
{\tt raphael.loubere@math.univ-toulouse.fr} (R.~Loub{\`e}re), 
 {\tt michael.dumbser@unitn.it} (M.~Dumbser). }

\begin{abstract}
This paper is concerned with the numerical solution of the \textit{unified} first order \textit{hyperbolic} formulation of 
\textit{continuum mechanics} proposed by Peshkov \& Romenski \cite{PeshRom2014} (HPR model), which is based on the theory of nonlinear 
\textit{hyperelasticity} of Godunov \& Romenski \cite{GodunovRomenski72,Godunov:2003a}.  
Notably, the governing PDE system is \textit{symmetric hyperbolic} and fully consistent with the first and the second principle of thermodynamics. The nonlinear system of 
governing equations of the HPR model is large and includes stiff source terms as well as non-conservative products.
In this paper we solve this model for the first time on \textit{moving unstructured meshes} in multiple space dimensions by employing high order 
accurate one-step ADER-WENO finite volume schemes in the context of cell-centered direct Arbitrary-Lagrangian-Eulerian (ALE) algorithms. 

The numerical method is based on a WENO polynomial reconstruction operator on moving unstructured meshes, a fully-discrete one-step ADER 
scheme that is able to deal with stiff sources \cite{DumbserEnauxToro}, a nodal solver with relaxation to determine the mesh motion, 
and a path-conservative technique of Castro \& Par\'es for the treatment of non-conservative products \cite{Pares2006,Castro2006}.  
We present numerical results obtained by solving the HPR model with ADER-WENO-ALE schemes in the stiff relaxation limit, 
showing that fluids (Euler or Navier-Stokes limit), as well as purely elastic or elasto-plastic solids can be simulated in the 
framework of nonlinear hyperelasticity with the \textit{same} system of governing PDE. 
The obtained results are in good agreement when compared to exact or numerical reference solutions available in the literature.
\end{abstract}

\keywords{
 high order direct Arbitrary-Lagrangian-Eulerian finite volume schemes, 
 hyperbolic Peskhov \& Romenski model (HPR model), 
 nonlinear hyperelasticity, 
 stiff source terms,
 non-conservative products,
 ADER-WENO schemes on unstructured meshes, 
 high order of accuracy in space and time, 
 hyperbolic conservation laws, 
 fluid mechanics and solid mechanics,
 continuum mechanics}

\maketitle
\vspace*{-1cm}
\setcounter{tocdepth}{2}


\section{Introduction} \label{sec:introduction}

The aim of this paper is the numerical solution of the \textit{unified} first order hyperbolic formulation of \textit{continuum mechanics} proposed by Peshkov \& Romenski \cite{PeshRom2014}, denoted as \textbf{HPR model} in the following, which is based on the theory of nonlinear \textbf{hyperelasticity} of Godunov and Romenski \cite{GodunovRomenski72,Godunov:2003a}, and which describes \textbf{fluid mechanics} and 
\textbf{solid mechanics} at the same time in \textit{one single system} of governing partial differential equations (PDE). In the HPR model the viscous stresses are computed from the so-called distortion tensor $\AAA$, which is one of the primary state variables in this first order system. The appealing property of the HPR model is its ability to describe \textit{within the same mathematical framework}  the behavior of inviscid and viscous compressible Newtonian and non-Newtonian fluids with heat conduction, and, at the same time, the behavior of elastic and elasto-plastic solids. In this model fluids as well as solids are modeled via a stiff source term that accounts for strain relaxation in the evolution equations of the distortion tensor. In addition, heat conduction is included using a first order hyperbolic evolution equation of the thermal impulse which allows the heat flux to be retrieved. The governing system of PDEs is \textit{symmetric hyperbolic} and fully consistent with the first and the  second principle of thermodynamics, as detailed in \cite{PeshRom2014,Dumbser_HPR_16}. However, this system has a large number of equations, is nonlinear and it includes stiff source terms and also non-conservative products.

Consequently, the numerical solution of such a large multi-dimensional system on moving meshes is a big challenge. For this purpose, in this work we propose to employ a high order accurate multi-dimensional ADER-WENO finite volume scheme in the context of direct Arbitrary-Lagrangian-Eulerian (ALE) algorithms. This scheme is constructed with a high order WENO polynomial reconstruction operator 
on unstructured meshes \cite{DumbserKaeser06b,Dumbser2007204}, a one-step space-time ADER integration \cite{titarevtoro,titarevtoro2,Toro:2006a} that is suitably extended for dealing with stiff 
sources \cite{DumbserEnauxToro,HidalgoDumbser}, a nodal solver with relaxation to determine the mesh motion \cite{MaireRezoning,Lagrange2D,Lagrange3D,LagrangeMHD}, and a path-conservative integration technique for the treatment of non-conservative products, 
following the ideas of Castro \& Par\'es \cite{Pares2006,Castro2006}, which have been recently extended to the  
moving-mesh framework in \cite{LagrangeNC,Lagrange3D,ALE-MOOD2}. The proper treatment of boundary conditions is of paramount importance for these simulations on moving meshes. We will pay special attention to them in this work.

In this paper we intend to show that, although the HPR model may seem to be more complex and difficult to solve than other classical ones (Euler \& Navier-Stokes equations, linear elasticity or 
nonlinear hypo-elasticity with plastic strain), the high order ADER-WENO-ALE schemes which allow for a proper treatment of non-conservative terms and stiff source terms \cite{LagrangeNC,Lagrange3D} 
are an appropriate candidate for this task. Therefore, we will present numerical results obtained with ADER-WENO-ALE schemes for the HPR model in the stiff relaxation limits showing that fluids 
(Euler or Navier-Stokes limits) as well as pure elastic and elasto-plastic solids can be simulated. In these different situations --- fluids, elastic and elasto-plastic solids --- which usually 
require a different mathematical model for each situation, we will numerically prove that the high order accurate ADER-WENO-ALE algorithm is able to reproduce existing exact or numerical 
reference solutions even for very demanding test cases. These test problems involve shocks (viscous or inviscid ones), contacts and rarefactions in fluids, along with reversible or irreversible 
deformations in elasto-plastic solids.

%
%
The rest of this paper is organized as follows. Section \ref{sec:HPR} introduces the unified first order hyperbolic Peshkov-Romenski (HPR) model of continuum mechanics, which 
is numerically solved in this paper. Section \ref{sec:framework} presents the high order accurate ADER-WENO-ALE schemes devoted to solve general hyperbolic systems of conservation laws 
with stiff source terms and non-conservative products. Boundary conditions are discussed in Section \ref{sec.BCs}, while numerical experiments are carried out in Section \ref{sec:numerics}, 
which also contains a detailed description of these test cases, as well as the obtained numerical results with associated comments. 
Note that the numerical experiments are designed so that the scheme solves two extreme limits of the HPR model, namely inviscid and viscous fluids (i.e. the compressible Euler 
equations for gasdynamics and the compressible Navier-Stokes equations) as well as elastic / elasto-plastic solids. Conclusions and perspectives are proposed in the last Section 
\ref{sec:conclusion}. 


\section{The HPR model: a unified first order hyperbolic approach to continuum mechanics} \label{sec:HPR}

In this work we consider the so-called \textbf{H}yperbolic \textbf{P}eshkov-\textbf{R}omenski (HPR) model 
\cite{PeshRom2014}, which is the first successful attempt to build a \textit{unified} formulation of \textit{continuum mechanics} 
under a first order \textit{symmetric hyperbolic} form that includes classical \textit{fluid mechanics} and \textit{solid mechanics} just as two special limiting cases of the same formulation.
We refer to the recent work of Dumbser et al. \cite{Dumbser_HPR_16}, where a detailed introduction to this model is given and where the HPR model has been 
solved numerically for the first time using high order accurate \textit{Eulerian} ADER-WENO and ADER-DG schemes on \textit{fixed grids}, and where many numerical examples 
have been provided.  
The HPR model also includes a hyperbolic formulation of heat conduction and it can be written under the form given in \cite{Dumbser_HPR_16} as follows: 
\begin{subequations}\label{eqn.HPR}
	\begin{align}
	& \frac{\partial \rho}{\partial t}+\frac{\partial \rho v_k}{\partial x_k}=0,\label{eqn.conti}\\[2mm]
	&\displaystyle\frac{\partial \rho v_i}{\partial t}+\frac{\partial \left(\rho v_i v_k + p \delta_{ik} - \sigma_{ik} \right)}{\partial x_k}=0, \label{eqn.momentum}\\[2mm]
	&\displaystyle\frac{\partial A_{i k}}{\partial t}+\frac{\partial A_{im} v_m}{\partial x_k}+v_j\left(\frac{\partial A_{ik}}{\partial x_j}-\frac{\partial A_{ij}}{\partial x_k}\right)
	=-\dfrac{ \psi_{ik} }{\theta_1(\tau_1)},\label{eqn.deformation}\\[2mm]
	&\displaystyle\frac{\partial \rho J_i}{\partial t}+\frac{\partial \left(\rho J_i v_k+ T \delta_{ik}\right)}{\partial x_k}=-\dfrac{\rho H_i}{\theta_2(\tau_2)}, \label{eqn.heatflux}\\[2mm]
	&\displaystyle\frac{\partial \rho s}{\partial t}+\frac{\partial \left(\rho s v_k + H_k \right)}{\partial x_k}=\dfrac{\rho}{\theta_1(\tau_1) T} \psi_{ik} \psi_{ik} + \dfrac{\rho}{\theta_2(\tau_2) T} H_i H_i \geq0. \label{eqn.entropy}
	\end{align}
\end{subequations}
The solutions of the above PDE system fulfill also the additional conservation of total energy  
\begin{equation}\label{eqn.energy}
\frac{\partial \rho  E}{\partial t}+\frac{\partial \left(v_k \rho  E + v_i (p \delta_{ik} - \sigma_{ik}) + q_k \right)}{\partial x_k}=0.
\end{equation}
At this point we emphasize that the system above is an \textit{overdetermined} system of PDE, hence in the numerical solution of the above model we solve the \textit{total energy conservation equation}  
\eqref{eqn.energy} and \textit{not} the \textit{entropy equation} \eqref{eqn.entropy}. Such a choice is mandatory for overdetermined systems. 
We use the following notation: 
 $\rho$ is the mass density,  
 $ [v_i]=\vv=(u,v,w) $ is the velocity vector, 
 $ [A_{ik}]=\AAA $ is the distortion tensor, 
 $ [J_i]=\JJ $ is the thermal impulse vector, 
 $ s $ is the entropy, 
 $ p = \rho^2E_{\rho} $ is the pressure, 
 $ E =E(\rho,s,\vv,\AAA,\JJ)$ is the total energy potential, 
 $ \delta_{ik} $ is the Kronecker delta, 
$ [\sigma_{ik}]=\boldsymbol{\sigma} = - [\rho A_{mi}E_{A_{mk}}] $ is the symmetric viscous shear stress tensor, 
$ T=E_s $ is the temperature, 
$[q_k] = \mathbf{q} = [E_s E_{J_k}]$  is the heat flux vector and 
$ \theta_1 = \theta_1(\tau_1) > 0$ and $ \theta_2 = \theta_2(\tau_2) > 0$ are positive scalar functions depending on the strain dissipation time $ \tau_1 > 0$ and the thermal impulse relaxation time $ \tau_2 > 0$, respectively. 
The dissipative terms $\psi_{ik}$ and $H_i$ on the right hand side  
of the evolution equations for $\AAA$, $\JJ$ and $s$ are defined as $[\psi_{ik}] = \boldsymbol{\psi} = [E_{A_{ik}}]$ and $[H_i] = \mathbf{H} = [E_{J_i}]$, respectively. 
Accordingly, the viscous stress tensor and the heat flux vector are directly related to the dissipative terms on the right hand side via 
$ \boldsymbol{\sigma} = - \rho \AAA^T \boldsymbol{\psi} $ and $ \mathbf{q} = T \, \mathbf{H}$.     
Note that $ E_\rho $, $ E_s $, $ E_{A_{ik}} $ and $ E_{J_i} $ denote the partial derivatives $ \partial E/\partial \rho$, $ \partial E/\partial s $,  
$ \partial E/\partial A_{ik}$ and $ \partial E/\partial J_i$; they are the \textit{energy gradients in the state space} or the \textit{thermodynamic forces}. 
The Einstein summation convention over repeated indices is implied throughout this paper.

These equations express 
the mass conservation (\ref{eqn.conti}),
 the momentum conservation~(\ref{eqn.momentum}),
 the time evolution for the distortion tensor~(\ref{eqn.deformation}),
 the time evolution for the thermal impulse~(\ref{eqn.heatflux}), 
 the time evolution for the entropy~(\ref{eqn.entropy}), 
and the total energy conservation~(\ref{eqn.energy}). 
The PDE governing the time evolution of the thermal impulse~\eqref{eqn.heatflux} looks similar to the momentum equation \eqref{eqn.momentum}, where the temperature $T$ takes the role of the pressure $p$.  Therefore we refer to this equation as the thermal momentum equation.

To close the above system, the total energy potential $ E(\rho,s,\vv,\AAA,\JJ) $ must be specified. 
This potential definition will then generate all constitutive fluxes (\textit{i.e.} non advective fluxes)  
and source terms by means of its partial derivatives with respect to the state variables. 
As a consequence the energy potential specification is fundamental for the model formulation.

In order to specify $ E $, following \cite{PeshRom2014,Dumbser_HPR_16} we note that there are three scales: 
the molecular scale, referred to as the  \textit{microscale};
 the scale of the material elements, called here \textit{mesoscale}; 
 and the main flow scale, that is the \textit{macroscale}. 
As a consequence it is assumed that the total energy $ E $ is decomposed into three terms,
each of them representing the energy distributed in its corresponding scale, that is: 
\begin{equation}\label{eq:total_energy}
E(\rho,s,\vv,\AAA,\JJ)=E_{1}(\rho,s)+E_{2}(\AAA,\JJ) + E_{3}(\vv).
\end{equation}
The specific \textit{kinetic energy} per unit mass $ E_3(\vv) =\dfrac{1}{2} v_iv_i$ refers to the macroscale part of the total energy. 
The \textit{internal energy} $ E_1(\rho,s) $ is related to the kinetic energy of the molecular motion and it is sometimes referred to as the \textit{equilibrium energy} because it is the only energy which does not disappear in the thermodynamic equilibrium when meso- and macro-scopic dynamics are absent, but only molecular dynamics remains. In this paper, for $ E_1 $ we will use either the \textit{ideal gas equation of state} 
\begin{equation}\label{eq:ideal_gas_eos}
E_1(\rho,s)=\frac{c_0^2}{\gamma(\gamma-1)},\ \ c_0^2=\gamma\rho^{\gamma-1}e^{s/c_v},
\end{equation}
%
or the \textit{Mie-Gr{\"u}neisen equation of state}
\begin{equation}\label{eq:mie_gruneisen_eos}
E_1(\rho,p)= \frac{p-\rho_0 c_0^2 \ f(\nu)}{\rho_0 \Gamma_0}, \ f(\nu) = \frac{(\nu-1)(\nu-\frac12\Gamma_0(\nu-1))}{(\nu - s(\nu-1))^2}, \ \nu=\frac{\rho}{\rho_0},
\end{equation}
where $c_0$ has the meaning of the adiabatic sound speed, $ c_v $ and $c_p$ are the specific heat capacities at constant volume and at constant pressure, respectively, 
which are related by the ratio of specific heats $ \gamma=c_p/c_v $. 
Moreover $ \rho_0 $ is the reference mass density and $ p_0 $ is the reference 
(atmospheric) pressure. For the mesoscopic, or \textit{non-equilibrium}, part of the total energy we adopt a simple quadratic form
\begin{equation}\label{eq:e_2}
E_2(\AAA,\JJ)=\dfrac{c_s^2}{4}G^{\rm TF}_{ij}G^{\rm TF}_{ij}+\frac{\alpha^2}{2}J_i J_i,
\end{equation}
with 
\begin{equation}
 [G_{ij}^{\rm TF}] = \dev(\GG) = \GG-\frac{1}{3} {\rm tr}(\GG) \II,  \qquad \textnormal{ and } \qquad \GG=\AAA^\mathsf{T}\AAA. 
\end{equation} 
Here, $[G_{ij}^{\rm TF}] = \dev(\GG)$ is the deviator, or the \textit{trace-free} part, of the tensor $\GG=\AAA^\mathsf{T}\AAA$ and ${\rm tr}(\GG)=G_{ii}$ is its trace, 
$ \II $ is the unit tensor and $ c_s $ is the characteristic velocity of propagation of transverse perturbations. 
In the following we shall refer to it as the \textit{shear sound velocity}. The characteristic velocity of heat wave propagation $c_h$ is related to the variable $\alpha$.

The fundamental \textit{frame invariance principle} implies that the total energy can only depend on vectors and tensors by means of their invariants. Hence, 
\[G_{ij}^{\rm TF} G_{ij}^{\rm TF}\equiv I_2-I_1^2/3,\]
where $ I_1={\rm tr}(\GG) $ and $ I_2={\rm tr}(\GG^2) $, and therefore $ E_2 $, as well as the total energy $ E $, are only a function of invariants of $ \AAA $ and $ \JJ $.



The algebraic source term on the right-hand side of equation~(\ref{eqn.deformation}) describes the shear strain dissipation due to material element rearrangements, and the source term 
on the right-hand side of \eqref{eqn.heatflux} describes the relaxation of the thermal impulse due to heat exchange between material elements. Once the total energy potential is specified, all fluxes and source terms have an explicit form. 
Thus, for the energy $ E_2(\AAA,\JJ) $ given by \eqref{eq:e_2}, 
we have $\boldsymbol{\psi}= E_{\AAA}= c_s^2 \AAA \dev(\GG)$, hence the shear stresses are 
explicitly given by
\begin{equation}
\label{eqn.stress} 
\boldsymbol{\sigma}= -\rho\AAA^\mathsf{T} \boldsymbol{\psi} = -\rho\AAA^\mathsf{T} E_{\AAA} = -\rho c_s^2 \GG \dev(\GG),\ \ \qquad {\rm tr}(\boldsymbol{\sigma})=0,
\end{equation}
and the strain dissipation source term becomes
\begin{equation}
\label{eqn.psi} 
-\dfrac{\boldsymbol{\psi}}{\theta_1(\tau_1)} = -\dfrac{E_{\AAA}}{\theta_1(\tau_1)}=-\dfrac{3}{\tau_1 } \left| \AAA \right|^{\frac{5}{3}} \AAA \dev(\GG),
\end{equation}
where we have chosen $ \theta_1(\tau_1) = \tau_1 c_s^2 / 3 \, |\AAA|^{-\frac{5}{3}} $, with $|\AAA|=\det(\AAA) > 0$ the determinant of $\AAA$ and $\tau_1$ being 
the strain relaxation time, also called the particle-settled-life (PSL) time in \cite{Frenkel1955,PeshRom2014}.
In other words, this time scale characterizes how long a material element is connected with its neighbor elements before  rearrangement occurs. The determinant of $\AAA$ must  satisfy the \textit{constraint}    
\begin{equation} 
\label{eqn.compatibility} 
 |\AAA| = \frac{\rho}{\rho_0},    
\end{equation} 
where $\rho_0$ is the density at the reference configuration, see \cite{PeshRom2014}. 
Furthermore, from the energy potential $E_2(\AAA,\JJ)$ the heat flux vector follows from $E_{\JJ} = \alpha^2 \JJ$
as
\begin{equation}
\label{eqn.hyp.heatflux}
 \mathbf{q} = T \, \mathbf{H} = E_s E_{\JJ} = \alpha^2 T \JJ.    
\end{equation} 
For the thermal impulse relaxation source term, we postulate that $ \theta_2=\tau_2\alpha^2 \frac{\rho}{\rho_0} \frac{T_0}{T}$ yielding
\begin{equation}
-\dfrac{\rho \mathbf{H}}{\theta_2(\tau_2)} = -\dfrac{\rho E_{\JJ}}{\theta_2(\tau_2)}= - \frac{T}{T_0} \frac{\rho_0}{\rho} \dfrac{\rho\JJ}{\tau_2}.
\end{equation}
The previous formula 
contains another characteristic relaxation time $\tau_2$ which is associated to heat conduction.  
The motivation for this particular choice of $\theta_1$ and $\theta_2$
is the connection with classical Navier-Stokes-Fourier theory in the stiff limit $\tau_1 \to 0$ and $\tau_2 \to 0$,
see  \cite{Dumbser_HPR_16} for details. 

As shown in \cite{PeshRom2014,Dumbser_HPR_16}, the HPR model is compatible with the first and second law of thermodynamics and it constitutes a \textit{hyperbolic} system of PDEs. 
For a detailed discussion of the hyperbolicity of nonlinear hyperelasticity, see \cite{Ndanou2014}. For a discussion on the symmetric hyperbolic structure, see 
\cite{Dumbser_HPR_16} and references therein. In other words the Cauchy problem for the system \eqref{eqn.HPR} is well-posed. A detailed discussion of the intrinsic nature of this model can be found in \cite{PeshRom2014,Dumbser_HPR_16} and we refer the interested reader to these references. Further work on nonlinear hyperelasticity can be found e.g. in 
\cite{GodunovRomenski72,Godunov:1995a,Godunov:2003a,Kluth2010,Pesh2010,Pesh2015,BartonRomenski2010,Barton2013,BartonRom2012}. 
In this paper we assume the model as given and our goal is to solve it numerically in an accurate, robust and efficient way on moving unstructured meshes using one of the most advanced  high order accurate ADER-WENO direct Arbitrary-Lagrangian-Eulerian schemes that it currently available \cite{Lagrange1D,Lagrange2D,LagrangeNC,Lagrange3D}.


\section{High order accurate direct ADER-WENO-ALE schemes for hyperbolic PDE} \label{sec:framework}

As already mentioned, the HPR model is a large nonlinear system of hyperbolic balance laws which contains
 non-conservative products and stiff source terms.
To solve this system we consider the arbitrary high order accurate ADER-WENO direct Arbitrary-Lagrangian-Eulerian (ALE) 
finite volume schemes derived in \cite{Lagrange2D,LagrangeNC,Lagrange3D} that we refer to as ADER-WENO-ALE in the rest of the paper.
The HPR model \eqref{eqn.HPR} can be cast into the following general formulation which holds in multiple space dimensions $d\in[2,3]$: 
\begin{equation}
\label{eqn.pde.nc} 
  \frac{\partial \Q}{\partial t} + \nabla \cdot \tens{\F}(\Q) + \B(\Q)\cdot\nabla\Q = \S(\Q) , \qquad \x \in \Omega \subset \mathds{R}^d, t \in \mathds{R}_0^+, 
\end{equation} 
where $\Q=(q_1,q_2,\cdots , q_v)$ is the vector of conserved variables, $\tens{\F} = (\mathbf{f}, \mathbf{g}, \mathbf{h})$ is the conservative nonlinear flux tensor, $\B=(\B_1,\B_2,\B_3)$ is the purely non-conservative part of the system written in block-matrix notation and $\S ( \Q )$ is the vector of algebraic source terms. We furthermore introduce the abbreviation $\P = \P(\Q, \nabla \Q) = \B(\Q)\cdot\nabla\Q$ to 
simplify the notation in some parts of the manuscript.

In our moving-mesh framework the computational domain $\Omega(t) \subset \mathds{R}^d$  is discretized at any time level $t^n$ by a set of moving and deforming simplexes $T^n_i$. $N_E$ denotes the total number of elements and the union of all elements is referred to as the mesh configuration $\mathcal{T}^n_{\Omega}$ of the domain: 
$
\mathcal{T}^n_{\Omega} = \bigcup \limits_{i=1}^{N_E}{T^n_i}. 
$
We assume that the computational domain continuously changes in time. Because of this fundamental
assumption we adopt the mapping between the physical element $T^n_i$ to the reference element $T_e$ defined in the reference coordinate system $\boldsymbol{\xi}=(\xi,\eta,\zeta)$. As usual, the reference element $T_e$ is taken to be the unit triangle in 2D or the unit tetrahedron in 3D, see \cite{Lagrange2D,Lagrange3D}. 

For any finite volume scheme, data are represented by piecewise constant cell averages both in space and time. As a consequence we define at each time level $t^n$ within the control volume $T^n_i$ the mean value of the state vector $\Q_i^n$ as
\begin{equation}
  \Q_i^n = \frac{1}{|T_i^n|} \int_{T^n_i} \Q(\mathbf{x},t^n) \ d\x,     
 \label{eqn.cellaverage}
\end{equation}  
where $|T_i^n|$ is the volume of element $T_i^n$. 
High order of accuracy in space is obtained by means of a polynomial reconstruction technique that provides piecewise high 
order WENO polynomials $\mathbf{w}_h(\x,t^n)$ from the known cell averages $\Q_i^n$ (see next Section \ref{sec.reconstruction}).
High order of accuracy in time is further achieved by applying a local space-time discontinuous Galerkin predictor method starting from the high accurate WENO reconstruction polynomials $\mathbf{w}_h(\x,t^n)$ (see Section \ref{sec.lst}). Both techniques are now introduced.

\subsection{Polynomial reconstruction} 
\label{sec.reconstruction} 

\subsubsection{Single stencil reconstruction.}
The reconstruction operator generates piecewise polynomials $\mathbf{w}_h(\x,t^n)$ of degree $M$ which are computed for each element $T^n_i$ considering the so-called reconstruction stencil $\mathcal{S}_i$ and its associated known cell averages. The reconstruction stencil $\mathcal{S}_i$ is composed of a number $n_e$ of neighbor elements of $T^n_i$, which is bigger than the smallest number 
\begin{equation} 
\mathcal{M} = \mathcal{M}(M,d) = \prod \limits_{k=1}^d \frac{(M+k)}{d!}, 
\end{equation} 
 needed to reach the nominal order of accuracy $M+1$ in $d$ space dimensions, according to \cite{StencilRec1990,Olliver2002,friedrich,kaeserjcp,DumbserKaeser06b}. As suggested in \cite{DumbserKaeser06b,Dumbser2007204}, for an unstructured mesh we usually take $n_e = d \cdot \mathcal{M}$, with $d\in[2,3]$ representing the number of space dimensions. 
The stencil called $\mathcal{S}_i$ is defined as
$\mathcal{S}_i = \bigcup \limits_{j=1}^{n_e} T^n_{m(j)}$,
where $1 \leq j \leq n_e$ is a local index counting the elements in the stencil and $m(j)$ is a mapping from the local index $j$ to the global index of the element in $\mathcal{T}^n_{\Omega}$. We rely on the orthogonal Dubiner-type basis functions $\psi_l(\xi,\eta,\zeta)$ \cite{Dubiner,orth-basis,CBS-book}, defined on the reference element $T_e$, to explicitly write the high order accurate reconstructed polynomial as
\begin{equation}
\label{eqn.recpolydef} 
\w_h(\x,t^n) = \sum \limits_{l=1}^\mathcal{M} \psi_l(\boldsymbol{\xi}) \hat \w^{n}_{l,i} := \psi_l(\boldsymbol{\xi}) \hat \w^{n}_{l,i},   
\end{equation}
where the mapping from $\x$ to the reference coordinate system $\boldsymbol{\xi}$ is considered and the $\hat \w^{n}_{l,i}$ denote the unknown degrees of freedom, also called expansion coefficients.
The procedure to determine the degrees of freedom demands the 
integral conservation for the reconstruction on each element $T_j^n$ belonging to stencil $\mathcal{S}_i$, 
that is
\begin{equation}
\label{intConsRec}
\frac{1}{|T^n_j|} \int \limits_{T^n_j} \psi_l(\boldsymbol{\xi}) \hat \w^{n}_{l,i} d\x = \Q^n_j, \qquad \forall T^n_j \in \mathcal{S}_i.     
\end{equation}
The above relations \eqref{intConsRec} yield an \textit{overdetermined} linear system of equations for the unknowns $\hat \w^{n}_{l,i}$ that can be solved using either a constrained 
least squares technique (LSQ), see \cite{DumbserKaeser06b}, or a more sophisticated singular value decomposition (SVD) algorithm \cite{DumbserKaeser06b,ADER_MOOD_14}.  

\subsubsection{WENO procedure.}
As stated by the Godunov theorem \cite{godunov}, linear monotone schemes are at most of order one and if the scheme is required to be high order accurate and non-oscillatory, 
it must be nonlinear. 
In this work we consider the pragmatic polynomial WENO approach that has also been adopted in 
\cite{friedrich,kaeserjcp,DumbserKaeser06b,Dumbser2007204,AboiyarIske,MixedWENO2D,MixedWENO3D,LagrangeNC,Lagrange2D,Lagrange3D,LagrangeNC,LagrangeQF,LagrangeMHD,LagrangeMDRS} 
to supplement the linear polynomial reconstruction procedure previously described with a nonlinearity. For optimal WENO schemes, see 
\cite{balsarashu,HuShuTri,shi,ZhangShu3D,WENOsubcell,Semplice2016,Cravero2015}. 
Seven or nine reconstruction stencils are first determined for $d=2$ and $d=3$, respectively, and they are further used to compute the associated different polynomials for each cell of the computational domain. These stencils are supposed to cover sufficiently enough \quotew{directions}  in order to \quotew{catch} local oscillatory phenomena. Next, these seven or nine polynomials are blended together using nonlinear weights to obtain the actual high order WENO polynomials $\mathbf{w}_h(\x,t^n)$.
This rather classical procedure has already been described in \cite{DumbserKaeser06b,Dumbser2007204,Lagrange2D,Lagrange3D} and in all the aforementioned references, consequently we omit the details 
in this paper. However we highly recommend the interested readers to consult these references.

\subsection{Local space-time Discontinuous Galerkin predictor on moving curved meshes} 
\label{sec.lst} 
The reconstructed polynomials $\w_h(\x,t^n)$ computed at time $t^n$ are then \textit{evolved} during one time step \textit{locally} within each element $T_i(t)$, without needing any neighbor information, but still solving the original PDEs \eqref{eqn.pde.nc}. As a result one obtains piecewise space-time polynomials of degree $M$, denoted by $\q_h(\x,t)$, that allow the scheme to achieve high order of accuracy even in time. An \textit{element-local} weak space-time formulation of the governing equations \eqref{eqn.pde.nc} is employed, following the approach developed in the Eulerian framework on fixed grids by Dumbser et al. in \cite{DumbserEnauxToro,USFORCE2,HidalgoDumbser}. According to \cite{DumbserEnauxToro,HidalgoDumbser,DumbserZanotti,Lagrange3D} we adopt the local space-time Discontinuous Galerkin predictor method due to the presence of stiff source terms in the governing equations \eqref{eqn.HPR}. Let $\mathbf{x}=(x,y,z)$ and $\boldsymbol{\xi}=(\xi,\eta,\zeta)$ be the spatial coordinate vectors defined in the physical and in the reference system, respectively, and let $\mathbf{\widetilde{x}}=(x,y,z,t)$ and $\boldsymbol{\widetilde{\xi}}=(\xi,\eta,\zeta,\tau)$ be the corresponding space-time coordinate vectors. Let furthermore $\theta_l=\theta_l(\boldsymbol{\widetilde{\xi}})=\theta_l(\xi,\eta,\zeta,\tau)$ be a space-time basis function defined by the Lagrange interpolation polynomials passing through the space-time nodes $\boldsymbol{\widetilde{\xi}}_m=(\xi_m,\eta_m,\zeta_m,\tau_m)$, which are defined by the tensor product of the spatial nodes of classical conforming high order finite elements in space and the Gauss-Legendre quadrature points in time.
Following \cite{DumbserPNPM}, the local solution $\q_h$, the fluxes $\F_h = (\f_{h}, \g_h, \h_h)$, the source term $\S_h$ and the non-conservative products $\P_h = \B(\q_h) \cdot \nabla \q_h$, are approximated within the space-time element $T_i(t) \times [t^n;t^{n+1}]$ with  
\begin{eqnarray}
\q_h=\q_h(\boldsymbol{\widetilde{\xi}}) = \theta_{l}(\boldsymbol{\widetilde{\xi}}) \, \widehat{\q}_{l,i}, \,
&& \F_h=\F_h(\boldsymbol{\widetilde{\xi}}) = \theta_{l}(\boldsymbol{\widetilde{\xi}}) \, \widehat{\F}_{l,i},  \nonumber \\
\S_h=\S_h(\boldsymbol{\widetilde{\xi}}) = \theta_{l}(\boldsymbol{\widetilde{\xi}}) \, \widehat{\S}_{l,i}, \,
&&\P_h=\P_h(\boldsymbol{\widetilde{\xi}}) = \theta_{l}(\boldsymbol{\widetilde{\xi}}) \, \widehat{\P}_{l,i}.
\label{thetaSol}
\end{eqnarray}
Since the Lagrange interpolation polynomials lead to a \textit{nodal} basis, we evaluate the degrees of freedom of $\F_h$, $\S_h$ and $\P_h$ from $\q_h$ in a pointwise
manner as 
\begin{equation}
  \widehat{\F}_{l,i} = \F(\widehat{\q}_{l,i}), \quad 
  \widehat{\S}_{l,i} = \S(\widehat{\q}_{l,i}), \quad 
  \widehat{\P}_{l,i} = \P(\widehat{\q}_{l,i},\nabla \widehat{\q}_{l,i}), \quad 
  \nabla \widehat{\q}_{l,i} = \nabla \theta_{m}(\boldsymbol{\widetilde{\xi}}_l) \widehat{\q}_{m,i},  
\end{equation}  
with $\nabla \widehat{\q}_{l,i}$ representing the gradient of $\q_h$ at node $\boldsymbol{\widetilde{\xi}}_l$.
An isoparametric approach is adopted, where the mapping between the physical space-time coordinate vector $\mathbf{\widetilde{x}}$ and the reference space-time coordinate vector $\boldsymbol{\widetilde{\xi}}$ is represented by the \textit{same} basis functions $\theta_l$ used for the discrete solution $\q_h$. Consequently we have
$ \x(\boldsymbol{\widetilde{\xi}}) = \theta_l(\boldsymbol{\widetilde{\xi}}) \, \widehat{\x}_{l,i}, $ and $ t(\boldsymbol{\widetilde{\xi}}) = \theta_l(\boldsymbol{\widetilde{\xi}}) \, \widehat{t}_l$,
where $\widehat{\mathbf{x}}_{l,i} = (\widehat{x}_{l,i},\widehat{y}_{l,i},\widehat{z}_{l,i})$ are the degrees of freedom of the spatial physical coordinates of the moving space-time control volume, which are unknown, while $\widehat{t}_l$ denote the \textit{known} degrees of freedom of the physical time at each space-time node $\widetilde{\x}_{l,i} = (\widehat{x}_{l,i}, \widehat{y}_{l,i}, \widehat{z}_{l,i}, \widehat{t}_l)$. The mapping in time is simply linear:
$t = t^n + \tau \, \Delta t$, then 
$\widehat{t}_l = t^n + \tau_l \, \Delta t$, 
with $t^n$ denoting the current time. $\Delta t$ is the time step and it is computed under a classical Courant-Friedrichs-Levy number (CFL) stability condition of the form
\begin{equation}
\Delta t = \textnormal{CFL} \, \min \limits_{T_i^n} \frac{d_i}{|\lambda_{\max,i}|}, \qquad \forall T_i^n \in \Omega^n, 
\label{eq:timestep}
\end{equation}
where $d_i$ is the insphere diameter of element $T_i^n$ and $|\lambda_{\max,i}|$ corresponds to the maximum absolute value of the eigenvalues computed from the solution $\Q_i^n$ in $T_i^n$. For the HPR model \eqref{eqn.HPR} the sound speed $c$ is computed according to \cite{PeshRom2014} as
\begin{equation}
  c = \sqrt{\frac{\gamma p}{\rho} + \frac{4}{3}c_s^2}.
\label{eqn.c}
\end{equation}
On unstructured meshes the CFL stability condition for explicit upwind schemes must satisfy the inequality $\textnormal{CFL} \leq \frac{1}{d}$.

We want the governing PDE formulation \eqref{eqn.pde.nc} to be written in the space-time reference system $\mathbf{\widetilde{x}}$, hence we first define the Jacobian of the space-time transformation from the physical to the reference element and its inverse:
\begin{equation}
J_{st} = \frac{\partial \mathbf{\widetilde{x}}}{\partial \boldsymbol{\widetilde{\xi}}} = \left( \begin{array}{cccc} x_{\xi} & x_{\eta} & x_{\zeta} & x_{\tau} \\ y_{\xi} & y_{\eta} & y_{\zeta} & y_{\tau} \\ z_{\xi} & z_{\eta} & z_{\zeta} & z_{\tau} \\ 0 & 0 & 0 & \Delta t \\ \end{array} \right) ,
 \quad \quad \quad
J_{st}^{-1} = \frac{\partial \boldsymbol{\widetilde{\xi}}}{\partial \mathbf{\widetilde{x}}} = \left( \begin{array}{cccc} \xi_{x} & \xi_{y} & \xi_{z} & \xi_{t} \\ \eta_{x} & \eta_{y} & \eta_{z} & \eta_{t} \\ \zeta_{x} & \zeta_{y} & \zeta_{z} & \zeta_{t} \\ 0 & 0 & 0 & \frac{1}{\Delta t} \\ \end{array} \right).
\label{JaciJac}
\end{equation}

Furthermore let us introduce the nabla operator $\nabla$ in the reference space $\boldsymbol{\xi}=(\xi,\eta,\zeta)$ and in the physical space $\mathbf{x}=(x,y,z)$ as:
\begin{equation}
 \nabla_{\xxi} = \left( \begin{array}{c} \frac{\partial}{\partial \xi} \\ \frac{\partial}{\partial \eta} \\ \frac{\partial}{\partial \zeta}  \end{array} \right), \qquad 
 \nabla        = \left( \begin{array}{c} \frac{\partial}{\partial x  } \\ \frac{\partial}{\partial y   } \\ \frac{\partial}{\partial z   }  \end{array} \right) = 
  \left( \begin{array}{ccc} \xi_x & \eta _x & \zeta _x \\ \xi_y & \eta_y & \zeta _y \\ \xi_z & \eta_z & \zeta _z \end{array} \right) 
  \left( \begin{array}{c} \frac{\partial}{\partial \xi} \\ \frac{\partial}{\partial \eta} \\ \frac{\partial}{\partial \zeta}  \end{array} \right)  = 
  \left( \frac{\partial \xxi}{\partial \x} \right)^T \nabla_{\xxi},
\label{not.ref.pde}
\end{equation}
and two integral operators
\begin{eqnarray}
\left[f,g\right]^{\tau} = \int \limits_{T_e} f(\xi,\eta,\zeta,\tau) g(\xi,\eta,\zeta,\tau)  \, d\boldsymbol{\xi}, \nonumber \quad \quad
\left\langle f,g \right\rangle = \int \limits_{0}^{1} \int \limits_{T_e} f(\xi,\eta,\zeta,\tau)g(\xi,\eta,\zeta,\tau)  \, d\boldsymbol{\xi} \, d\tau,  
\label{intOperators}
\end{eqnarray}
that denote the scalar products of two functions $f$ and $g$ over the spatial reference element $T_e$ at time $\tau$ and over the space-time reference element $T_e\times \left[0,1\right]$, respectively.

The system of balance laws \eqref{eqn.pde.nc} is then reformulated in the reference coordinate system $\mathbf{\widetilde{x}}$ with the following compact notation
\begin{equation}
\frac{\partial \Q}{\partial \tau} + \Delta t \; \mathbf{H} = \Delta t \; \mathbf{S}(\Q),
\label{PDECGsimple}
\end{equation}
where we have introduced the unified term $ \mathbf{H} = \frac{\partial \Q}{\partial \xxi} \cdot \frac{\partial \xxi}{\partial t} + \left( \frac{\partial \xxi}{\partial \x} \right)^T \nabla_{\xxi} \cdot \F  + \B(\Q) \cdot \left( \frac{\partial \xxi}{\partial \x} \right)^T \nabla_{\xxi} \Q $ by using the inverse of the associated Jacobian matrix \eqref{JaciJac} and the gradient notation \eqref{not.ref.pde}. The numerical approximation of $\mathbf{H}$ is computed by the same isoparametric approach \eqref{thetaSol}, i.e.
$\mathbf{H}_h = \theta_{l}(\boldsymbol{\widetilde{\xi}}) \, \widehat{\mathbf{H}}_{l,i}$.
Inserting this approximation and \eqref{thetaSol} into \eqref{PDECGsimple}, then multiplying \eqref{PDECGsimple} with a space-time test function $\theta_k(\xxi)$ and further integrating the resulting equation over the space-time reference element $T_e \times [0,1]$, one obtains a weak formulation of the original governing system \eqref{eqn.pde.nc}:
\begin{equation}
\left\langle \theta_k,\frac{\partial \theta_l}{\partial \tau} \right\rangle \widehat{\q}_{l,i}  
=  \left\langle \theta_k,\theta_l \right\rangle \Delta t  \; \left( \widehat{\mathbf{S}}_{l,i} - \widehat{\mathbf{H}}_{l,i} \right). \nonumber\\ 
\label{eqn.weak.lag} 
\end{equation}
The term on the left hand side can be integrated by parts in time considering the initial condition 
of the local Cauchy problem $\w^n_h$, yielding
\begin{equation}
 \left[ \theta_k(\xxi,1), \theta_l(\xxi,1)\right]^1 \widehat{\q}_{l,i} - \left\langle \frac{\partial \theta_k}{\partial \tau}, \theta_l \right\rangle \widehat{\q}_{l,i}  
= \left[ \theta_k(\xxi, 0), \psi_l(\xxi) \right]^0 \hat \w^n_{l,i} + \left\langle \theta_k,\theta_l \right\rangle \Delta t  \; \left( \widehat{\mathbf{S}}_{l,i}  - \widehat{\mathbf{H}}_{l,i}\right),
\label{LagrSTPDECG}
\end{equation}
that simplifies to 
\begin{equation}
 \mathbf{K}_1 \widehat{\q}_{l,i} = \F_0 \hat \w^n_{l,i} + \Delta t \; \mathbf{M}  \left( \widehat{\mathbf{S}}_{l,i} - \widehat{\mathbf{H}}_{l,i} \right),
 \label{DGsystem}
\end{equation}
with the following more compact matrix-vector notation:
\begin{equation}
\mathbf{K}_{1} = \left[ \theta_k(\xxi,1), \theta_l(\xxi,1)\right]^1 - \left\langle \frac{\partial \theta_k}{\partial \tau}, \theta_l \right\rangle, \quad 
\F_0   = \left[ \theta_k(\xxi, 0), \psi_l(\xxi) \right], \quad 
\mathbf{M} = \left\langle \theta_k,\theta_l \right\rangle.
\label{eqn.DGterms}
\end{equation} 
\textit{De facto} equation \eqref{DGsystem} constitutes an element-local nonlinear system of algebraic equations for the unknown space-time expansion coefficients $\widehat{\q}_{l,i}$ \footnote{
  This system is solved using the following iterative scheme 
\vspace{-0.2cm}
$$
\vspace{-0.2cm}
 \widehat{\q}_{l,i}^{r+1} - \Delta t  \; \mathbf{K}_1^{-1} \mathbf{M} \, \widehat{\mathbf{S}}^{r+1}_{l,i}  = \mathbf{K}_1^{-1} \left( \F_0 \hat \w^n_{l,i} - \Delta t  \; \mathbf{M} \widehat{\mathbf{H}}_{l,i}^r \right), 
 $$
where $r$ denotes the iteration number. Stiff algebraic source terms $\mathbf{S}$ are implicitly discretized, see \cite{DumbserEnauxToro,DumbserZanotti,HidalgoDumbser}.} .

Together with the solution, we have to evolve the geometry of the space-time control volume which moves in time. The motion of the nodes of element $T^n_i$ is described by the ODE system 
\begin{equation}
\frac{d \mathbf{x}}{dt} = \mathbf{V}(\Q,\x,t),
\label{ODEmesh}
\end{equation}
with $\mathbf{V}=\mathbf{V}(\Q,\x,t)$ denoting the local mesh velocity. Our direct Arbitrary-Lagrangian-Eulerian (ALE) method allows the mesh velocity to be chosen independently from the fluid velocity.
Following the same philosophy as for the solution, the velocity inside element $T_i(t)$ is also expressed in terms of the  space-time basis functions $\theta_{l}$ as $\; \mathbf{V}_h\;=\; \theta_{l}(\xxi,\tau) \widehat{\mathbf{V}}_{l,i} \;$, with the notation $\widehat{\mathbf{V}}_{l,i} = \mathbf{V}(\mathbf{\hat \q}_{l,i}, \hat{\x}_{l,i}, \hat t_l)$.
The local space-time DG method is used again to solve \eqref{ODEmesh} for the unknown coordinate vector $\widehat{\mathbf{x}}_l=(x_l,y_l,z_l)$, according to \cite{Lagrange2D,LagrangeNC}, hence 
\begin{equation}
\K_1 \widehat{\mathbf{x}}_{l,i} = \left[ \theta_k(\xxi,0), \x(\xxi,t^n) \right]^0 + \Delta t \; \mathbf{M} \, \widehat{\mathbf{V}}_{l,i},
\label{VCG}
\end{equation}
where $\x(\xxi,t^n)$ is given by the mapping based on the known vertex coordinates of simplex $T_i^n$ at time $t^n$. The above system is iteratively solved together with \eqref{DGsystem}.

Once the above procedure is performed for all cells, an element-local predictor for the numerical solution $\q_h$, for the fluxes $\mathbf{F}_h=(\f_h,\g_h,\h_h)$, for the non-conservative products $\P_h$, for the source term $\S_h$ and also for the mesh velocity $\mathbf{V}_h$ is available. This procedure is carried out \textit{locally} for each cell, consequently it remains to update the mesh motion \textit{globally}, by assigning a unique velocity vector to each node. To address this issue, in the next section a local nodal solver algorithm for the velocity together with an embedded rezoning technique are presented.

\subsection{Mesh motion}  \label{sec.meshMot}

The aim of any ALE scheme is to follow as closely as possible the material motion. 
This motion can generate highly deformed cells specifically when fluids or gases are considered. That may drastically reduce the admissible timestep, or, worse, may lead to tangled elements. 
In order to guarantee good resolution properties for contact waves and material interfaces together with a good geometrical mesh quality, the mesh velocity must be chosen carefully. 
When natural evidences emanate from the motion of the material boundary conditions, such a 
mesh velocity can be inferred. However in the general case, specifically for fluids and gases, we
adopt a  suitable Lagrangian nodal solver technique \cite{Despres2005,Maire2009,chengshu1,chengshu2} 
to assign a unique velocity vector to each node accurately representing the \quotew{true} material velocity. 
Notice that since we are dealing with a direct ALE formulation the mesh velocity is a degree of freedom. 
As a consequence we could run our ALE code in a pure Eulerian regime by setting 
the mesh velocity to zero, or in an almost Lagrangian regime by setting the velocity to an
local average of the computed Lagrangian velocities. We could also force any sort of intermediate
or artificial mesh motion leading \textit{de facto} to a so-called ALE motion. 
In this work the simple nodal solver of Cheng and Shu is used \cite{chengshu1,chengshu2} and the rezoning strategy exposed in \cite{MaireRezoning,Lagrange3D} is employed to locally improve the mesh quality. The final mesh configuration, i.e. the vertex coordinates at the new time level $t^{n+1}$ are then computed relying on the relaxation algorithm presented in \cite{MaireRezoning}.
\subsection{Finite volume scheme}
\label{sec.SolAlg}

The same approach already developed in two and three space dimensions discussed in \cite{Lagrange2D,LagrangeNC,Lagrange3D} is briefly summarized here. To begin with, the governing PDE \eqref{eqn.pde.nc} is more compactly reformulated using a space-time divergence operator $\widetilde \nabla$:
\begin{equation}
\widetilde \nabla \cdot \widetilde{\F} + \widetilde \B(\Q) \cdot \widetilde \nabla \Q = \mathbf{S}(\Q),  \qquad \widetilde \nabla  = \left( \frac{\partial}{\partial x}, \, \frac{\partial}{\partial y}, \, \frac{\partial}{\partial z}, \, \frac{\partial}{\partial t} \right)^T,
\label{eqn.st.pde}
\end{equation}
where the space-time flux tensor $\widetilde{\F}$ and the system matrix $\widetilde{\B}$ are given by $\widetilde{\F}  = \left( \mathbf{f}, \, \mathbf{g}, \, \mathbf{h}, \, \Q \right)$ and $\widetilde{\B}  = ( \B_1, \B_2, \B_3, 0)$.  
For the computation of the state vector at the new time level $\Q^{n+1}$, the balance law \eqref{eqn.st.pde} is integrated over a four-dimensional space-time control volume $\mathcal{C}^n_i = T_i(t) \times \left[t^{n}; t^{n+1}\right]$, which after the application of the theorem of Gauss yields
\begin{equation}
\int \limits_{\partial \mathcal{C}^{n}_i} \widetilde{\F} \cdot \ \mathbf{\widetilde n} \, \, dS +  
\int\limits_{\mathcal{C}^n_i} \widetilde{\B}(\Q) \cdot \widetilde \nabla \Q \, d\mathbf{x} \, dt = 
\int\limits_{\mathcal{C}^n_i} \S(\Q) \, d\mathbf{x} dt.
\label{STPDEgauss}
\end{equation} 

The non-conservative products are treated with the path-conservative approach  of Castro and Par\'es, see \cite{Toumi1992,Pares2006,Castro2006,Castro2008,Rhebergen2008,ADERNC,USFORCE2,OsherNC,LagrangeNC}, for a non-exhaustive overview, hence leading to 
\begin{equation}
\int \limits_{\partial \mathcal{C}^{n}_i} \left( \widetilde{\F} + \widetilde{\D} \right) \cdot \ \mathbf{\widetilde n} \, \,  dS + \! \! 
\int \limits_{\mathcal{C}^n_i \backslash \partial \mathcal{C}^n_i} \! \! \! \widetilde{\B}(\Q) \cdot \widetilde \nabla \Q \, \,  d\mathbf{x} \, dt = 
\int\limits_{\mathcal{C}^n_i} \S(\Q) \, d\mathbf{x} dt,   
\label{I1}
\end{equation}
where a new term $\widetilde{\D}$ has been introduced in order to take into account the jumps of the solution $\Q$ on the space-time element boundaries $\partial \mathcal{C}^n_i$. This term is computed by the path integral
\begin{equation}
 \widetilde{\D} \cdot  \mathbf{\widetilde n}  = \halb \int \limits_0^1 \widetilde{\B}\left(\Path (\Q^-,\Q^+,s)\right) \cdot \mathbf{\widetilde n} \, \frac{\partial \Path}{\partial s} \, ds = \halb  \left( \int \limits_0^1 \widetilde{\B}\left(\Path(\Q^-,\Q^+,s)\right) \cdot \mathbf{\widetilde n} \, ds \right) \left( \Q^+ - \Q^- \right), 
 \label{eqn.pathint} 
\end{equation}
where the integration path $\Path$ in \eqref{eqn.pathint} is chosen according to \cite{Pares2006,Castro2006,USFORCE2,OsherNC} to be a simple straight-line segment, i.e. $\Path(\Q^-,\Q^+,s) = \Q^- + s (\Q^+ - \Q^-)$, and $\left(\Q^-,\Q^+\right)$ are the conserved variables in element $T_i^n$ and its direct neighbor $T_j^n$, respectively. 
Moreover $\mathbf{\widetilde n} = (\widetilde n_x,\widetilde n_y,\widetilde n_z,\widetilde n_t)$ denotes the outward pointing space-time unit normal vector on the varying space-time volume $\partial C^n_i$.

Let $\mathcal{N}_i$ denote the \textit{Neumann neighborhood} of simplex $T_i(t)$, which is the set of directly adjacent neighbors $T_j(t)$ that share a common face $\partial T_{ij}(t)$ with $T_i(t)$. The space-time volume $\partial C^n_i$ is composed of $d+1$ space-time sub-volumes $\partial C^n_{ij}$, each of them defined for each face of $T_i(t)$, and two more space-time sub-volumes, $T_i^{n}$ and $T_i^{n+1}$, that represent the simplex configuration at times $t^n$ and $t^{n+1}$, respectively (see \cite{Lagrange3D} for details). Therefore the space-time volume $\partial C^n_i$ involves overall a total number of $2+d+1$ space-time sub-volumes, i.e.
\begin{equation}
\partial C^n_i = \left( \bigcup \limits_{T_j(t) \in \mathcal{N}_i} \partial C^n_{ij} \right) 
\,\, \cup \,\, T_i^{n} \,\, \cup \,\, T_i^{n+1}.  
\label{dCi}
\end{equation} 
Each of the space-time sub-volumes is mapped to a reference element in order to simplify the integral computation. For the configurations at the current and at the new time level, $T_i^{n}$ and $T_i^{n+1}$, we use the mapping from the physical to the reference element. The space-time unit normal vectors simply read $\mathbf{\widetilde n} = (0,0,0,-1)$ for $T_i^{n}$ and $\mathbf{\widetilde n} = (0,0,0,1)$ for $T_i^{n+1}$, since these volumes are orthogonal to the time coordinate. For the lateral sub-volumes $\partial C^n_{ij}$ we adopt a linear parametrization to map the physical volume to a $d+1$-dimensional space-time reference prism \cite{Lagrange3D}. 

Starting from the old vertex coordinates $\mathbf{X}_{ik}^n$ and the new ones $\mathbf{X}_{ik}^{n+1}$, that are \textit{known} from the mesh motion algorithm described in Section \ref{sec.meshMot}, the lateral sub-volumes are parametrized using a set of linear basis functions $\beta_k(\chi_1,\chi_2,\tau)$ that are defined on a local reference system $\boldsymbol{\chi}=(\chi_1,\chi_2,\tau)$ which is oriented orthogonally w.r.t. the face $\partial T_{ij}(t)$ of $T_i^n$, e.g. the reference time coordinate $\tau$ is orthogonal to the reference space coordinates $(\chi_1,\chi_2)$ that lie on $\partial T_{ij}(t)$. The temporal mapping is simply given by $t = t^n + \tau \, \Delta t$, hence $t_{\chi_1} = t_{\chi_2} = 0$ and $t_\tau = \Delta t$. The lateral space-time volume $\partial C_{ij}^n$ is defined by six vertices of physical coordinates $\mathbf{\widetilde{X}}_{ij,k}^n$. The first three vectors $(\mathbf{X}^n_{ij,1},\mathbf{X}^n_{ij,2},\mathbf{X}^n_{ij,3})$ are the nodes defining the common face $\partial T_{ij}(t^n)$ at time $t^n$, while the same procedure applies at the new time level $t^{n+1}$. Therefore the six vectors $\mathbf{\widetilde{X}}_{ij,k}^n$ are given by
\begin{eqnarray}
\mathbf{\widetilde{X}}_{ij,k}^n     = \left( \mathbf{X}^n_{ij,k}, t^n \right), \quad \text{ and } \quad
\mathbf{\widetilde{X}}_{ij,k+3}^n = \left( \mathbf{X}^{n+1}_{ij,k}, t^{n+1} \right), \quad k=1,2,3 .
\label{eqn.lateralnodes} 
\end{eqnarray}
The parametrization for $\partial C_{ij}^n$ reads
\begin{equation}
\partial C_{ij}^n = \mathbf{\widetilde{x}} \left( \chi_1,\chi_2,\tau \right) = 
 \sum\limits_{k=1}^{6}{\beta_k(\chi_1,\chi_2,\tau) \, \mathbf{\widetilde{X}}_{ij,k}^n } ,
\label{eqn.DOFbeta}
\end{equation}
with $0 \leq \chi_1 \leq 1$, $0 \leq \chi_2 \leq 1-\chi_1$ and $0 \leq \tau \leq 1$ and the linear basis functions $\beta_k(\chi_1,\chi_2,\tau)$ given by 
\begin{eqnarray}
\beta_1(\chi_1,\chi_2,\tau) = (1-\chi_1-\chi_2)(1-\tau),&& \quad \beta_4(\chi_1,\chi_2,\tau) = (1-\chi_1-\chi_2)(\tau) \nonumber \\
\beta_2(\chi_1,\chi_2,\tau) = \chi_1(1-\tau),&    &         \quad \beta_5(\chi_1,\chi_2,\tau) = \chi_1\tau, \nonumber \\
\beta_3(\chi_1,\chi_2,\tau) = \chi_2(1-\tau),&     &        \quad \beta_6(\chi_1,\chi_2,\tau) = \chi_2\tau.
\label{eq:BetaBaseFunc}
\end{eqnarray}
The coordinate transformation is associated with a matrix $\mathcal{T}$ that reads
\begin{equation}
\mathcal{T} = \left( \hat{\mathbf{e}}, \frac{\partial \mathbf{\widetilde{x}}}{\partial \chi_1}, \frac{\partial \mathbf{\widetilde{x}}}{\partial \chi_2}, \frac{\partial \mathbf{\widetilde{x}}}{\partial \tau} \right)^T,
\label{eqn.JacSTreference}
\end{equation}
with $\hat{\mathbf{e}}=(\hat{\mathbf{e}}_1,\hat{\mathbf{e}}_2,\hat{\mathbf{e}}_3,\hat{\mathbf{e}}_4)$. Let $\hat{\mathbf{e}}_p$ represent the unit vector aligned with the $p$-th axis of the physical coordinate system $(x,y,z,t)$ and let $\widetilde{x}_q$ denote the $q$-th component of vector $\mathbf{\widetilde{x}}$. The determinant of $\mathcal{T}$ produces at the same time the quantity $| \partial C_{ij}^n|$ of the space-time sub-volume $\partial C_{ij}^n$ and the space-time normal vector $\mathbf{\widetilde n}_{ij}$, as     
\begin{equation}
 \mathbf{\widetilde n}_{ij} = \left( \epsilon_{pqrs} \, \hat{\mathbf{e}}_p \, \frac{\partial {\widetilde{x}_q}}{\partial \chi_1} \, \frac{\partial {\widetilde{x}_r}}{\partial \chi_2} \, \frac{\partial {\widetilde{x}_s}}{\partial \tau} \right)/ | \partial C_{ij}^n|,
 \label{n_lateral1}
\end{equation}
where the \textit{Levi-Civita} symbol has been used according to the usual definition
\begin{equation}
\epsilon_{pqrs} = \left\{ \begin{array}{l} +1, \quad \textnormal{if $(p,q,r,s)$ is an \textit{even} permutation of $(1,2,3,4)$}, \\
																					 -1, \quad \textnormal{if $(p,q,r,s)$ is an \textit{odd} permutation of $(1,2,3,4)$}, \\
																					 0,  \quad \textnormal{otherwise,}
                      \end{array}  \right.  
\label{eqn.LeviCivita}
\end{equation}
and with 
\begin{equation*} 
| \partial C_{ij}^n| = \left\| \epsilon_{pqrs} \, \hat{\mathbf{e}}_p \, \frac{\partial {\widetilde{x}_q}}{\partial \chi_1} \, \frac{\partial {\widetilde{x}_r}}{\partial \chi_2} \, \frac{\partial {\widetilde{x}_s}}{\partial \tau} \right\|.   
\end{equation*} 

The final one-step direct ALE ADER-WENO finite volume scheme takes the following form: 
\begin{equation}
|T_i^{n+1}| \, \Q_i^{n+1} = |T_i^n| \, \Q_i^n - \sum \limits_{T_j \in \mathcal{N}_i} \,\, {\int \limits_0^1 \int \limits_0^1 \int \limits_{0}^{1-\chi_1} 
| \partial C_{ij}^n| \widetilde{\G}_{ij} \cdot \mathbf{\widetilde n}_{ij} \, d\chi_2 \, d\chi_1 \, d\tau}
+ \int \limits_{\mathcal{C}_i^n \backslash \partial \mathcal{C}_i^n}  \left( \S_h - \P_h \right) \, d\mathbf{x} \, dt, 
\label{PDEfinal}
\end{equation}
where in the term $\widetilde{\G}_{ij} \cdot \mathbf{\widetilde n}_{ij}$ the Arbitrary-Lagrangian-Eulerian numerical flux function is embedded, as well as the path-conservative jump term, 
which allows the discontinuity of the predictor solution $\mathbf{q}_h$ that occurs at the space-time boundary $\partial C_{ij}^n$ to be properly resolved also in the presence of 
non-conservative products. The volume integrals in \eqref{PDEfinal} are approximated using multidimensional Gaussian quadrature rules \cite{stroud} of suitable order of accuracy and 
the term $\widetilde{\G}_{ij}$ is evaluated relying on a simple ALE Rusanov-type scheme \cite{Lagrange1D,Lagrange2D,Lagrange3D} as 
\begin{equation}
  \widetilde{\G}_{ij} =  
  \frac{1}{2} \left( \widetilde{\F}(\q_h^+) + \widetilde{\F}(\q_h^-)  \right) \cdot \mathbf{\widetilde n}_{ij} +  
  \frac{1}{2} \left( \int \limits_0^1 \widetilde{\B}(\Path)\cdot \mathbf{\widetilde n} \ ds - |\lambda_{\max}| \mathbf{I} \right) \left( \q_h^+ - \q_h^- \right), 
  \label{eqn.rusanov} 
\end{equation}
where $\q_h^-$ and $\q_h^+$ are the local space-time predictor solution inside element $T_i(t)$ and the neighbor $T_j(t)$, respectively, and $|\lambda_{\max}|$ denotes the maximum absolute value of the eigenvalues of the matrix $\widetilde{\A} \cdot \mathbf{\widetilde n}$ in space-time normal direction. Using the normal mesh velocity $\mathbf{V} \cdot \mathbf{n}$, matrix $\widetilde{\A}$ reads   
\begin{equation} 
\widetilde{\A} \cdot \mathbf{\widetilde n} = \left( \sqrt{\widetilde n_x^2 + \widetilde n_y^2 + \widetilde n_z^2} \, \right) \left[ \left( \frac{\partial \mathbf{F}}{\partial \Q} + \B \right) \cdot \mathbf{n} - (\mathbf{V} \cdot \mathbf{n}) \,  \mathbf{I} \right],  
\end{equation}
with $\mathbf{I}$ denoting the $\nu \times \nu$ identity matrix, $\A = \partial \F / \partial \Q + \B$ representing the classical Eulerian system matrix and $\mathbf{n}$ being the spatial unit normal vector given by $
\mathbf{n} = \frac{(\widetilde n_x, \widetilde n_y, \widetilde n_z)^T}{\sqrt{\widetilde n_x^2 + \widetilde n_y^2 + \widetilde n_z^2}}$.


Finally we remark that the integration over a closed space-time control volume, as done in this scheme, automatically respects the geometric conservation law (GCL), since application of Gauss' theorem yields
$
 \int_{\partial \mathcal{C}_i^n} \mathbf{\widetilde n} \, dS = 0.
$
As already pointed out in \cite{Lagrange3D,ALE-MOOD2} the numerical method allows a mass flux even for \quotew{Lagrangian} motion. Consequently there is no associated pure Lagrangian scheme \textit{in sensu stricto} to this numerical method. Nonetheless, very accurate results can still be achieved with this high order accurate ALE scheme , see 
\cite{Lagrange2D,LagrangeNC,LagrangeMHD,LagrangeQF,Lagrange3D}.

%
%
\subsection{Timestep constraint}
\label{sec.timestep}
The timestep $\Delta t$, which is needed for the discretization of the governing equations \eqref{PDEfinal}, is computed taking into account two different criteria, namely a classical CFL stability condition and a user-defined geometrical limitation.
The Courant-Friedrichs-Levy (CFL) stability condition is given by \eqref{eq:timestep}, while the second criterion is based on the limitation of the rate of change of the element volume within one timestep, i.e. the volume of each cell $T_i^n$ is not allowed either to increase more than a certain threshold which is provided by the user at the beginning of the computation, see \cite{LagrangeNC,Lagrange3D,ALE-MOOD2} for details.


%
%
\section{Boundary conditions}  \label{sec.BCs}
In this section we design appropriate boundary conditions for the HPR model and the ALE ADER-WENO finite volume schemes employed in this work. From the practical viewpoint of implementation, the boundary conditions setting assigns a suitable \textit{boundary state} $\Q_g$ for the ghost neighbor $T_g$ of element $T_i$, which lies on boundary of the domain given its state $\Q_i$. The set of boundary conditions needed to run the test cases reported in Section \ref{sec:numerics} are the following ones:

\begin{itemize}
	\item \textit{Transmissive} boundary conditions are adopted to let the fluid flow across the domain boundary. The flow is governed by the internal state, hence yielding the simple setting $\Q_g = \Q_i$;
	\item \textit{Wall} (or reflective) boundary conditions are used for the treatment of wall boundaries. In this case the normal flux across the domain boundary is zero, therefore we first set $\Q_g = \Q_i$ and then the velocity vector $\v_g=(u_g,v_g)$ for the boundary state is computed as
	  \begin{equation}
		\v_g = \v_i - 2 \, \left( \v_i \cdot \mathbf{n} \right) \, \mathbf{n}, 
		\label{eqn.WallBC}
		\end{equation}
		where $\mathbf{n}$ denotes as usual the outward pointing unit normal vector on the boundary edge of element $T_i$ and $\v_i$ represents the velocity vector of the internal state $\Q_i$. This treatment is also called \textit{no-slip wall} boundary condition and, for \textit{inviscid} flows, the fluid is still allowed to flow along the boundary, i.e. tangential to the boundary edge;
	\item \textit{Free traction} boundary conditions are normally employed in the context of solid mechanics, where the viscous stress tensor components are set to zero in order to discard the stresses at boundaries. In the HPR model \eqref{eqn.HPR} we only have a control on the distortion tensor $\AAA$ and it is not possible to derive an analytical function of the type $\AAA=f(\BS)$. Therefore the following simple strategy has been designed: we compute the value of the distortion tensor at the boundary $\AAA^{*}$ via a \textit{stiff relaxation} to the stress-free boundary state, using the same source terms as in the original governing PDE system, but with a different relaxation time $\tau_1' \to 0$. Hence, for the free traction boundaries we solve the ODE 
		\begin{equation}
		\frac{d \AAA}{dt} = -\frac{\boldsymbol{\psi}(\AAA)}{\theta_1(\tau_1')},
		\label{eqn.A_g}
		\end{equation}
		with a simple \textit{implicit} backward Euler scheme, which yields the following nonlinear algebraic equation for the unknown tensor $\AAA^{*}$   
		\begin{equation}
		\AAA^{*} + \frac{\Delta t}{\theta_1(\tau_1')} \boldsymbol{\psi}(\AAA^{*}) = \AAA_i, 
		\label{eqn.A_g2} 
		\end{equation}
		that can be easily solved with a standard Newton method. Here, $\AAA_i$ is the known distortion tensor of the boundary element $T_i$ and $\AAA^{*}$ is the  distortion tensor on the boundary  	
		edge of element $T_i$. Note again that the source term on the right hand side of the ODE \eqref{eqn.A_g} is the same as the strain relaxation term given by the governing PDE \eqref{eqn.HPR}, 
		but with smaller relaxation time. The solution of Eqn. \eqref{eqn.A_g2} provides the sought boundary ghost distortion tensor $\AAA_g$ as 
		\begin{equation}
		\AAA_g = 2 \, \AAA^{*} - \AAA_i. 
		\label{eqn.A_g-final}
		\end{equation}
		Moreover, we also require the hydrodynamic part of the pressure to vanish at the free surface boundary, hence setting $p_g=-p_i$, while the remaining variables are copied from the 
		internal state $\Q_i$; 
	\item \textit{moving} boundary conditions impose a prescribed velocity vector $\v_b$ on the boundary, hence they are classically treated 
by imposing
	 \begin{equation}
		\v_g = 2 \, \left( \v_b \cdot \mathbf{n} \right) \, \mathbf{n} - \v_i,
		\label{eqn.MovingBC}
		\end{equation}
		after setting $\Q_g = \Q_i$ for the remaining variables.
\end{itemize}
We underline that for finite volume schemes no \quotew{canonical} procedure is available to specify the boundary conditions. Thus, different ways are possible and, in principle, equally appropriate.

%
%
\section{Numerical experiments} \label{sec:numerics}

The aim of this section is to describe and show the numerical results for a list of representative test cases for the HPR model \eqref{eqn.HPR}. The numerical solution is provided by the direct ALE ADER-WENO finite volume schemes presented in this paper, employing piecewise polynomial reconstructions of degree $M=1,2,3$. The CFL number is generally set to $0.5$, if not stated otherwise, 
and all tests are run on unstructured meshes made of $N_E$ triangular elements. The computational grids are automatically generated by an external software and the characteristic mesh size 
is denoted in the following by $h$. 

Since the HPR model \eqref{eqn.HPR} is able to handle in \textit{one single PDE system} both \textbf{fluid mechanics} \textit{and} \textbf{solid mechanics}, our methodology of validation and verification involves those two branches of continuum mechanics. We clearly state that physical units are based on the $[m,kg,s]$ unit system for fluid mechanics, while we rely on the $[cm,g,\mu s]$ system for solid mechanics with the stresses measured in $[Mbar]$. The ideal gas \eqref{eq:ideal_gas_eos} equation of state is employed for fluids, whereas the Mie-Gr{\"u}neisen EOS \eqref{eq:mie_gruneisen_eos} is used for solids as usually done \cite{Kluth2010,Burton2013,Maire_elasto_13}.    

Unless explicitly given, for each test case simulated in the following the thermal impulse vector is set to zero, i.e. $\JJ=\mathbf{0}$ with $\tau_2 \to \infty$, and the mesh velocity is chosen to be equal to the local fluid velocity computed with the nodal solver of Cheng and Shu \cite{LagrangeMHD}, hence achieving a Lagrangian-like behavior of our direct ALE scheme.

According to \cite{PeshRom2014,Dumbser_HPR_16}, in the case of fluid mechanics the relation between the relaxation time $\tau_1$ and the dynamic viscosity coefficient $\mu$ is given by 
\begin{equation}
\mu=\frac{1}{6}\tau_1 \rho_0 c_s^2,
\label{eqn.mu}
\end{equation}
which allows us to set either the relaxation time $\tau_1$ or the viscosity coefficient $\mu$ as parameter of the HPR model. For \textit{inviscid} fluids we simply set $\tau_1 \to 0$ as discussed in \cite{PeshRom2014,Dumbser_HPR_16}.

Regarding solid mechanics, if $\tau_1 \to \infty$, we can describe the governing equations of \textit{pure elastic solids}, while for general \textit{elasto-plastic solids} we compute the relaxation time $\tau_1$ following \cite{BartonRomenski2010} as a power law function, i.e.
\begin{equation}
\tau_1 = \tau_0 \left( \frac{\sigma_0}{\sigma_I}\right)^n
\label{eqn.tauS}
\end{equation} 
where $\tau_0$, $\sigma_0$ and $n$ are material specific constants and the shear stress intensity $\sigma_I$ is evaluated by  
\begin{equation}
\sigma_I = \sqrt{\frac{1}{2} \left[ (\sigma_{11}-\sigma_{22})^2 + (\sigma_{22}-\sigma_{33})^2 + (\sigma_{33}-\sigma_{11})^2 + 6 \, (\sigma_{12}^2+\sigma_{13}^2+\sigma_{23}^2) \right]}.
\label{eqn.Y}
\end{equation}
Note that the parameter $\sigma_0$ corresponds to the so-called Yield stress of the material under quasistatic loading and the generic quantity $\sigma_{ik}$ is a component of the viscous shear stress  tensor $\boldsymbol{\sigma}$ given by \eqref{eqn.stress}. In Table \ref{tab:dataEOS} we report some mechanical constants as well as the parameters needed in the Mie-Gr{\"u}neisen EOS for the materials considered in the test cases for solid mechanics presented in this paper.  

\begin{table}
\caption{Material parameters: reference density $\rho_0$, reference (atmospheric) pressure $\p_0$, adiabatic sound speed $c_0$, shear wave speed $c_s$, Yield stress $\sigma_0$ and the coefficients $\Gamma_0$ and $s$ appearing in the Mie-Gr{\"u}neisen equation of state \eqref{eq:mie_gruneisen_eos}.}
\centering
\begin{tabular}{|c||ccccccc|}
\hline
						& $\rho_0$	& $p_0$ 	& $c_0$ 	& $c_s$ 	& $\sigma_0$ 	& $\Gamma_0$ 	& s \\ 
\hline
Copper      & $8.930$ 	& $0.0$ 	& $0.394$ & $0.219$ & $0.004$  		& $2.00$ 	& $1.480$ \\
Beryllium  	& $1.845$ 	& $0.0$ 	& $1.287$ & $0.905$ & $1$					& $1.11$ 	& $1.124$ \\
Aluminum 		& $2.785$ 	& $0.0$ 	& $0.533$ & $0.305$ & $0.003$			& $2.00$ 	& $1.338$ \\
\hline
\end{tabular}
\label{tab:dataEOS}
\end{table}

\subsection{Numerical convergence results} \label{ssec:num_conv}

As fully detailed in \cite{Dumbser_HPR_16}, a zeroth order approximation of the HPR model can be obtained in the stiff limit $\tau_1 \to 0$ because the viscous stresses vanish, therefore retrieving the compressible Euler equations which govern an inviscid fluid. In this way, we can use the smooth isentropic vortex test problem presented in \cite{HuShuTri} to study the numerical convergence of our finite volume schemes. The initial computational domain is the square $\Omega(0)=[0;10]\times[0;10]$ with periodic boundaries everywhere. The initial condition is given in terms of primitive variables and it reads
\begin{equation}
(\rho, u, v, p) = (1+\delta \rho, 1+\delta u, 1+\delta v, 1+\delta p),
\label{eqn.iniSV}
\end{equation}
where the symbol $\delta$ represents the perturbations superimposed to a homogeneous background field. Since the vortex is isentropic, the entropy perturbation is assumed to be zero, i.e. $S=\frac{p}{\rho^\gamma}=0$, and the perturbations for density and pressure are 
\begin{equation}
\delta \rho = (1+\delta T)^{\frac{1}{\gamma-1}}-1, \quad \delta p = (1+\delta T)^{\frac{\gamma}{\gamma-1}}-1 \qquad \textnormal{ with } \quad \delta T = -\frac{(\gamma-1)\epsilon^2}{8\gamma\pi^2}e^{1-r^2}.
\label{eqn.drho-dp}
\end{equation} 
The generic radial coordinate is $r=\sqrt{(x-5)^2+(y-5)^2}$ and the velocity perturbation is given by
\begin{equation}
\left(\begin{array}{c} \delta u \\ \delta v \end{array}\right) = \frac{\epsilon}{2\pi}e^{\frac{1-r^2}{2}} \left(\begin{array}{c} -(y-5) \\ \phantom{-}(x-5) \end{array}\right),
\label{eqn.du}
\end{equation}
with $\epsilon=5$ denoting the vortex strength. The initial distortion tensor is set to $\AAA=\sqrt[3]{\rho} \, \mathbf{I}$ and the final time of the simulation is taken to be $t_f=1.0$. The parameters for the HPR model are $\gamma=1.4$, $c_v=2.5$, $\rho_0=1$, $c_s=0.5$ and the relaxation time is $\tau_1=10^{-12}$, which corresponds to the stiff inviscid limit $\tau_1 \to 0$. We run our direct ADER-WENO+ALE finite volume schemes on a series of successively refined grids up to fourth order of accuracy in space and time. The reference solution $\Q_e$ is given by the exact solution of the compressible Euler equations and it can be simply computed as the time-shifted initial condition, e.g. $\Q_e(\x,t_f)=\Q(\x-\v_c t_f,0)$, with the convective mean velocity $\v_c=(1,1)$. The error is measured at time $t_f$ using the continuous $L_2$ norm and the resulting convergence rates are listed in Table \ref{tab.convHPR}, confirming clearly that the proposed numerical method is able to achieve its designed order of
accuracy for smooth problems in the stiff relaxation limit $\tau_1 \to 0$.

\begin{table}  
\caption{Numerical convergence results for ALE ADER-WENO finite volume schemes applied to the HPR model in the stiff inviscid limit. The error norms refer to the variable $\rho$ (density) at time $t_f=1.0$ for first up to fourth order of accuracy and the exact solution is given by the inviscid compressible Euler equations.} 
\begin{center} 
\renewcommand{\arraystretch}{1.0}
\begin{tabular}{|ccc||ccc|} 
\hline
  $h(\Omega(t_f))$ & $\epsilon_{L_2}$ & $\mathcal{O}(L_2)$ & $h(\Omega,t_f)$ & $\epsilon_{L_2}$ & $\mathcal{O}(L_2)$ \\
\hline
  \multicolumn{3}{|c||}{
  \textbf{1st order ADER-WENO-ALE}
  } & \multicolumn{3}{c|}{
  \textbf{2nd order ADER-WENO-ALE}} \\
\hline
3.40E-01 & 3.084E-01 & -            & 3.70E-01 & 7.880E-02 & -       \\ 
2.48E-01 & 2.556E-01 & 0.6          & 2.48E-01 & 5.907E-02 & 0.7    \\ 
1.71E-01 & 1.921E-01 & 0.8          & 1.73E-01 & 2.542E-02 & 2.3    \\ 
1.33E-01 & 1.533E-01 & 0.9          & 1.28E-01 & 1.443E-02 & 1.9    \\ 
\hline 
  \multicolumn{3}{|c||}{
  \textbf{3rd order ADER-WENO-ALE}} & \multicolumn{3}{c|}{
  \textbf{4th order ADER-WENO-ALE}}   \\
\hline
3.37E-01  & 4.861E-02 & -           & 3.28E-01 & 1.746E-02 & -         \\ 
2.51E-01  & 2.806E-02 & 1.9         & 2.51E-01 & 6.416E-03 & 3.8   \\ 
1.68E-01  & 1.090E-02 & 2.3         & 1.68E-01 & 1.238E-04 & 4.1   \\ 
1.28E-01  & 5.052E-03 & 2.8         & 1.28E-01 & 3.728E-04 & 4.4   \\ 
\hline 
\end{tabular}
\end{center}
\label{tab.convHPR}
\end{table}

\subsection{2D Taylor-Green vortex} \label{sssec:vortex}
\begin{figure}[!htbp]
\begin{center}
\begin{tabular}{cc} 
\includegraphics[width=0.47\textwidth]{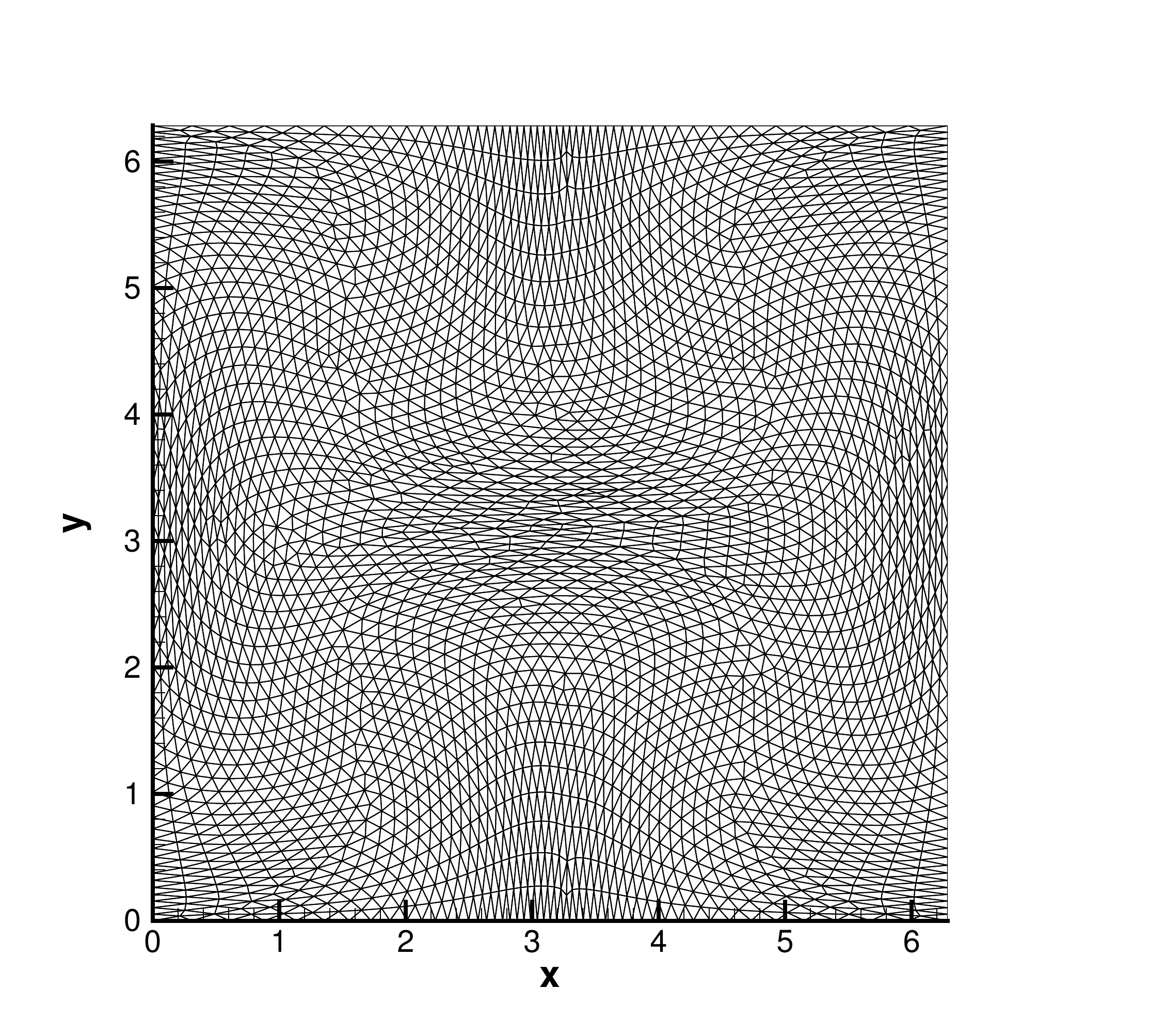} & 
\includegraphics[width=0.47\textwidth]{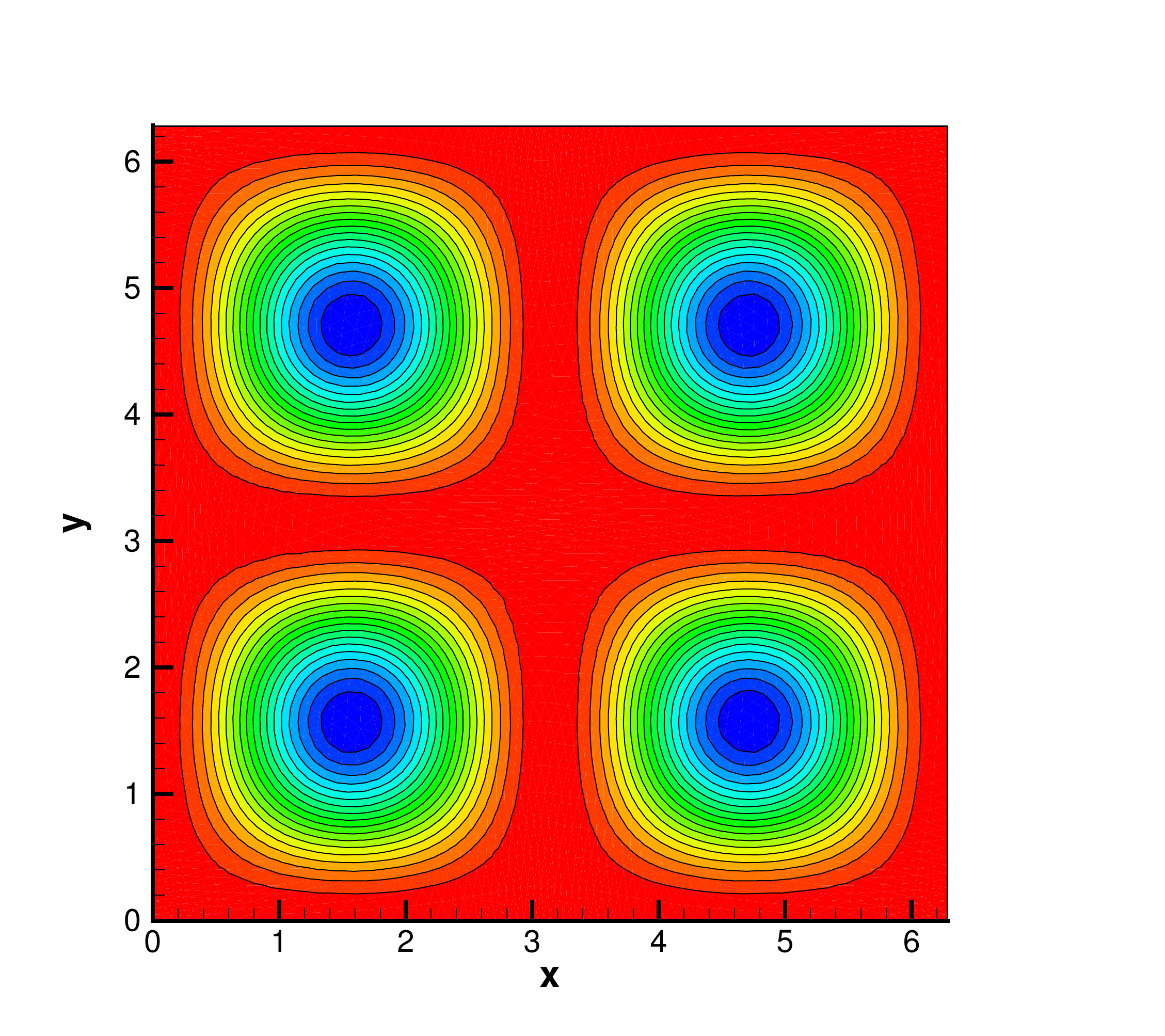}  \\  
\includegraphics[width=0.47\textwidth]{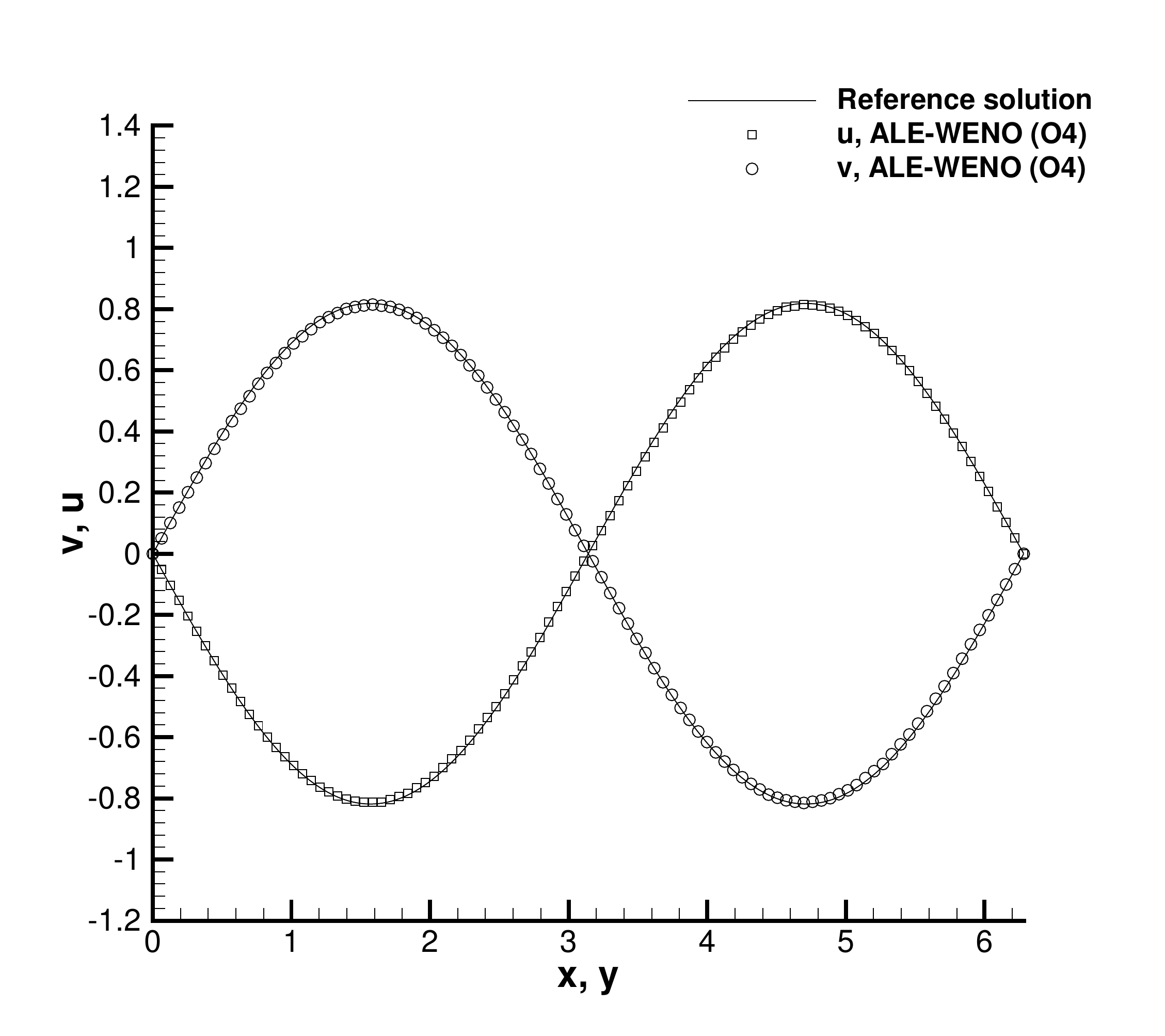}   & 
\includegraphics[width=0.47\textwidth]{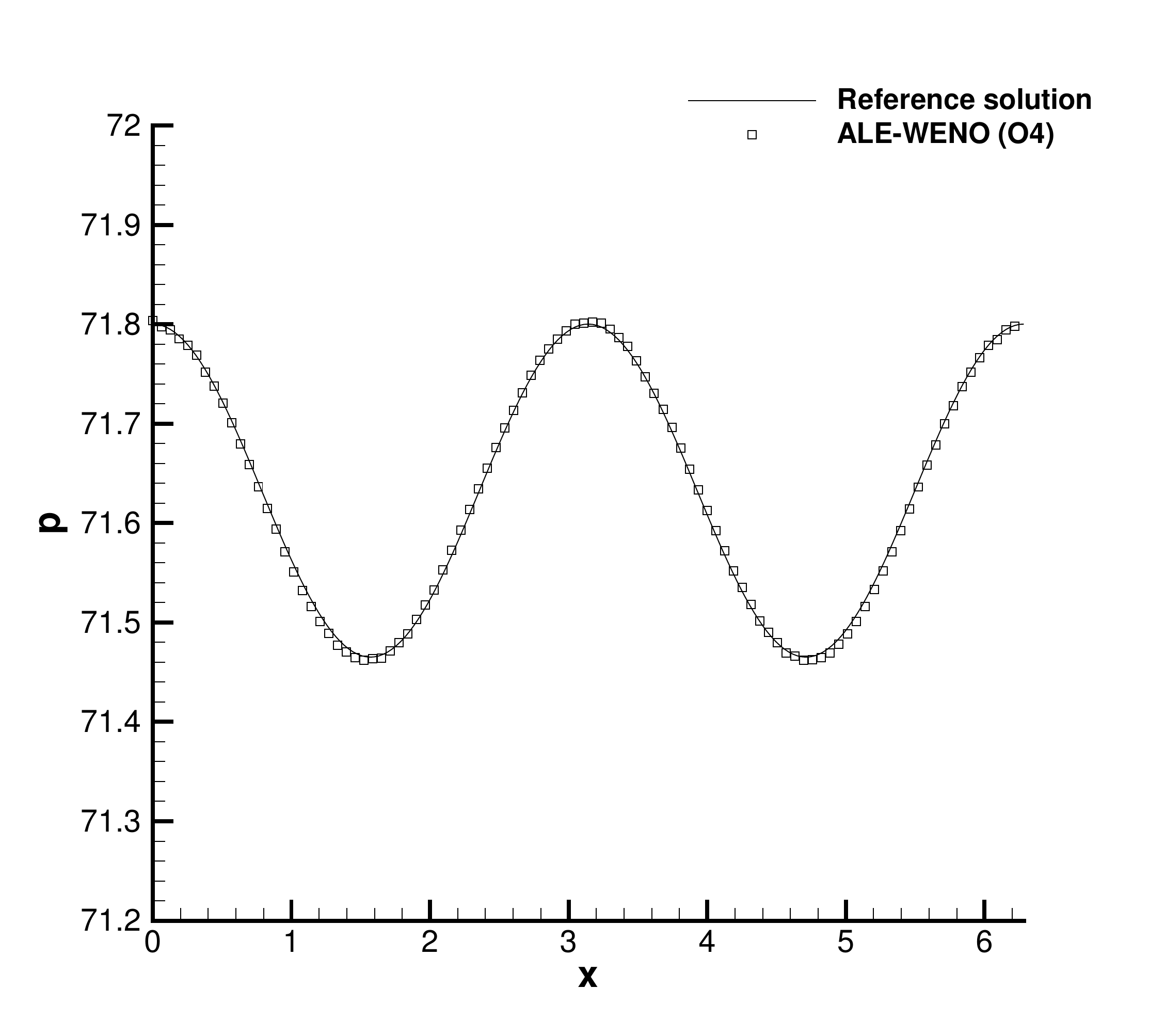}    
\end{tabular} 
\caption{Taylor-Green vortex with a viscosity of $\mu=10^{-1}$: exact solution of the Navier-Stokes equations and numerical solution for the hyperbolic model (HPR) at a final time of $t_f=1.0$ obtained with the direct ALE ADER-WENO fourth order scheme. Top: mesh configuration (left) and distortion tensor 
component $A_{11}$. Bottom: 1D cuts along the $x$ and the $y$ axis for velocity components $u$ and $v$ (left) and for the pressure $p$.  } 
\label{fig.tgv}
\end{center}
\end{figure}

A typical test problem used for the verification of numerical methods for the incompressible Navier-Stokes equations is the Taylor-Green vortex problem. An exact solution is available in two space dimensions, which is
\begin{eqnarray}
    u(x,y,t)&=&\sin(x)\cos(y)e^{-2\nu t},  \nonumber \\
		v(x,y,t)&=&-\cos(x)\sin(y)e^{-2\nu t}, \nonumber \\
    p(x,y,t)&=& C + \frac{1}{4}(\cos(2x)+\cos(2y))e^{-4\nu t},
\label{eq:TG_ini}
\end{eqnarray}
where $\nu=\frac{\mu}{\rho}$ represents the kinematic viscosity. The initial additive constant for the pressure field is given by $C=100/\gamma$ with the ratio of specific heats $\gamma=1.4$. The other parameters are chosen to be $\rho_0=1$, $c_v=1$, $c_s=10$ and the dynamic viscosity coefficient is set to $\mu=10^{-1}$. The computational domain is given by $\Omega(0)=[0;2\pi]^2$ with periodic boundaries imposed on each side and it is discretized with a total number of $N_E=5630$ triangles. The initial condition for velocity and pressure is given by \eqref{eq:TG_ini}, while the initial density and the distortion tensor are $\rho=\rho_0$ and $\mathbf{A}=\mathbf{I}$, respectively. The fourth order accurate numerical results are depicted in Figure \ref{fig.tgv} at the final time of the simulation $t_f=1.0$. An excellent agreement between the HPR model in the low Mach number regime and the exact solution of the incompressible Navier-Stokes equations can be observed, both for velocity and pressure.  We also plot the distortion tensor component $A_{11}$ which provides a useful and intuitive visualization of the flow. Moreover one can note that the mesh is adapted to the vortex structure of this problem. 

\subsection{Viscous shock problem} \label{sssec:problem_viscous_shock}
\begin{figure}[!htbp]
	\begin{center}
	\begin{tabular}{c}
	\includegraphics[width=0.85\textwidth]{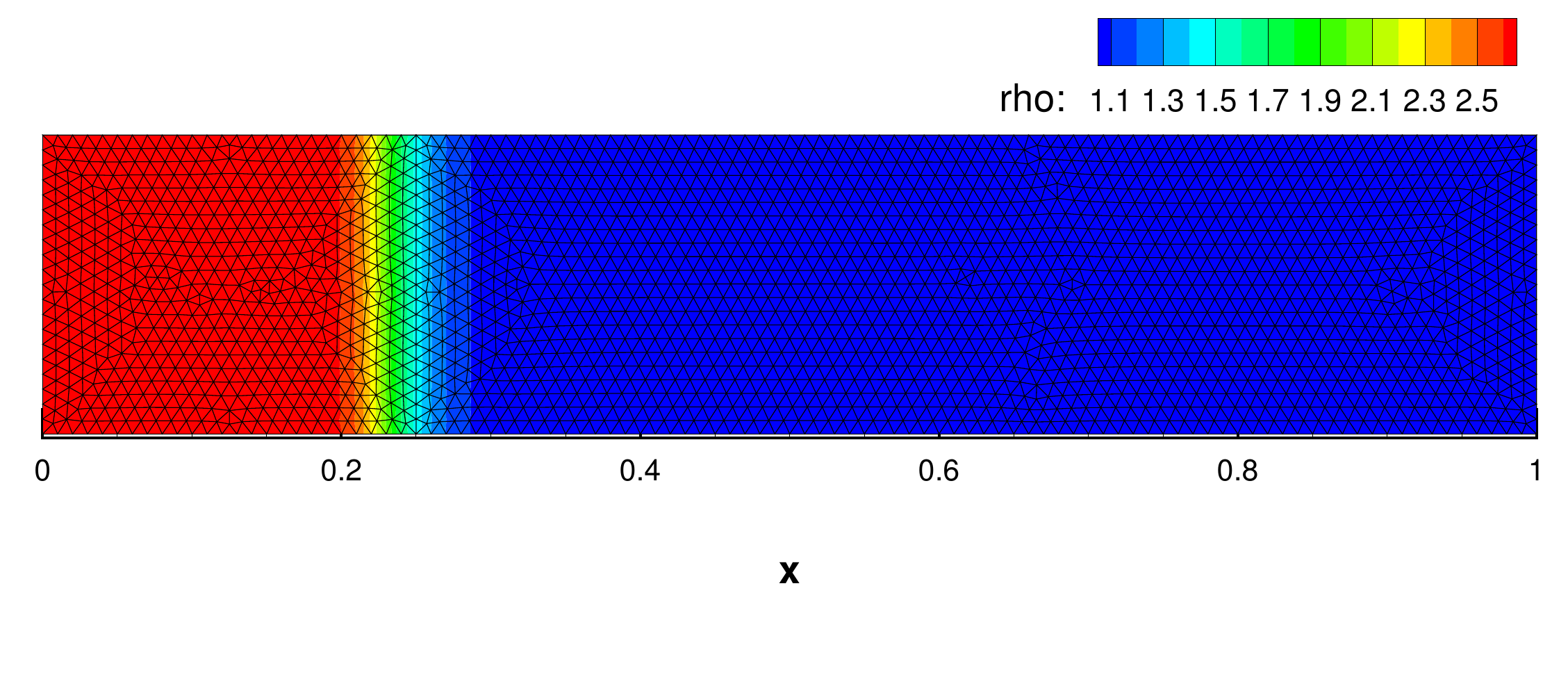} \\
	\includegraphics[width=0.85\textwidth]{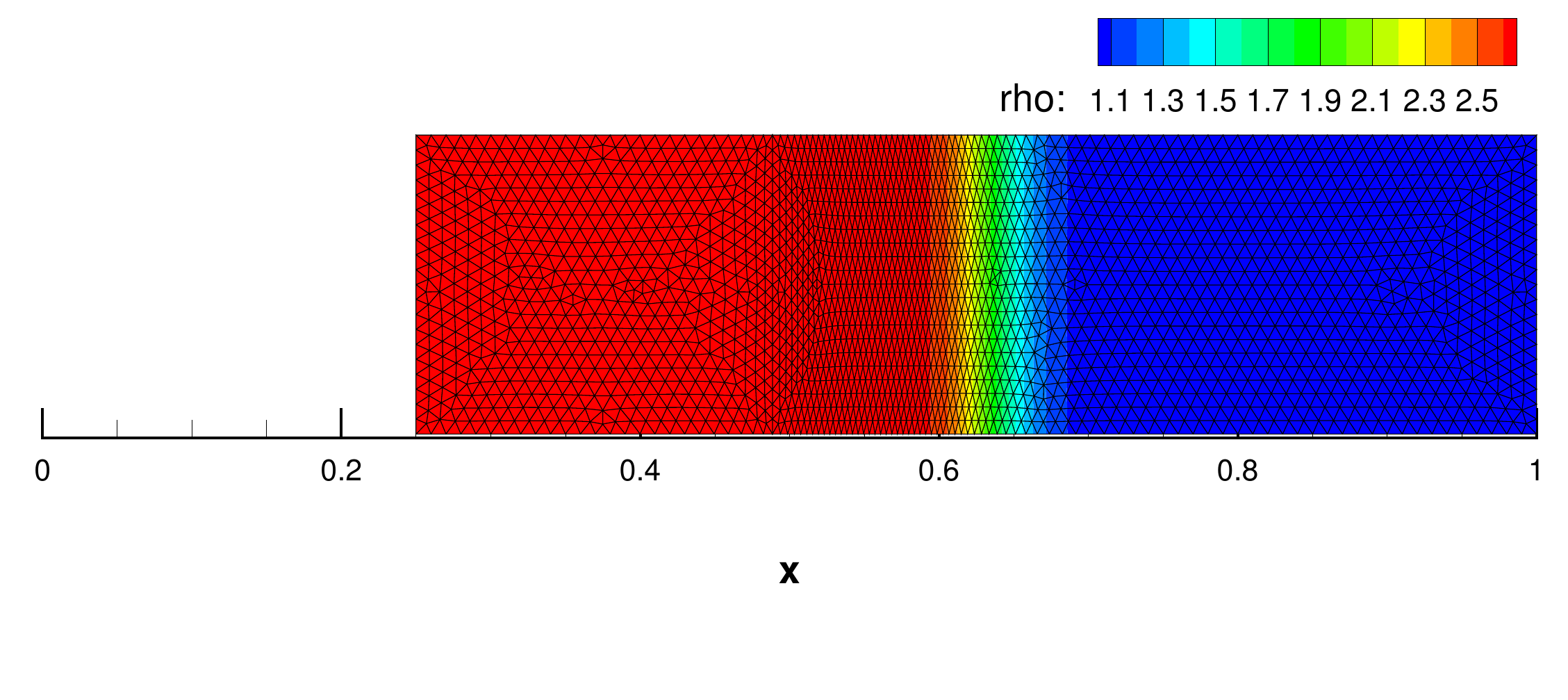}
	\end{tabular} 
	\caption{Viscous shock problem with shock Mach number $M_s=2$ and Prandtl number $Pr=0.75$. Initial (left) an final (right) mesh configuration and density distribution.} 
	\label{fig:viscous_shock_mesh}
	\end{center}
\end{figure}

\begin{figure}[!htbp]
	\begin{center}
	\begin{tabular}{cc}
	\includegraphics[width=0.47\textwidth]{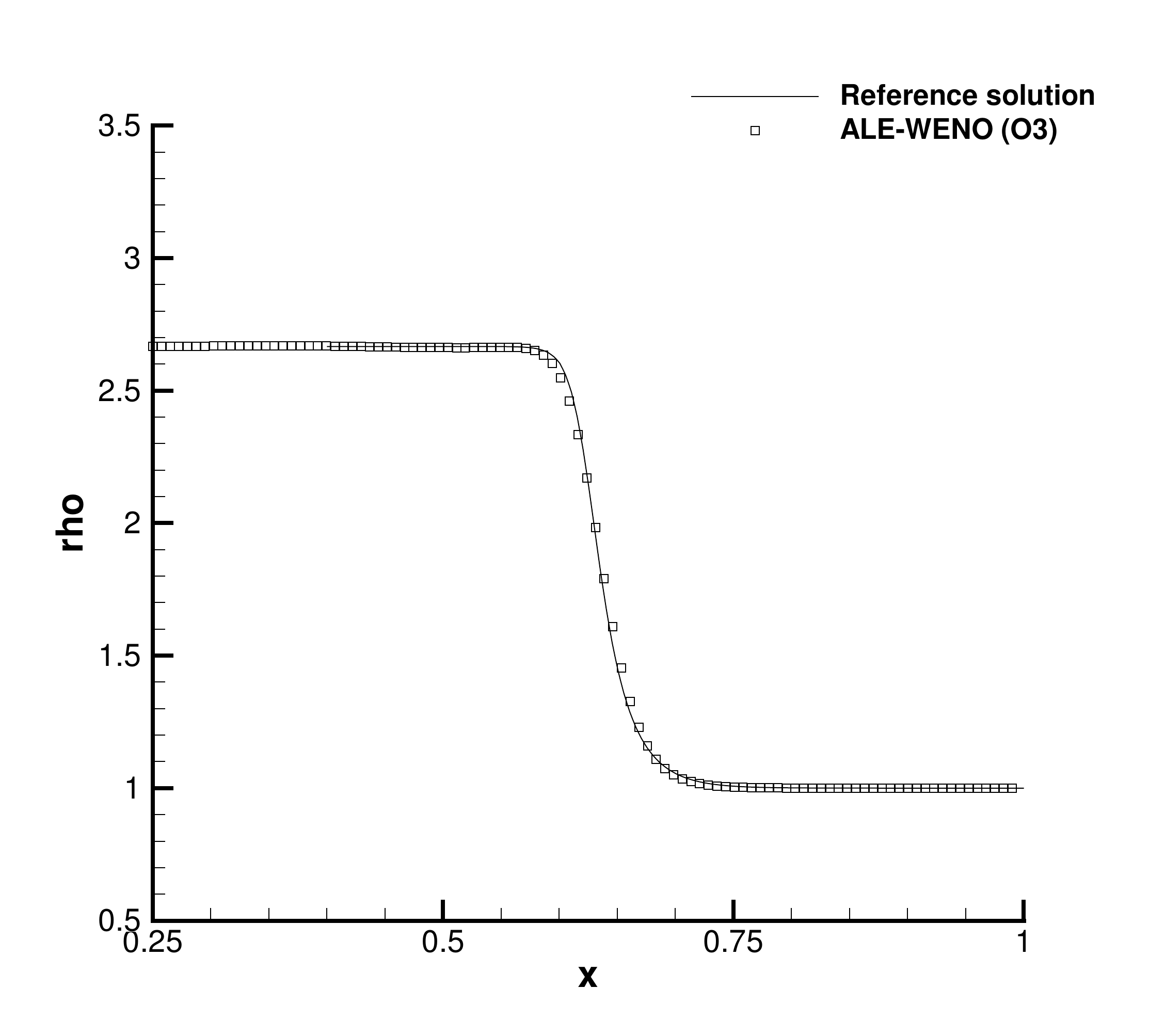} &
	\includegraphics[width=0.47\textwidth]{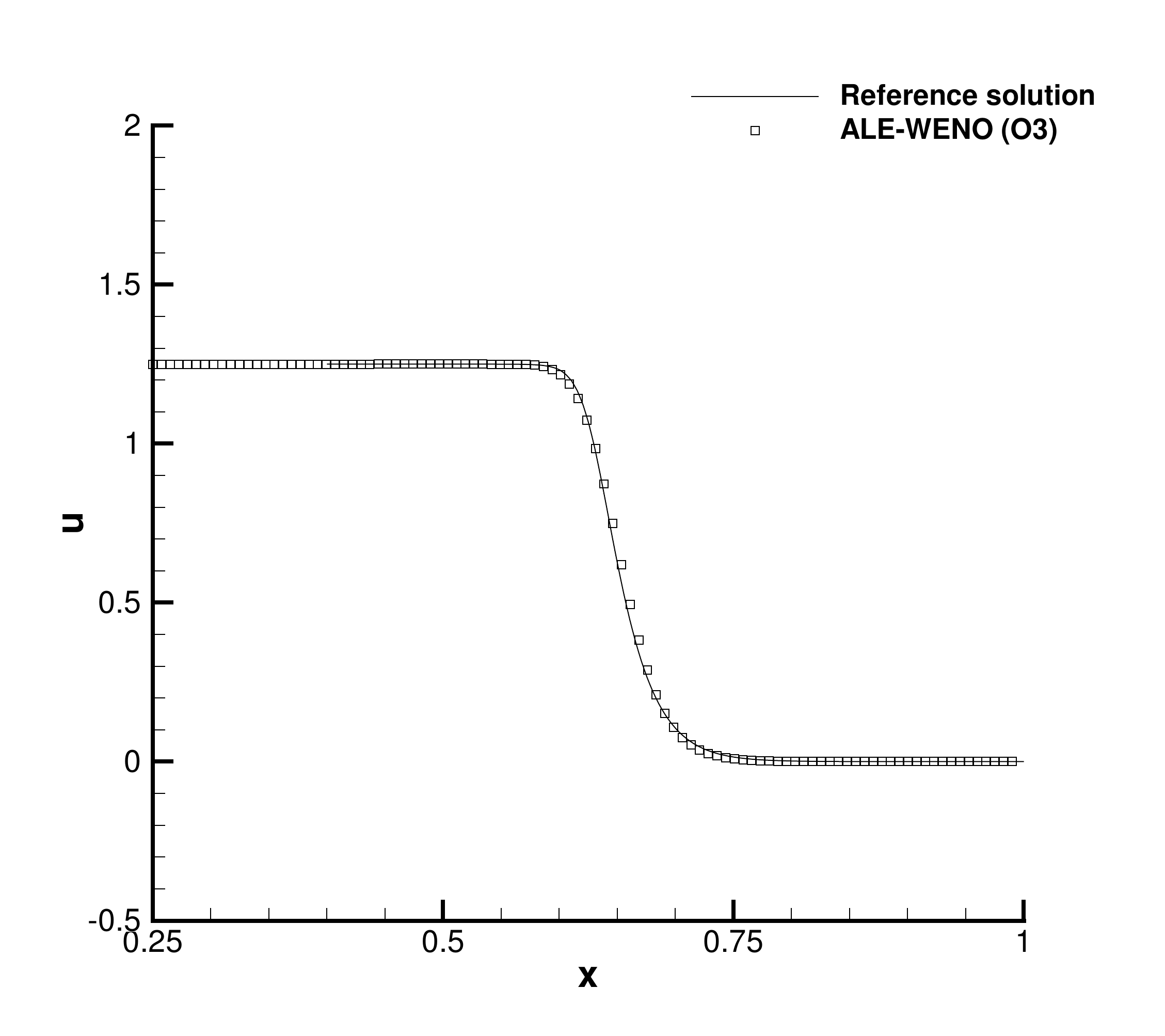} \\
	\includegraphics[width=0.47\textwidth]{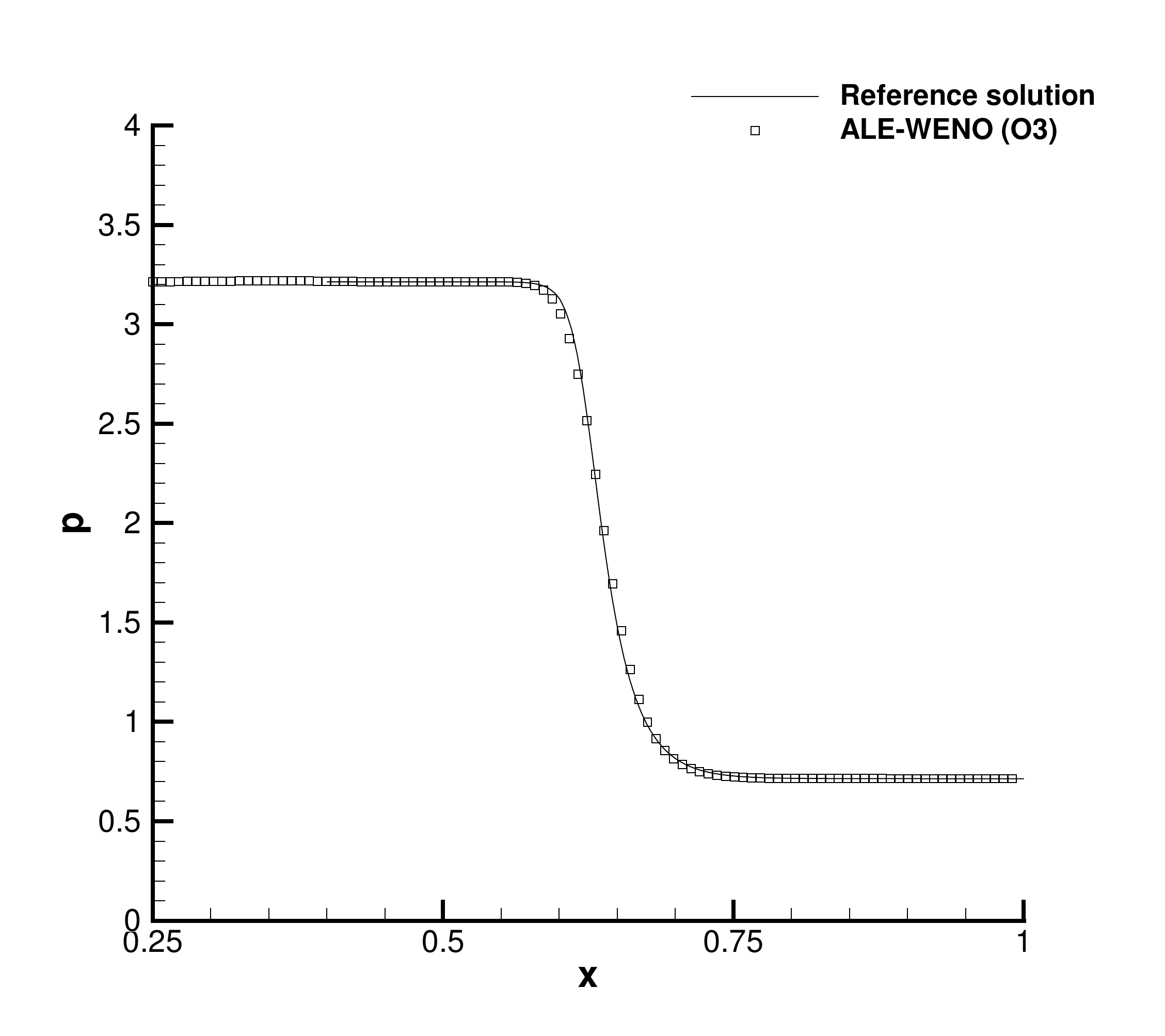} &
	\includegraphics[width=0.47\textwidth]{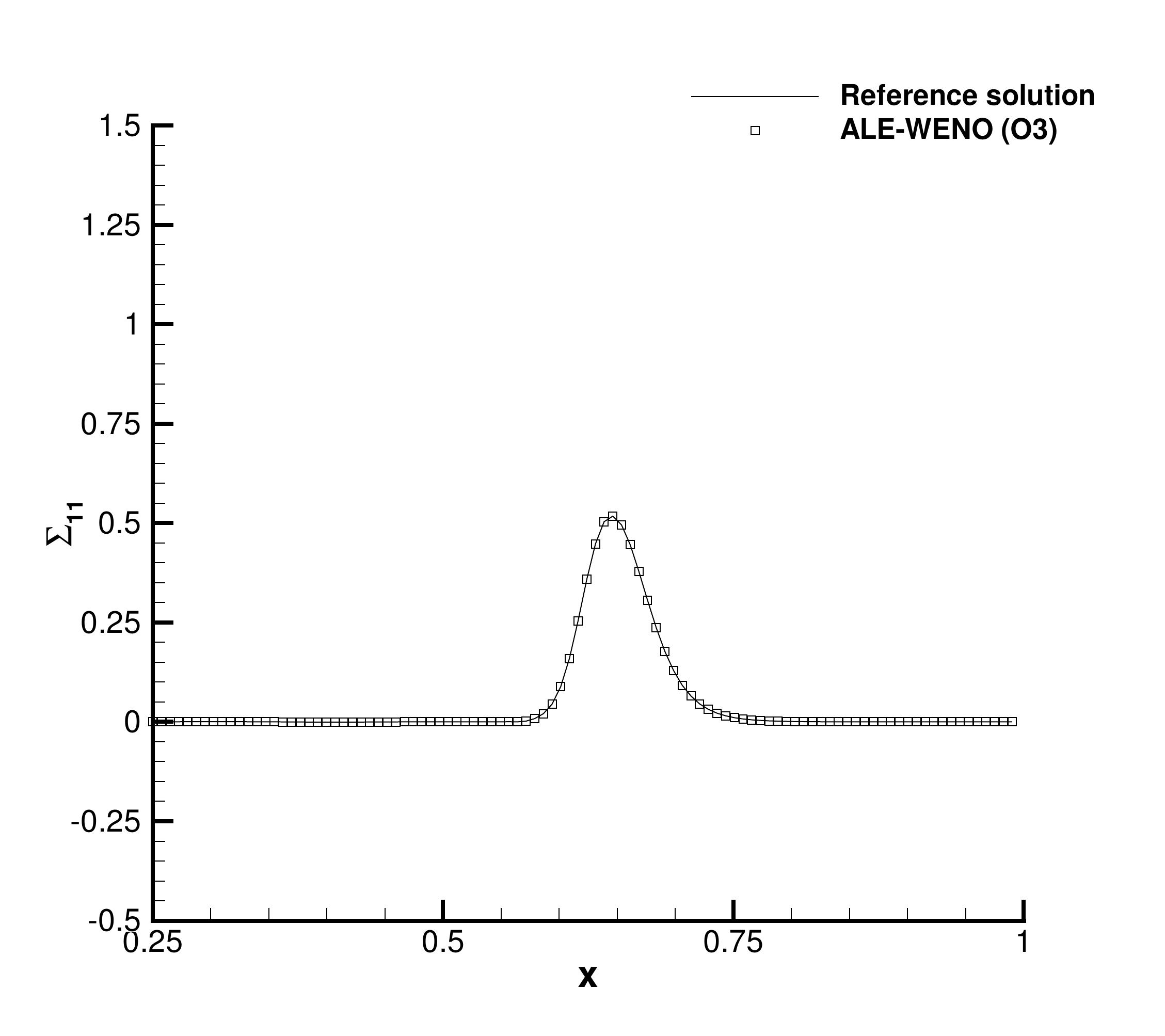}
	\end{tabular} 
	\caption{Viscous shock problem with shock Mach number $M_s=2$ and Prandtl number $Pr=0.75$. Comparison of the exact solution of the compressible Navier-Stokes equations according to Becker \cite{Becker1923} with the HPR model: Density, $x$ component of the velocity, pressure and viscous stress tensor component $\sigma_{11}$.} 
	\label{fig:viscous_shock_var}
	\end{center}
\end{figure}

The HPR model can also deal with supersonic viscous flows, therefore we propose to solve the problem of an isolated viscous shock wave which is traveling into a medium at rest with a shock Mach number of $M_s=2$. In \cite{Becker1923} an analytical solution for the compressible Navier-Stokes equations is derived for the special case of a stationary shock wave at Prandtl number $Pr= 0.75$ with constant viscosity. As done in \cite{Dumbser_HPR_16}, we superimpose a constant velocity field $u = M_s c_0$ to the previous stationary shock wave solution to obtain an unsteady shock wave traveling into a medium at rest. The computational domain is initially the rectangular box $\Omega(0)=[0;1]\times[0;0.2]$ which is paved with a set of non-overlapping triangles with characteristic mesh size $h=1/100$. No-slip wall boundary conditions are imposed everywhere, except on the left side of the domain where we let the piston move with the local fluid velocity. The initial condition involves a shock wave centered at $x=0.25$ propagating at Mach $M_s=2$ from left to right with a Reynolds number of $Re=100$. The upstream shock state is defined by $\rho_0=1$, $u_0=v_0=0$, $p_0=1/\gamma$ and $c_0=1$. and the parameters of the HPR model are $\gamma=1.4$, $c_v=2.5$, $c_s=50 $, $\mu=2\times 10^{-2} $. In this case we also consider the heat flux, hence setting initially $\JJ=\mathbf{0}$ with $\alpha=50 $, $T_0=1 $, and $\kappa=9/3 \times 10^{-2} $. The distortion tensor is initialized to $\A =\sqrt[3]{\rho} \II $ and the final time of the simulation is $t_f=0.2$ with the shock front located at $x=0.65$. Figure \ref{fig:viscous_shock_mesh} depicts the mesh configuration and the density distribution at the initial and at the final time, while in Figure \ref{fig:viscous_shock_var} one can note an excellent agreement of the third order ADER-WENO-ALE solution with the analytical solution of the compressible Navier-Stokes equations \cite{Becker1923}. We compare the exact solution and the
numerical density, $x$ component of the velocity, pressure and viscous stress tensor component $\sigma_{11}$ from top-left panel to bottom-right one

\subsection{Cylindrical explosion problem} \label{sssec:viscous_EP}
\begin{figure}[!htbp]
	\begin{center}
	\begin{tabular}{cc}
	\includegraphics[width=0.47\textwidth]{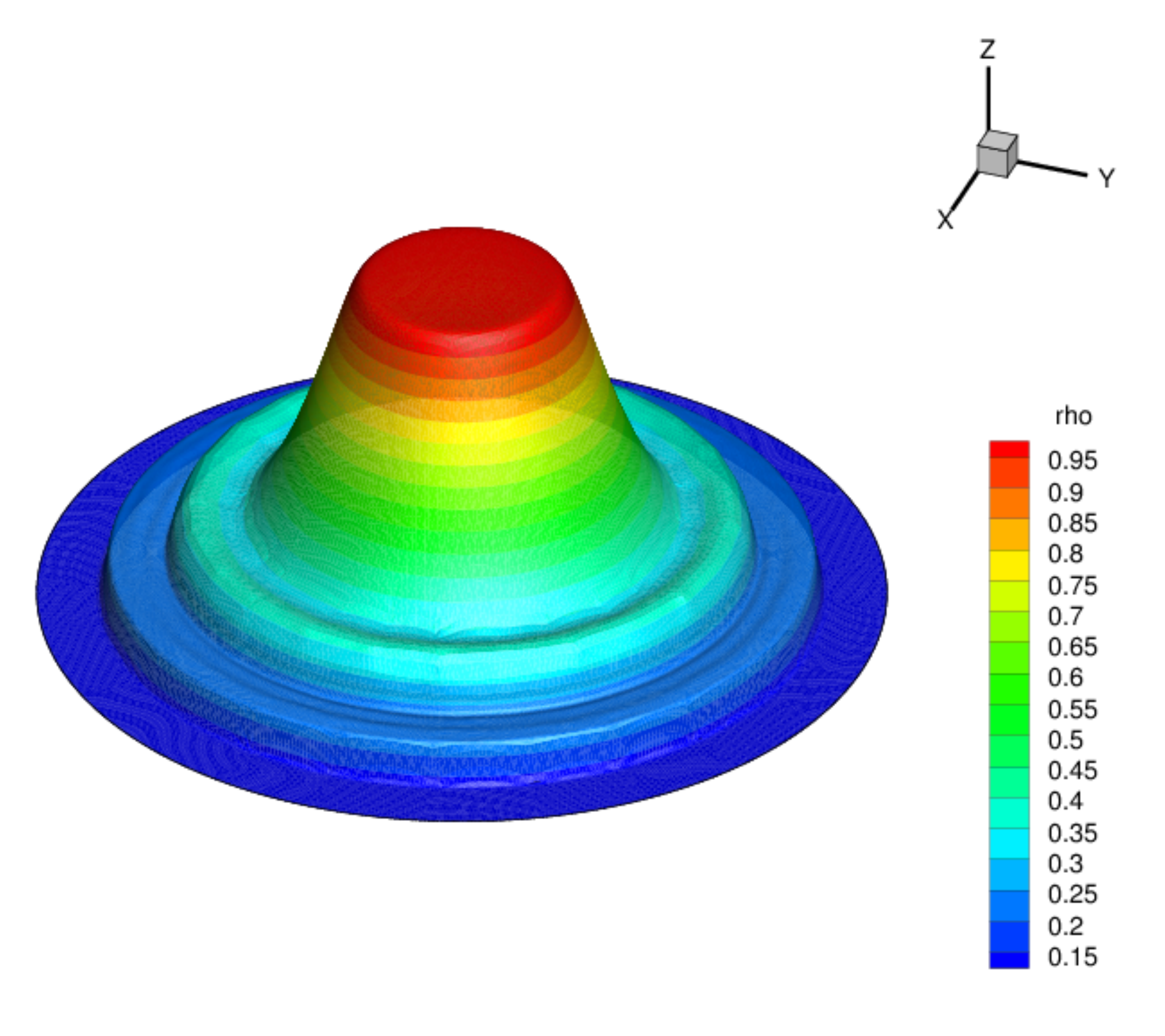} &
	\includegraphics[width=0.47\textwidth]{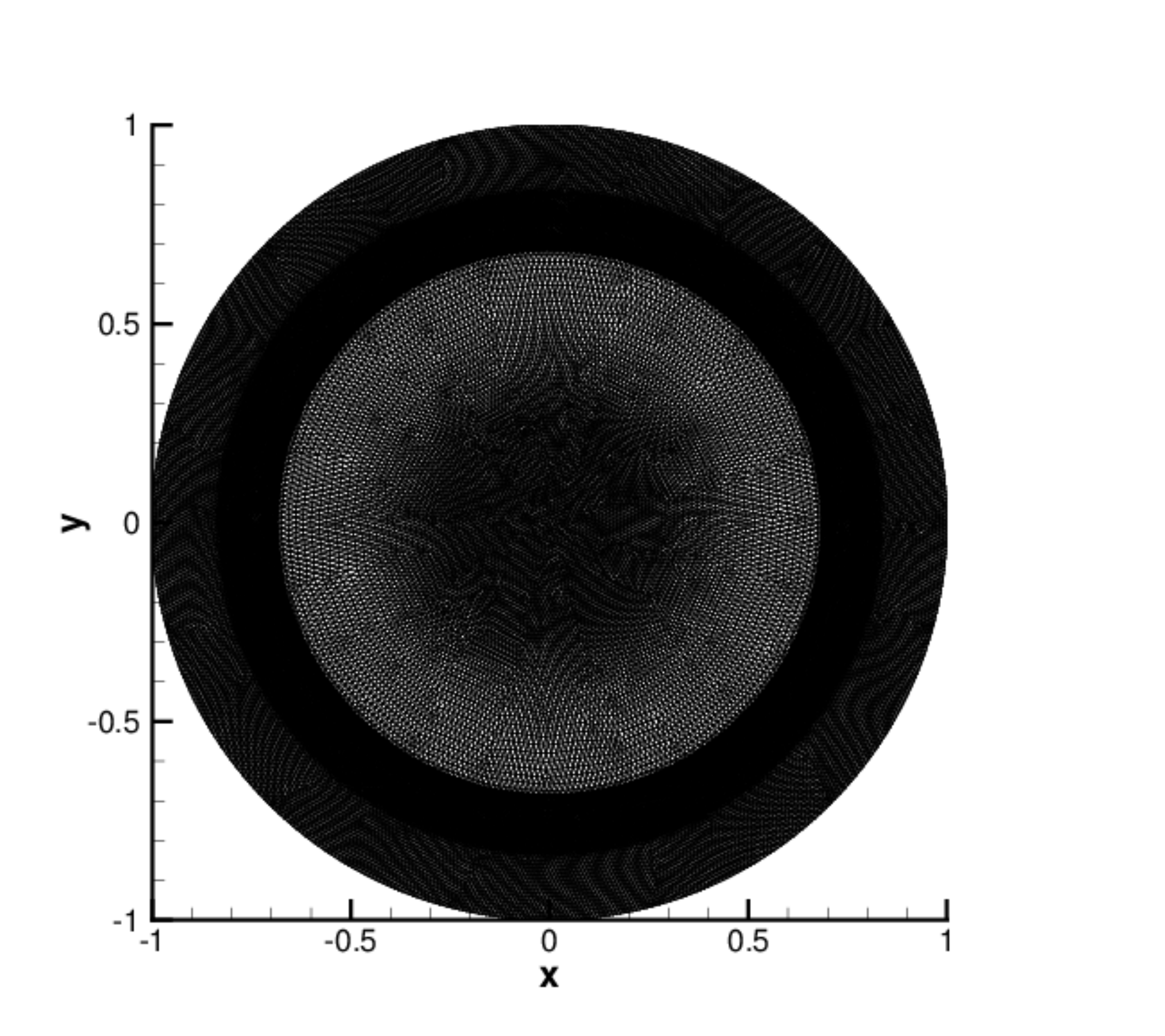} \\
	\includegraphics[width=0.47\textwidth]{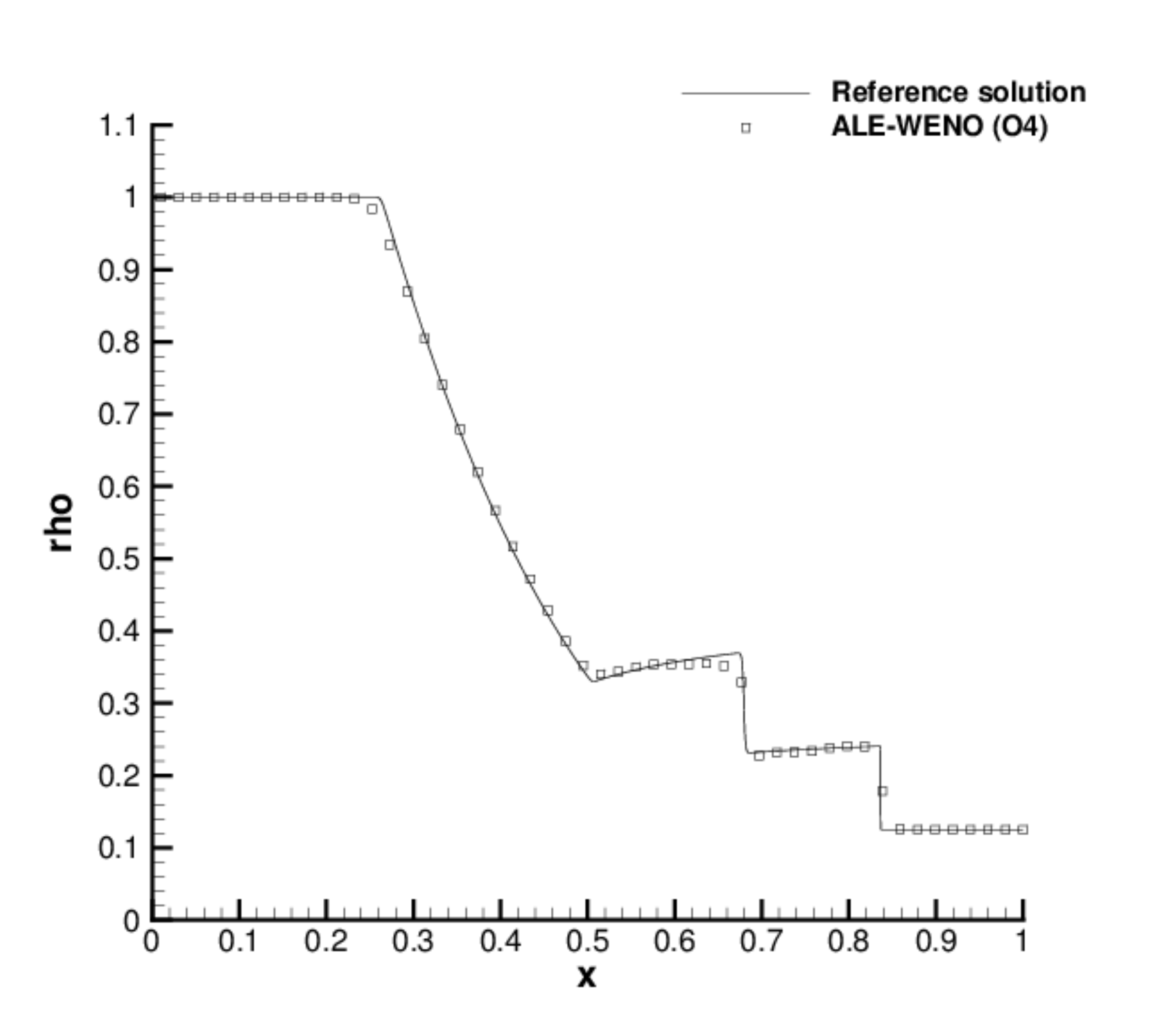} &
	\includegraphics[width=0.47\textwidth]{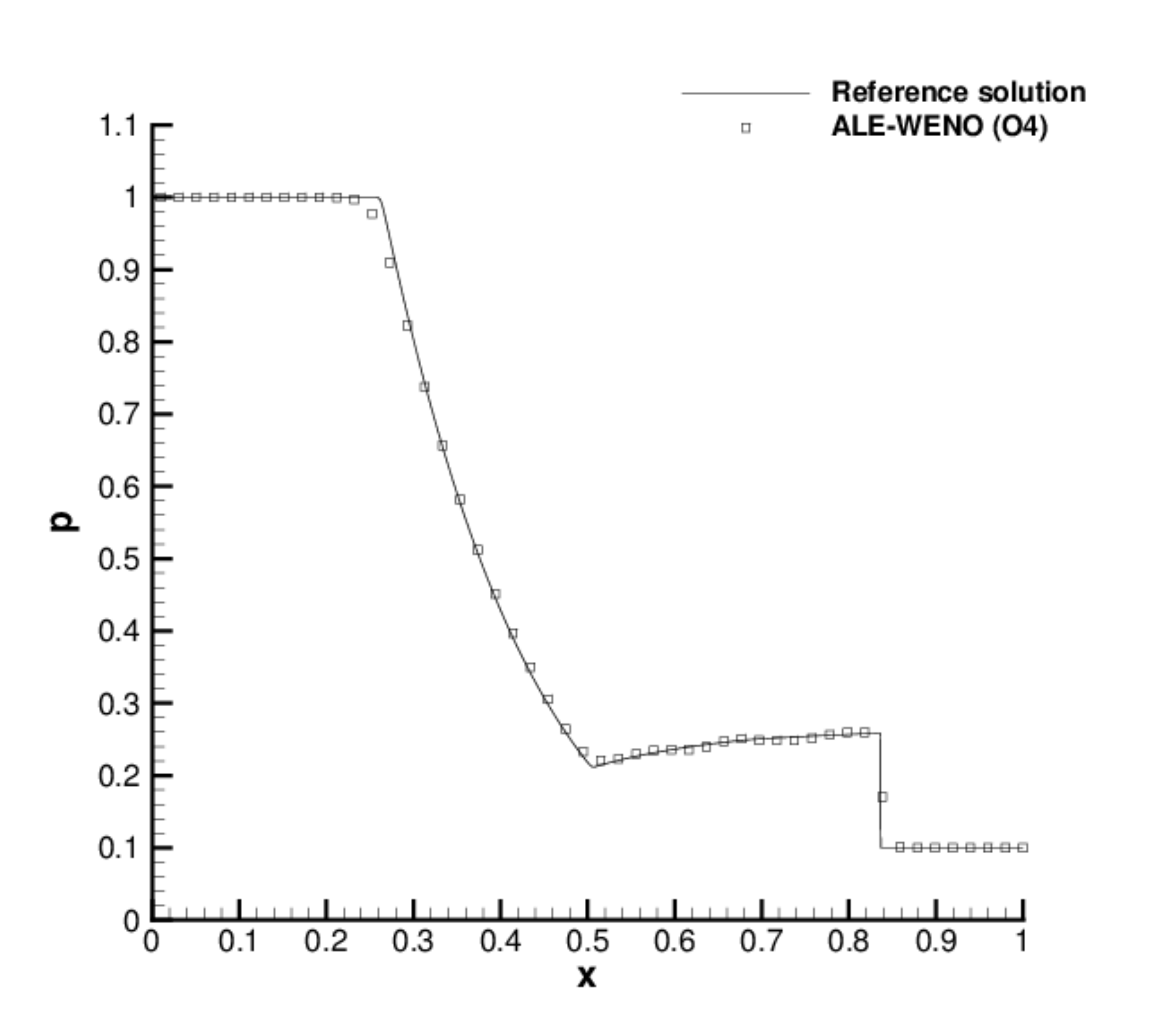}
	\end{tabular} 
	\caption{Cylindrical explosion problem. Results obtained with the HPR model at the final time $t=0.2$ obtained with a fourth order ADER-WENO-ALE scheme with $\mu=k=10^{-4}$. Top: three-dimensional density    distribution and final mesh configuration. Bottom: cut along the $x$-axis for density and pressure compared with the inviscid reference solution (compressible Euler equations).} 
	\label{fig:EP2D}
	\end{center}
\end{figure}

Here we present numerical results for a cylindrical explosion problem solved with the HPR model. The initial computational domain $\Omega(0)$ is the circle of radius $R=1.0$ and the initial condition is given by two constant states separated by a discontinuity located at radius $R_s=0.5$. Therefore the fluid is initially assigned with the corresponding primitive 
state vector $\mathbf{V}=(\rho, u, v, A_{11}, A_{12}, A_{21}, A_{22}, A_{33}, p, J_{1}, J_{2})$ which reads
\begin{equation}
  \mathbf{V}(\x,0) = \left\{ \begin{array}{ccc} (1,0,0,1,0,0,1,1,1,0,0) & \textnormal{ if } & r \leq R_s, \\ 
                                        (0.125,0,0,0.5,0,0,0.5,0.5,0.1,0,0) & \textnormal{ if } & r > R_s,        
                      \end{array}  \right. 
\end{equation}
where the generic radial coordinate is $r=\sqrt{x^2+y^2}$. The initial distortion tensor has been set to $\AAA=\sqrt[3]{\rho}$, while the initial thermal impulse vector is $\mathbf{J}=\mathbf{0}$. Transmissive boundary conditions are imposed on the external boundary and the mesh is composed by $N_E=68324$ triangles. The final time of the simulation is chosen to be $t_f=0.2$ and the parameters for the HPR model are $\gamma=1.4$, $c_v=2.5$, $c_s=0.5$, $\rho_0=1$, $\alpha=0.5$, $\mu=k=10^{-4}$. The reference solution can be computed by solving the one-dimensional compressible Euler equations with a geometric source term that takes into account the cylindrical geometry, as fully detailed in \cite{ToroBook,Lagrange2D}. We use a second order MUSCL scheme with the Rusanov flux on a one--dimensional mesh of 15000 points in the radial interval $r \in [0;1]$ to solve the inhomogeneous system and this solution is assumed to be our reference solution. We run a fourth order scheme to obtain the numerical results depicted in Figure \ref{fig:EP2D}, where one can note a good agreement with the reference solution for the 1D cut along the $x$-axis representing density and pressure. Furthermore we plot also the final mesh configuration which highlights the strongly compressed cells at the shock location and the stretched elements crossed by the rarefaction wave traveling towards the center of the domain.

\subsection{Purely elastic Riemann problems} \label{sssec:elasticRP}

In this section we consider the equations of nonlinear elasticity \cite{GodunovRomenski72,Godunov:1995a,Godunov:2003a}, which can be retrieved by the HPR model in the limit $\tau_1 \to \infty$ with $\sigma_0 > 0$. We set up two shock tube problems on a 2D domain $\Omega(0)=[0;1]\times[0;0.1]$ where periodic
boundary conditions are applied in $y$-direction and transmissive boundaries are imposed in $x$-direction. The initial discontinuity located at $x=0.5$ separates the two initial states, given in terms of primitive variables and entropy in Table \ref{tab:data_elast}. The material is copper, described by the parameters given in Table \ref{tab:dataEOS}, and the equation of state considered in this case is a complicated function of the internal energy and the distortion tensor, explicitly detailed in \cite{TitarevRomenskiToro,DumbserPNPM}. 

\begin{table}
\caption{Initial condition for the left state (L) and the right state (R) for the Riemann problems of nonlinear elasticity solved with the HPR model.}
\renewcommand{\arraystretch}{1.0}
\begin{center}
\begin{tabular}{ccccccccc}
\hline
			& $u$ 	& $v$		& $A_{11}$ 	& $A_{12}$ 	& $A_{21}$ 	& $A_{22}$ 	& $A_{33}$ 	& $s$  \\
\hline 
\multicolumn{1}{l}{\textbf{RP1 \cite{DumbserPNPM}:} } &  & & & & & & \\
\hline 
L 		& 0.0   & 0.0   & 0.95  		& 0.0 			& 0.0   		&  0.0 			& 1.0 			& 0.001 \\
R 		& 0.0   & 0.0   & 1.0  			& 0.0 			& 0.0   		&  1.0 			& 1.0 			& 0.0    \\
\hline 
\multicolumn{1}{l}{\textbf{RP2 \cite{DumbserPNPM}:} } &  & & & & & & \\
\hline 
L 		& 0.0   & 1.0   & 0.95  		& 0.0 			& 0.05  		&  1.0 			& 1.0 			& 0.001 \\
R 		& 0.0   & 0.0   & 1.0  			& 0.0 			& 0.0   		&  1.0 			& 1.0 			& 0.0    \\
\hline 
\end{tabular}
\end{center}
\label{tab:data_elast}
\end{table}

\begin{figure}[!htbp]
	\begin{center}
	\begin{tabular}{cc}
	\includegraphics[width=0.47\textwidth]{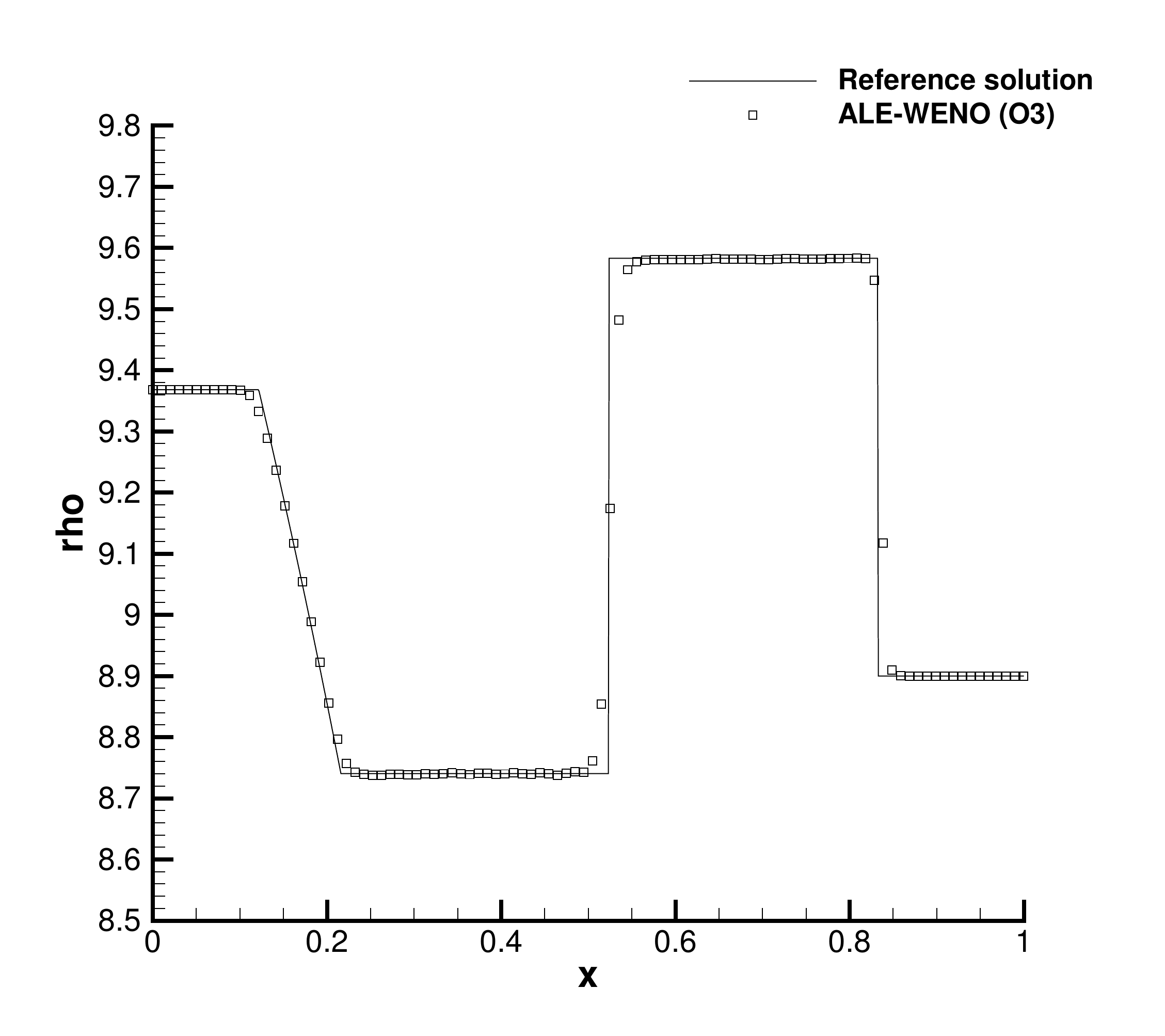} &
	\includegraphics[width=0.47\textwidth]{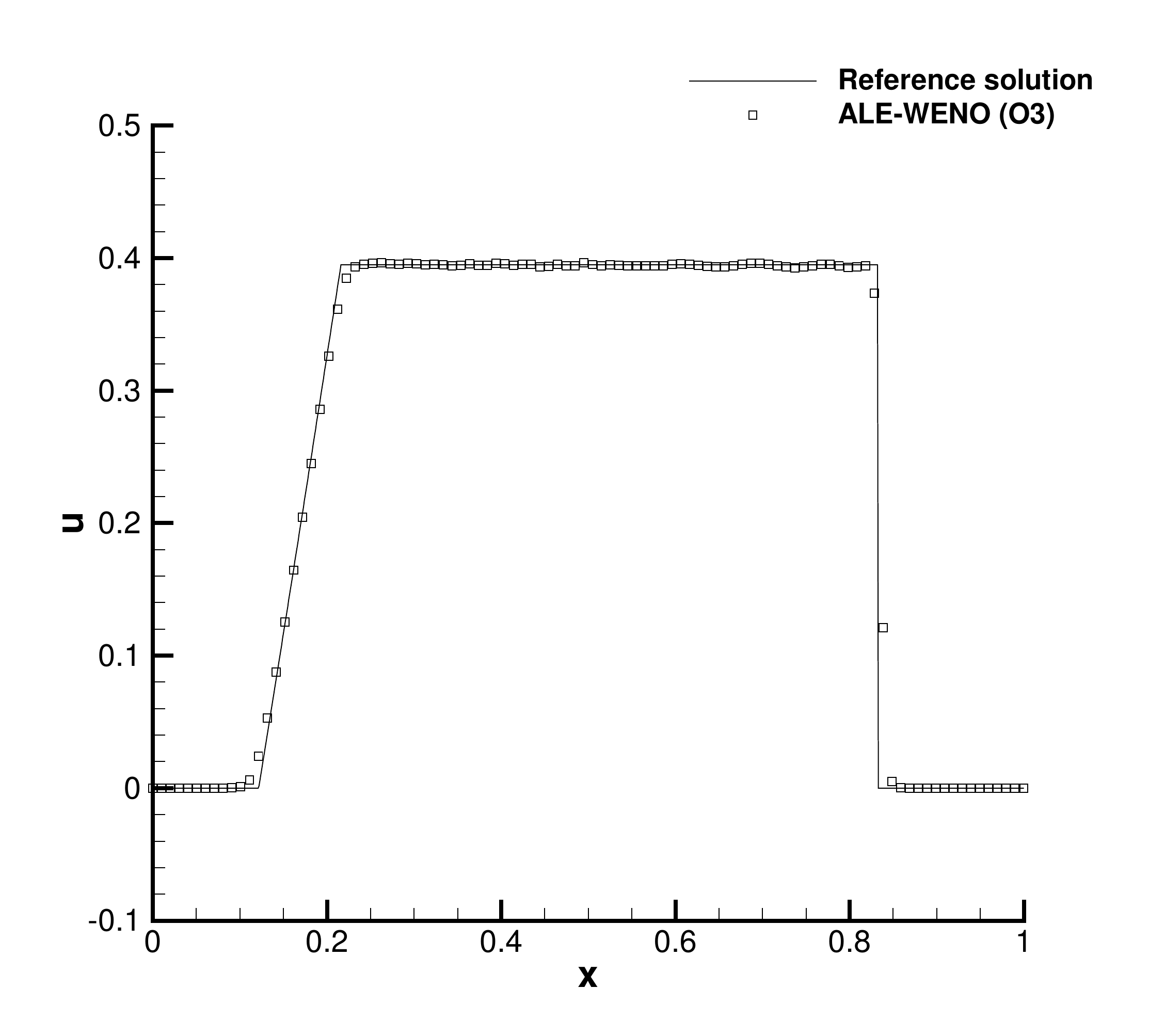}   \\
	\includegraphics[width=0.47\textwidth]{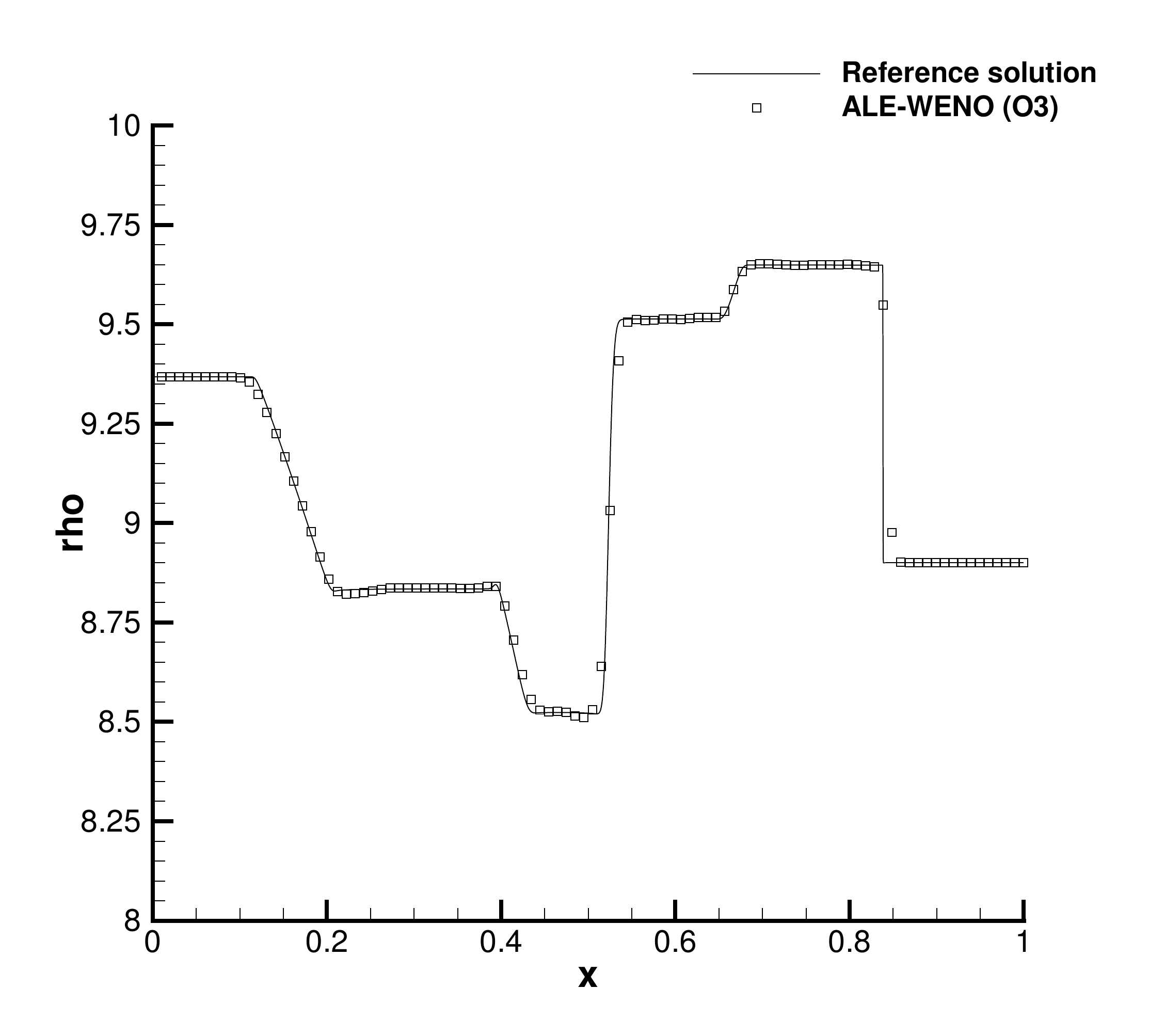} &
	\includegraphics[width=0.47\textwidth]{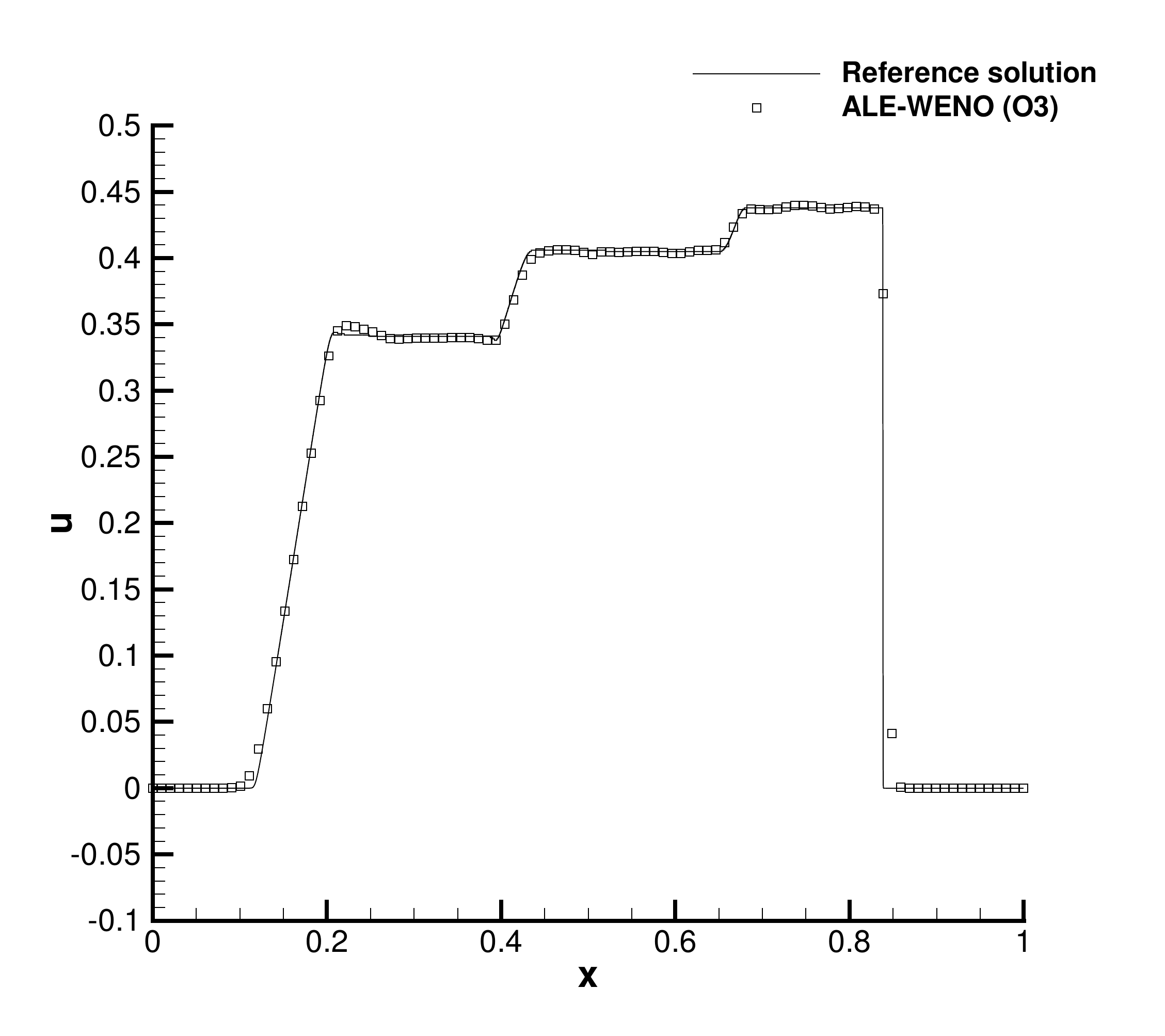} 
	\end{tabular} 
	\caption{Purely elastic Riemann problem $1$ (top) and $2$ (bottom). A 1D cut through the numerical solution at $y=0.025$ is plotted for density (left) and horizontal velocity component $u$ (right).} 
	\label{fig:NLE}
	\end{center}
\end{figure}

\begin{figure}[!htbp]
	\begin{center}
	\begin{tabular}{cc}
	\includegraphics[width=0.47\textwidth]{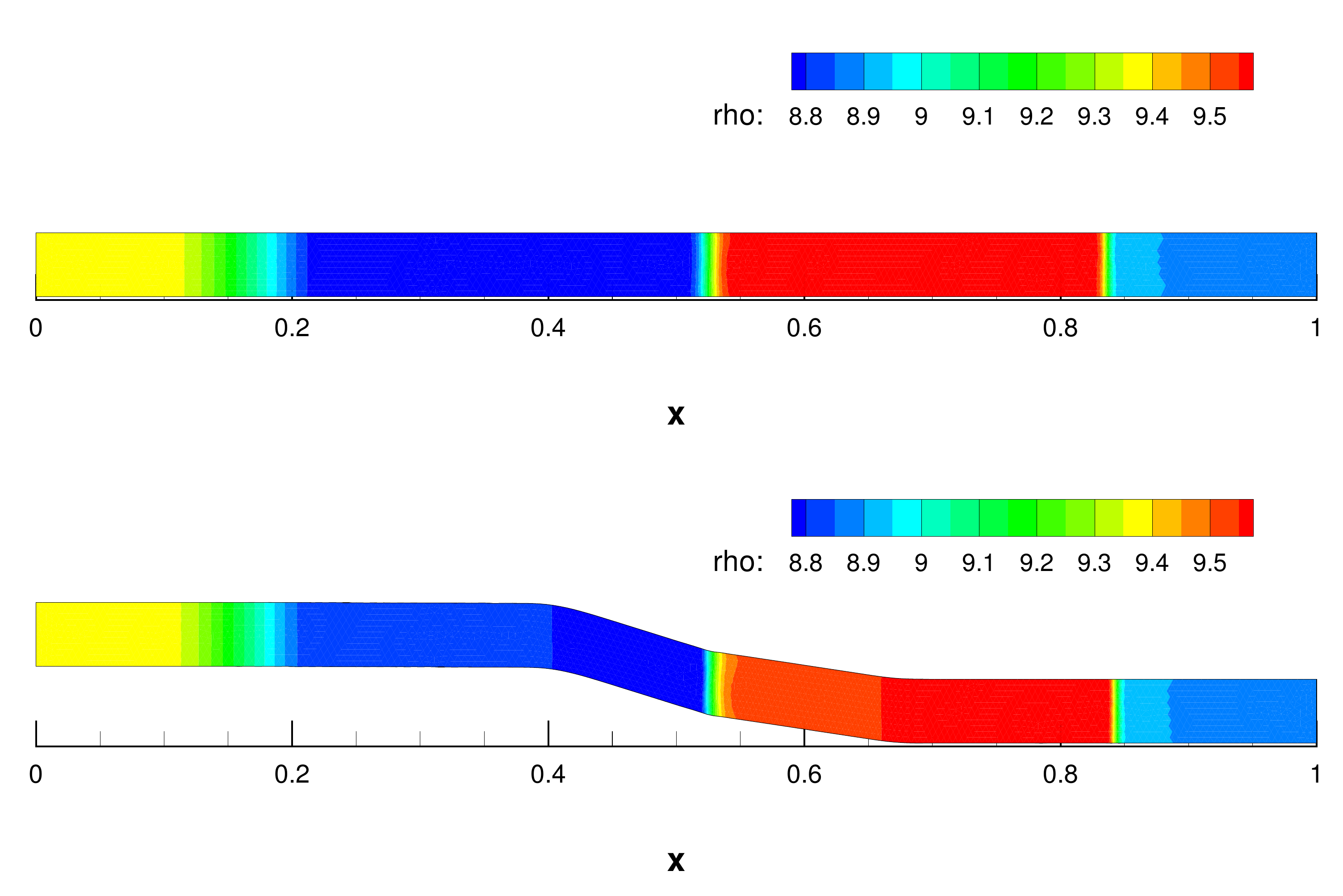} &
	\includegraphics[width=0.47\textwidth]{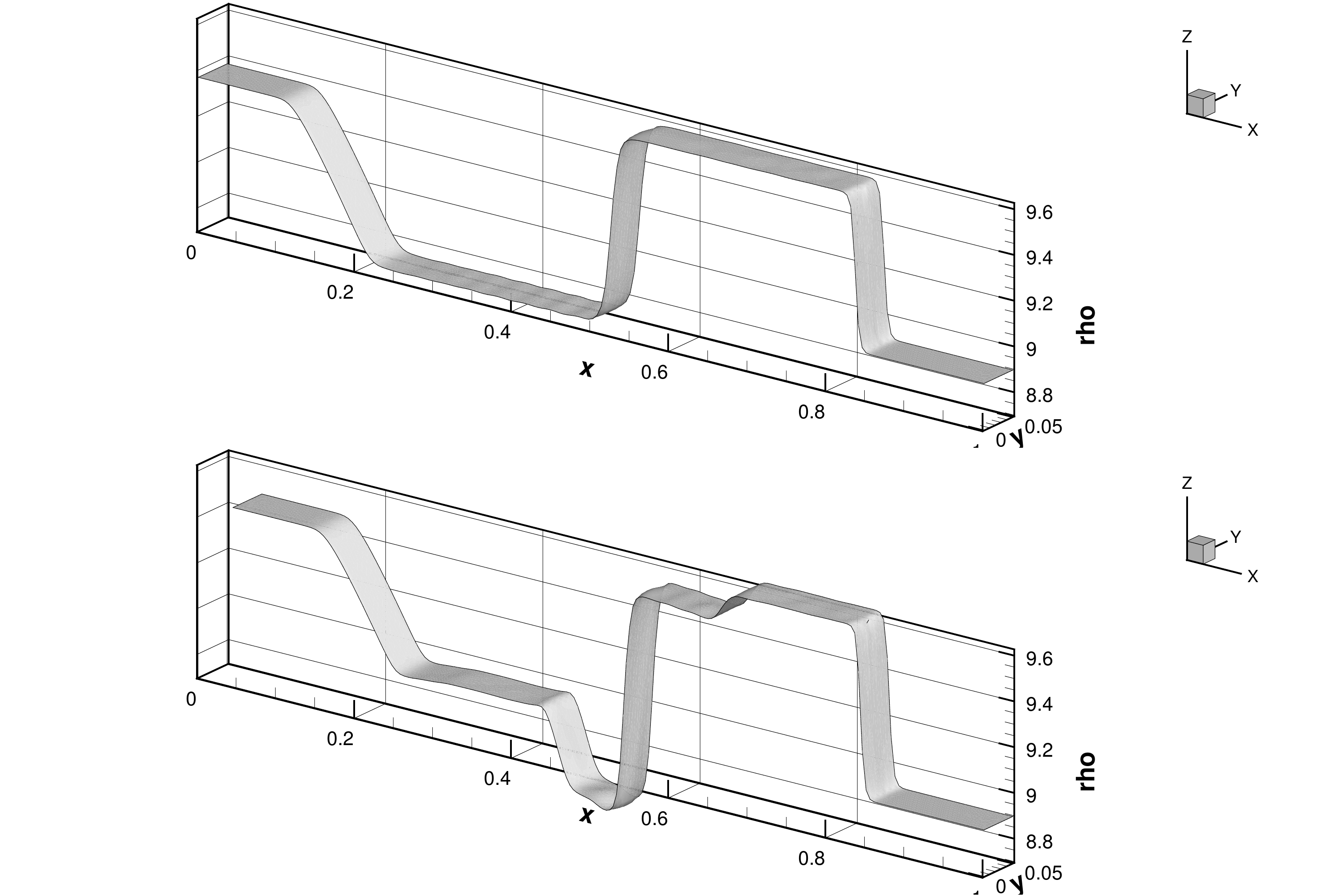} \\
	\end{tabular} 
	\caption{Purely elastic Riemann problem $1$ (top) and $2$ (bottom). Density distribution (left) and density elevation (right).} 
	\label{fig:NLE2}
	\end{center}
\end{figure}

The initial density is the reference density for copper, i.e. $\rho=8.930$, and we impose $c_v=0.4 \cdot 10^{-3}$. The final time of both simulations is $t_f=0.06$ and the numerical results are shown in Figures \ref{fig:NLE} and \ref{fig:NLE2}. Riemann problem 1 (RP1) corresponds to the three-wave shock tube problem, while RP2 considers a five-wave shock tube problem, originally proposed in 
\cite{TitarevRomenskiToro}. One can note a very good agreement between the numerical results obtained with a third order ADER-WENO-ALE scheme and the analytical solution of the nonlinear 
hyperelasticity model provided in \cite{TitarevRomenskiToro,Barton2009}. In Figure \ref{fig:NLE2} we can see for RP2 that the $y$ motion of the domain is not uniform but the waves are 
accurately maintained in their 1D shape. This may be a concern when using a moving mesh technique, but our approach seems to properly deal with this situation. 

\subsection{Elastic-plastic piston problem} \label{sssec:prec}
\begin{figure}[!htbp]
	\begin{center}
	\begin{tabular}{cc}
	\includegraphics[width=0.47\textwidth]{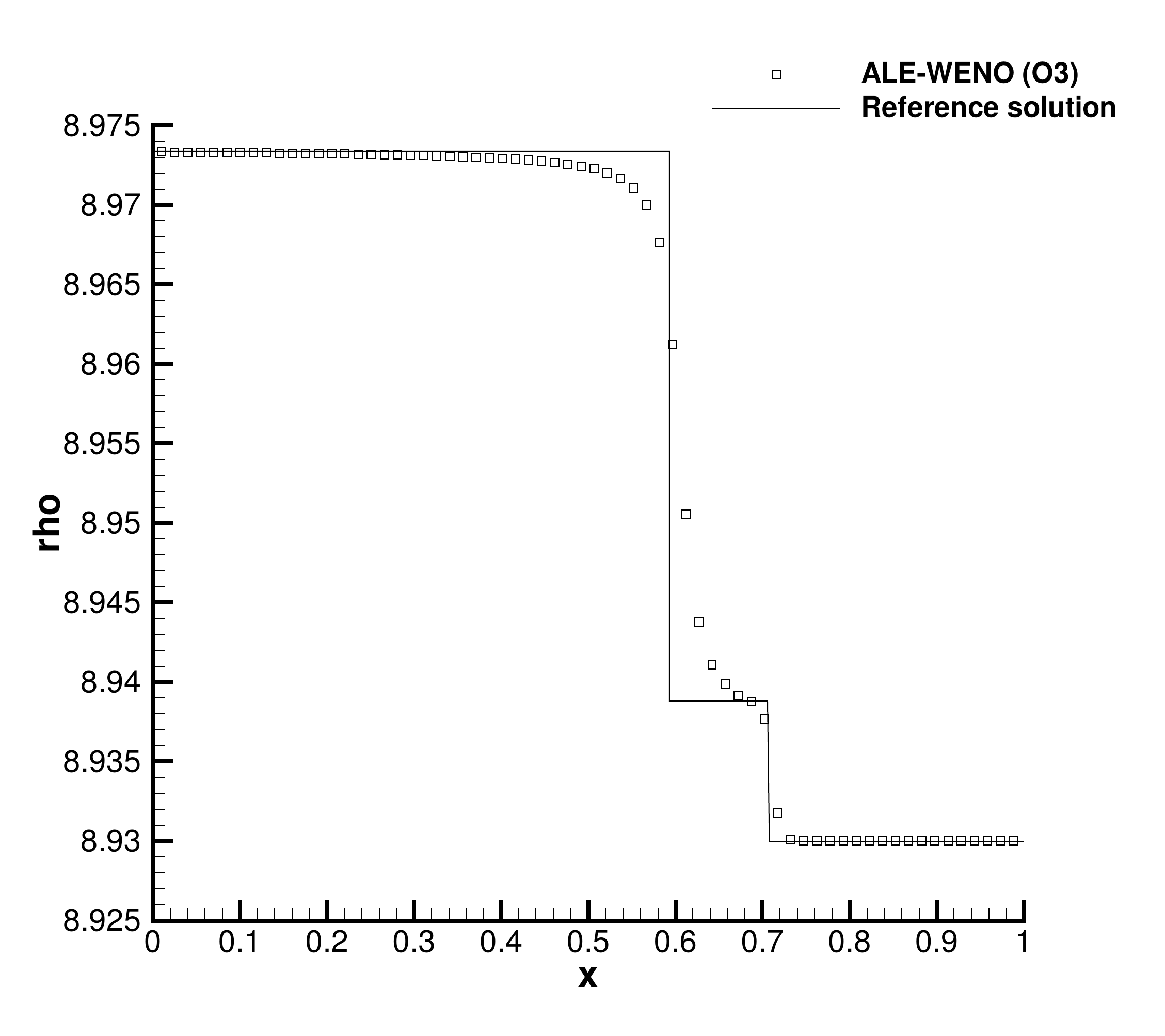} &
	\includegraphics[width=0.47\textwidth]{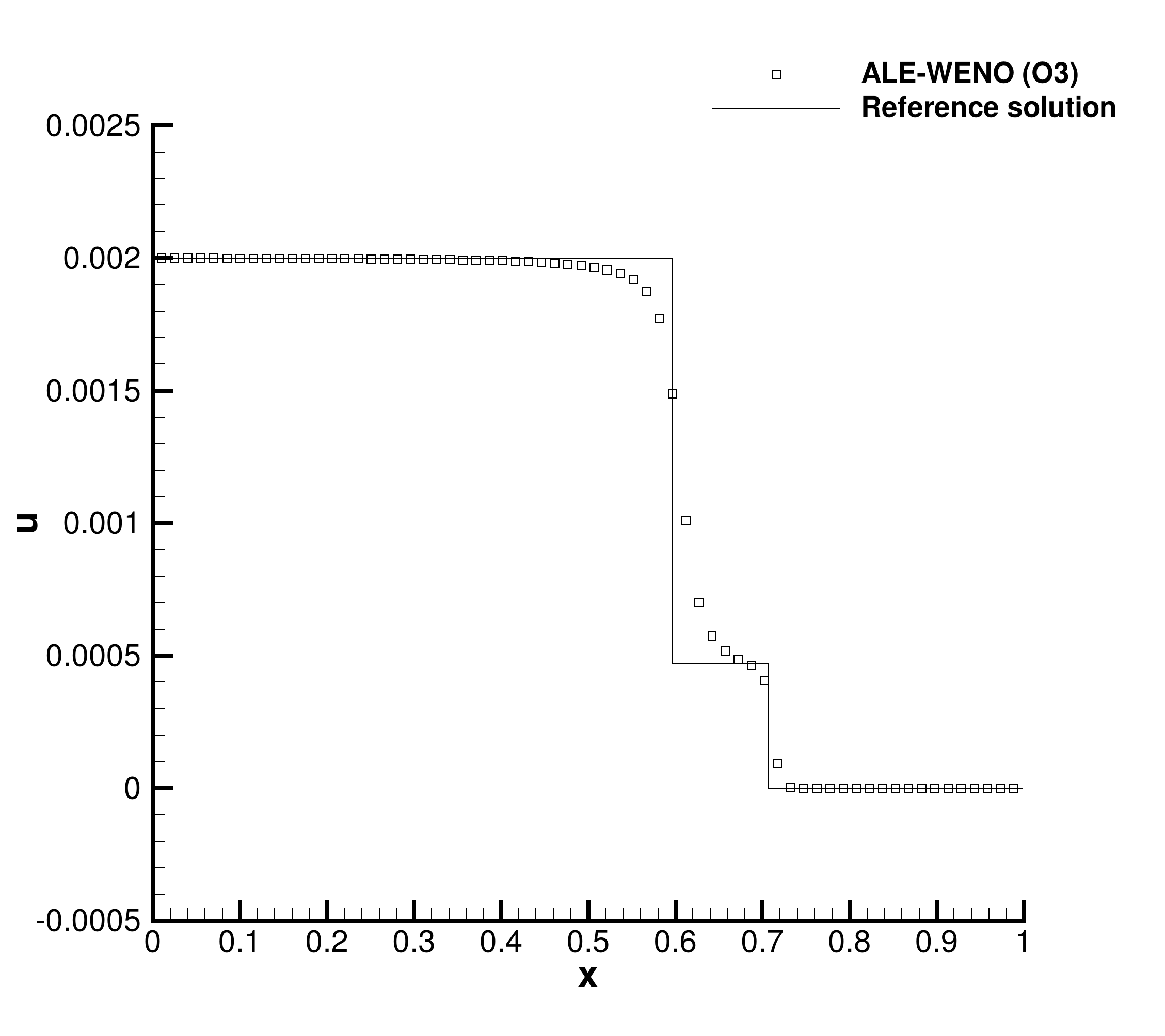}   \\
	\end{tabular} 
	\caption{Elastic-plastic piston problem at $t_f=1.5$. Density (left) and horizontal velocity (right) profiles obtained with the HPR model compared against the exact solution of the 
	ideal plasticity model  (straight line) using a third order accurate ADER-WENO-ALE scheme.} 
	\label{fig:EPP}
	\end{center}
\end{figure}

This test case is a one-dimensional flow characterized by a slope discontinuity which yields a two-wave structure with a first elastic shock wave, typically called the elastic precursor, followed by a plastic shock wave \cite{zeldovich67,Udaykumar2003}. An analytical solution is available and we refer the reader to \cite{Maire_elasto_13} for an exhaustive description. The material employed in this test case is copper modeled by the Mie-Gr{\"u}neisen equation of state with the parameters given in Table \ref{tab:dataEOS}. In this case the Yield stress is set to $\sigma_0=9\cdot10^{-4}$ and we consider $c_v=1.0$. The loading behavior of the material is described by the relation \eqref{eqn.tauS} with $\tau_0=0.1$ and $n=10$. The initial density and pressure correspond to the reference values and the initial velocity field is zero, while the distortion tensor is simply set to $\AAA=\II$. The computational domain is initially given by $\Omega(0)=[0;1.5]\times[0;0.1]$ and it is discretized with a characteristic mesh size of $h=1/300$ with $N_E=13248$ triangles. The left boundary condition is a piston of velocity $\v_c=(0.002,0)$ while the other boundaries are treated as no-slip walls. The final computational time is $t_f=1.5$ and we use a third order accurate ADER-WENO-ALE scheme to obtain the results plotted in Figure \ref{fig:EPP}, where we compare the density and the horizontal velocity profiles 
against an available exact solution for the model of ideal plasticity. 
Note that for this test problem we do not have an exact solution of the HPR model used for the numerical simulation, but there is only an exact solution available for the model of ideal plasticity, 
with rate independent yield stress.  It is therefore not easy to make a direct comparison, because only elastic precursors are discontinuous in both approaches, while the plastic wave 
is continuous in the HPR model, see e.g. the paragraph 17 in the book \cite{Godunov:2003a} for a more detailed discussion on this topic. 
However, from the results presented in Figure \ref{fig:EPP} we can observe that the main waves and plateaus are still well reproduced, despite the use of different mathematical models 
in the numerical simulation and for the exact solution. 

\subsection{Elastic vibrations of a beryllium plate} \label{sssec:beryllium}
\begin{figure}[!htbp]
	\begin{center}
	\begin{tabular}{ccc} 
	\includegraphics[width=0.33\textwidth]{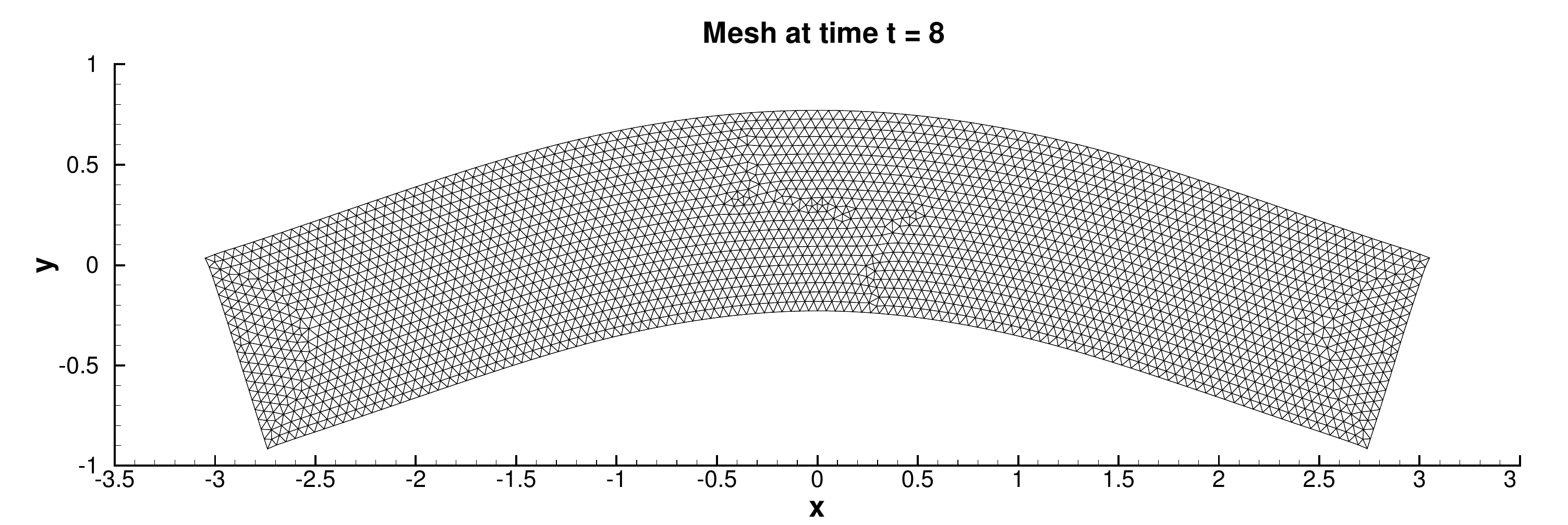} &
	\includegraphics[width=0.33\textwidth]{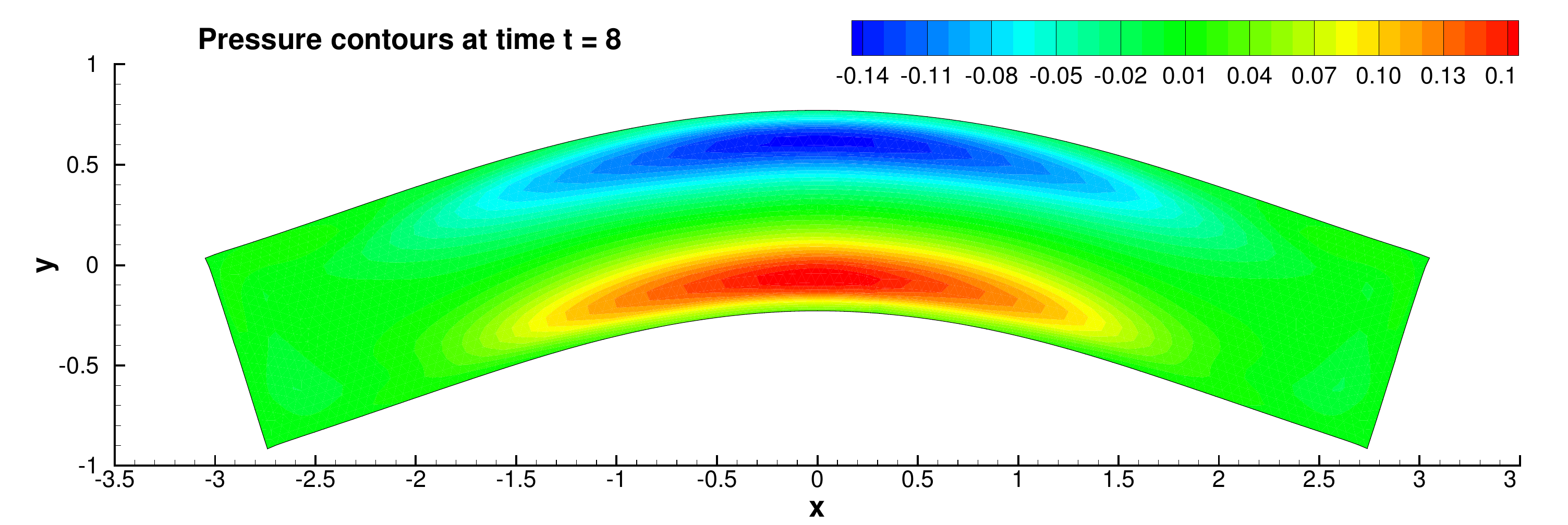} &
	\includegraphics[width=0.33\textwidth]{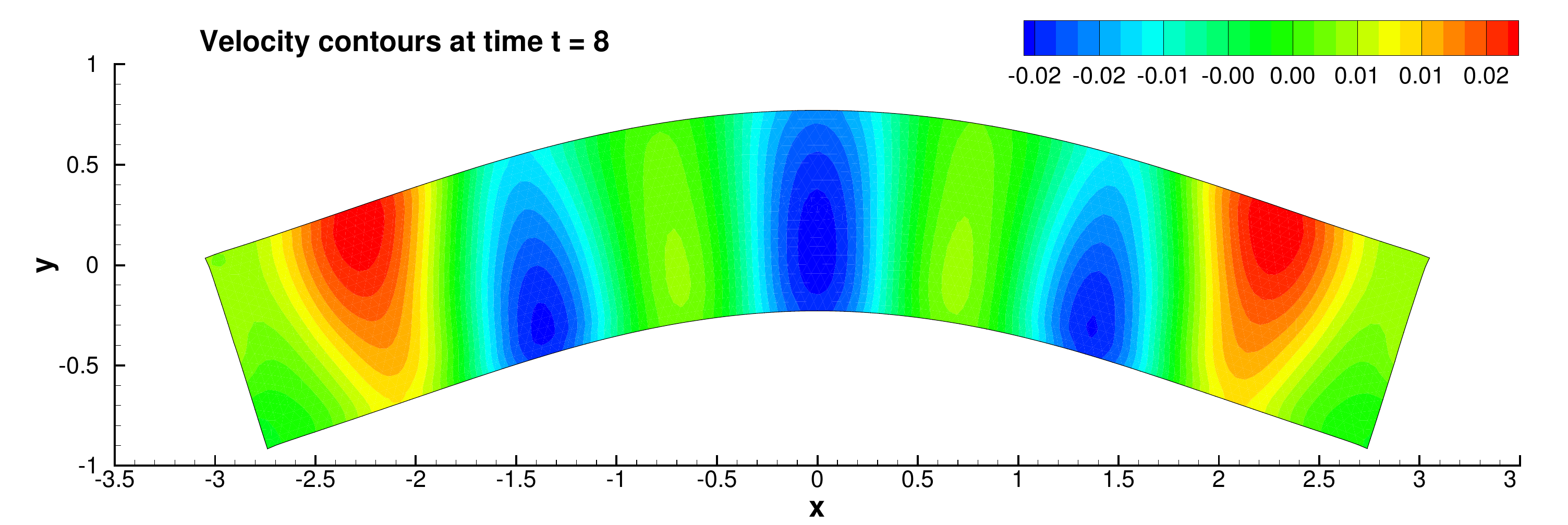} \\
	\includegraphics[width=0.33\textwidth]{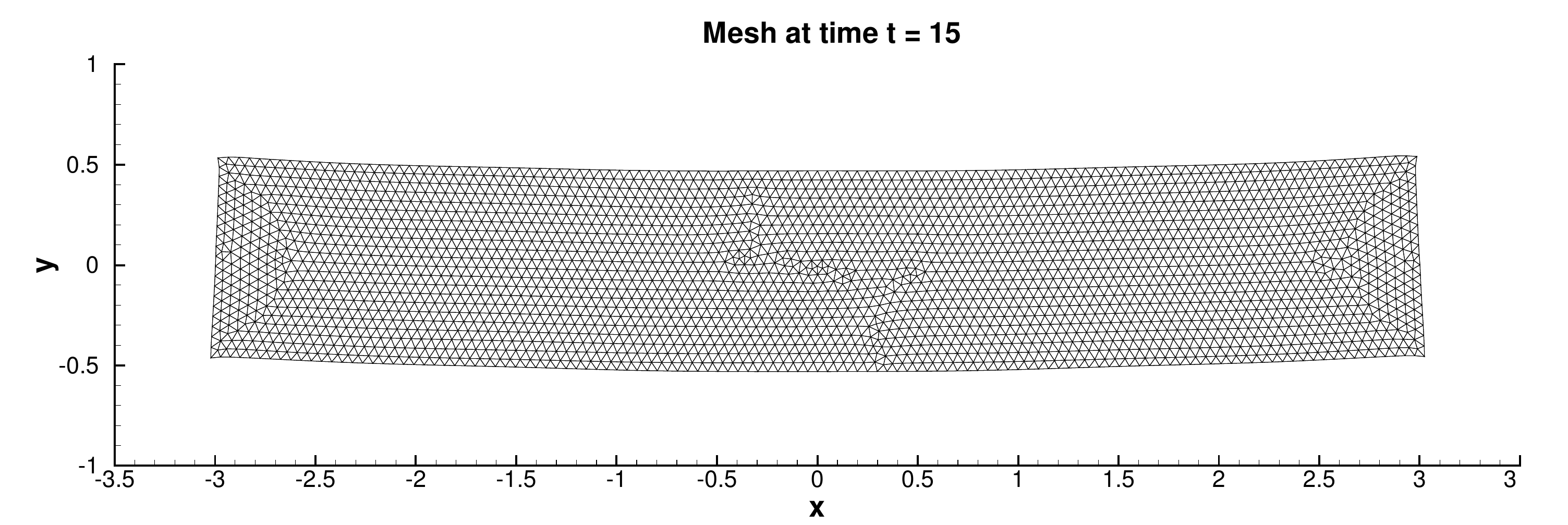}&
	\includegraphics[width=0.33\textwidth]{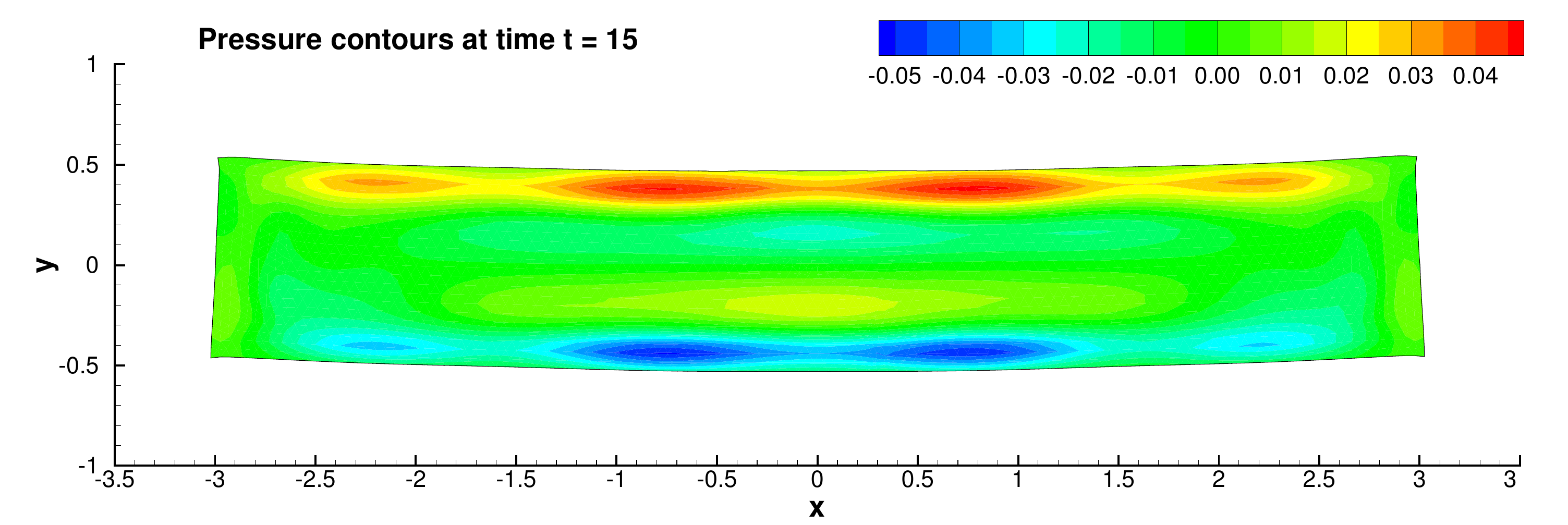}&
	\includegraphics[width=0.33\textwidth]{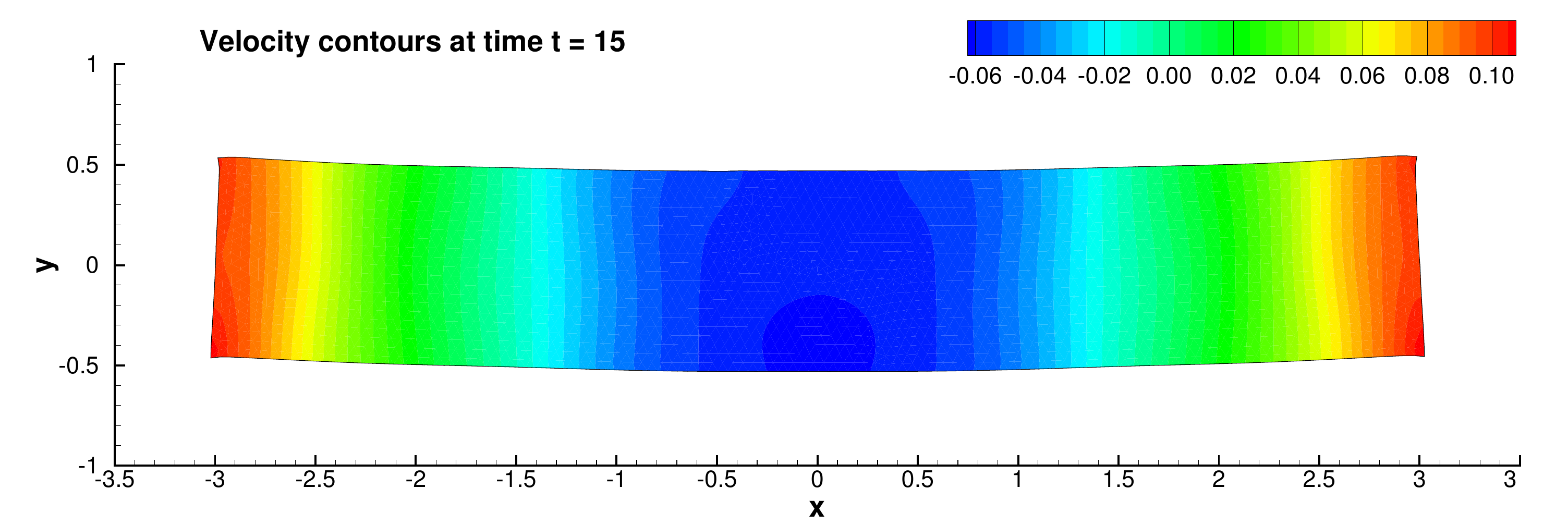}\\
	\includegraphics[width=0.33\textwidth]{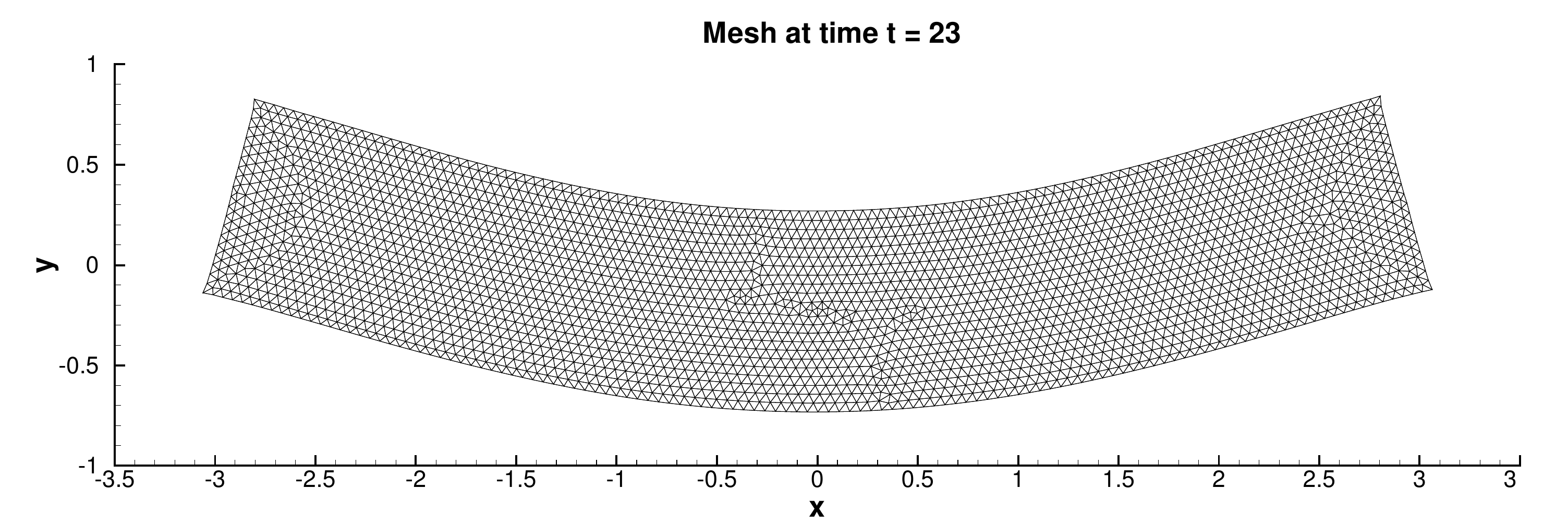} &
	\includegraphics[width=0.33\textwidth]{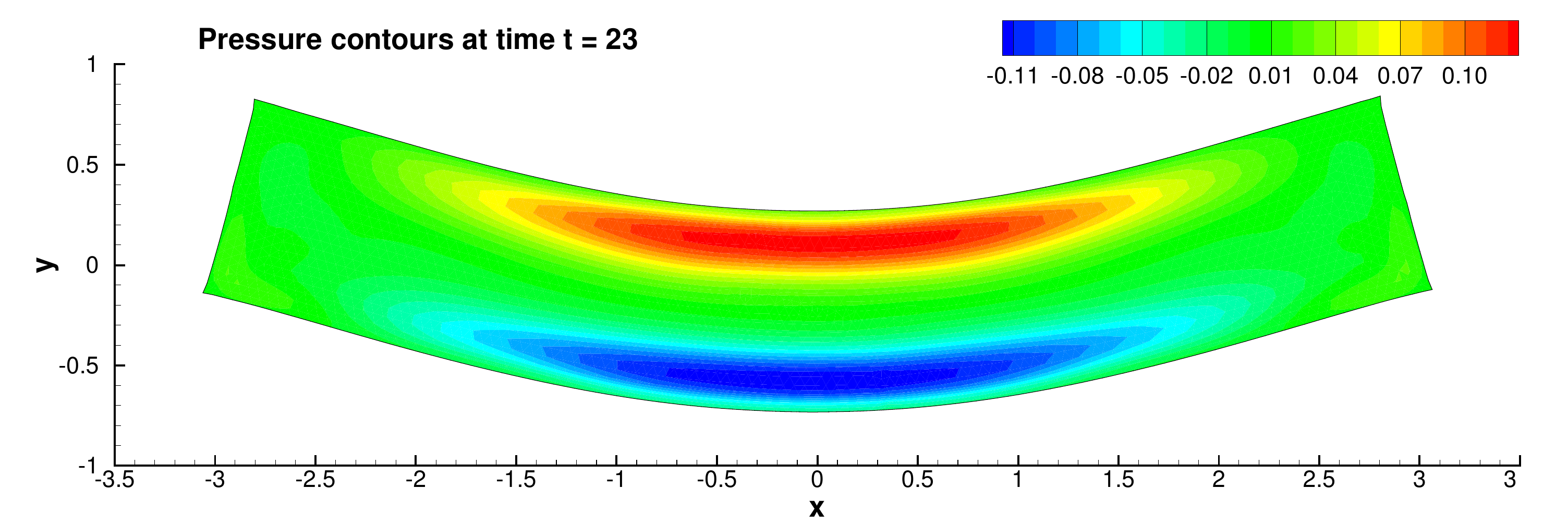} &
	\includegraphics[width=0.33\textwidth]{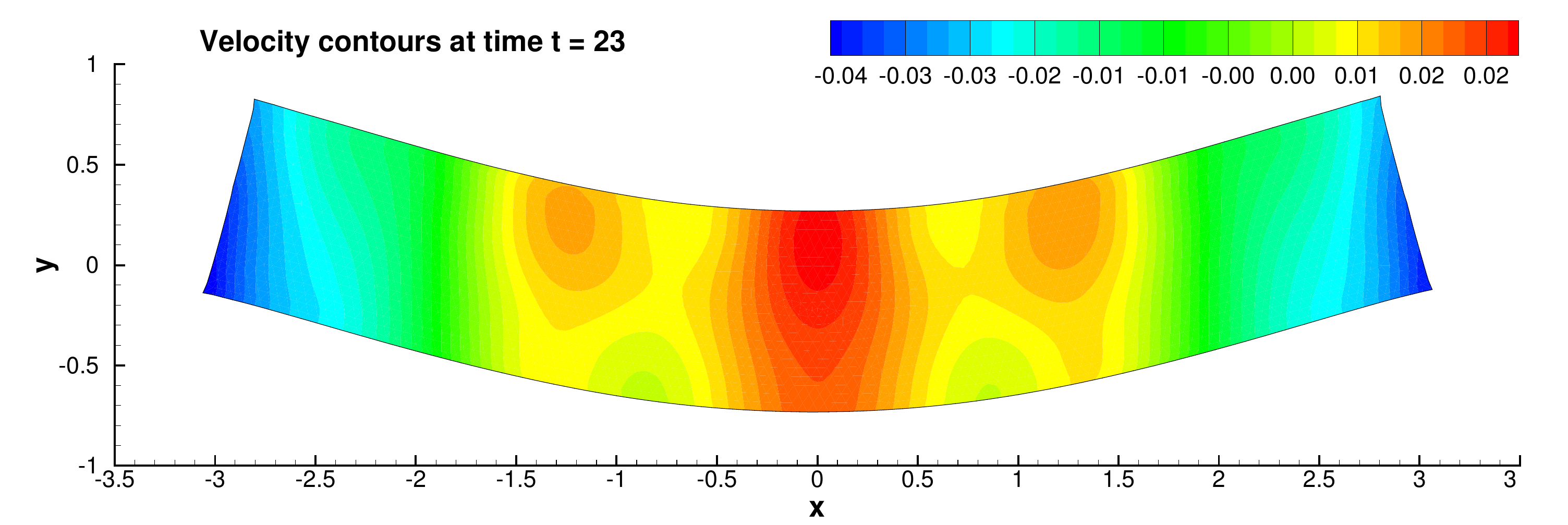} \\
	\includegraphics[width=0.33\textwidth]{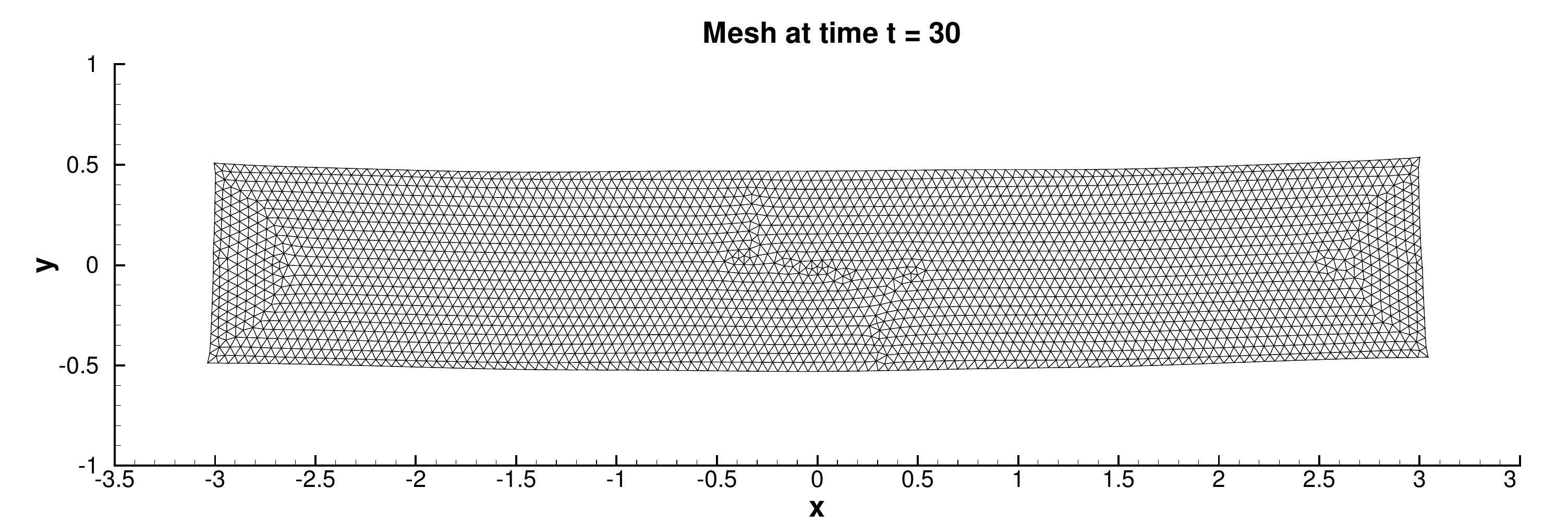} &
	\includegraphics[width=0.33\textwidth]{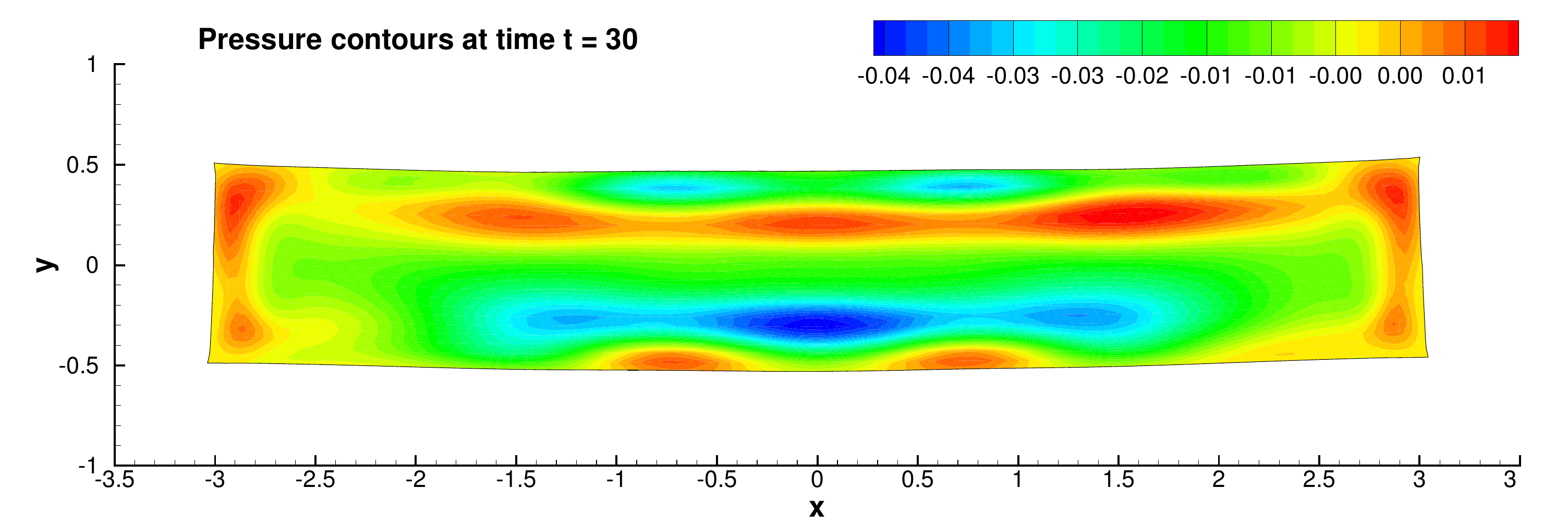} &
	\includegraphics[width=0.33\textwidth]{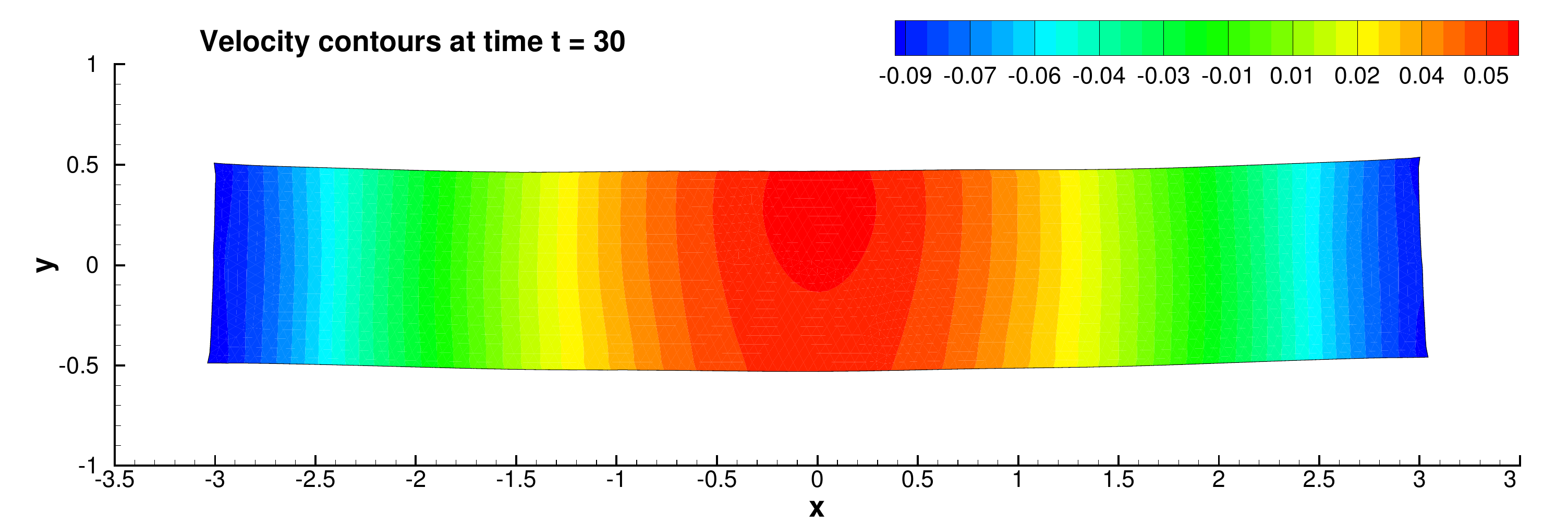} 
	\end{tabular}
	\caption{Results for the beryllium bar test case at output times $t=8$, $t=15$, $t=23$ and $t=30$ (from top to bottom). Left: mesh configuration. Middle: pressure distribution. Right: vertical component of the velocity.} 
	\label{fig:Be}
	\end{center}
\end{figure}

This problem simulates the elastic (reversible) vibrations of a beryllium bar after an initial velocity impulse \cite{Sambasivan_presentation}. The beryllium plate is characterized by the constants given in Table \ref{tab:dataEOS}. The initial computational domain is $\Omega(0)=[-3;3]\times[-0.5;0.5]$ and the computational grid counts $N_E=5338$ control volumes with $h=0.005$. Free traction boundary conditions are imposed everywhere as explained in Section \ref{sec.BCs} and the bar is initially assigned with the reference density and pressure, see Table \ref{tab:dataEOS}, with the distortion tensor $\AAA=\II$ and the velocity field $\v=(0,v(x))$, where the initial vertical velocity $v(x)$ is given by
\begin{equation}
v(x) = A \omega\left\{  
   C_1 \left( \sinh(\Omega(x+3))+\sin(\Omega(x+3)) \right)
 - S_1\left( \cosh(\Omega(x+3))+\cos(\Omega(x+3)) \right)
    \right\},
\end{equation}
with $\Omega=0.7883401241$, $\omega=0.2359739922$, $A=0.004336850425$, $S_1=57.64552048$ and $C_1=56.53585154$. The final time is set to $t_f=53.25$ according to \cite{Burton2015} such that it corresponds to two complete flexural periods $\omega$. At this time, the bar returns back to its original position for the second time. Furthermore the deformation should not generate any irreversible plastic transition in the beryllium, that means that the Yield stress must never be exceeded throughout the entire computation. The parameters for evaluating the relaxation time $\tau_1$ in \eqref{eqn.tauS} are $\tau_0=10$ and $n=1$. In Figure \ref{fig:Be} we present the mesh configuration, the pressure and the vertical velocity component respectively on left, middle and right panels for intermediate times $t=8$, $t=15$, $t=23$ and $t=30$ which cover approximately one flexural period. Please note that the color scales for the pressure are different depending if the bar is back to its original position or not. Qualitatively the bar is behaving as expected and these third order accurate results visually compare well against known results from other Lagrangian schemes \cite{Sambasivan_13,Burton2015}.

\begin{figure}[!htbp]
	\begin{center}
	\includegraphics[width=0.45\textwidth]{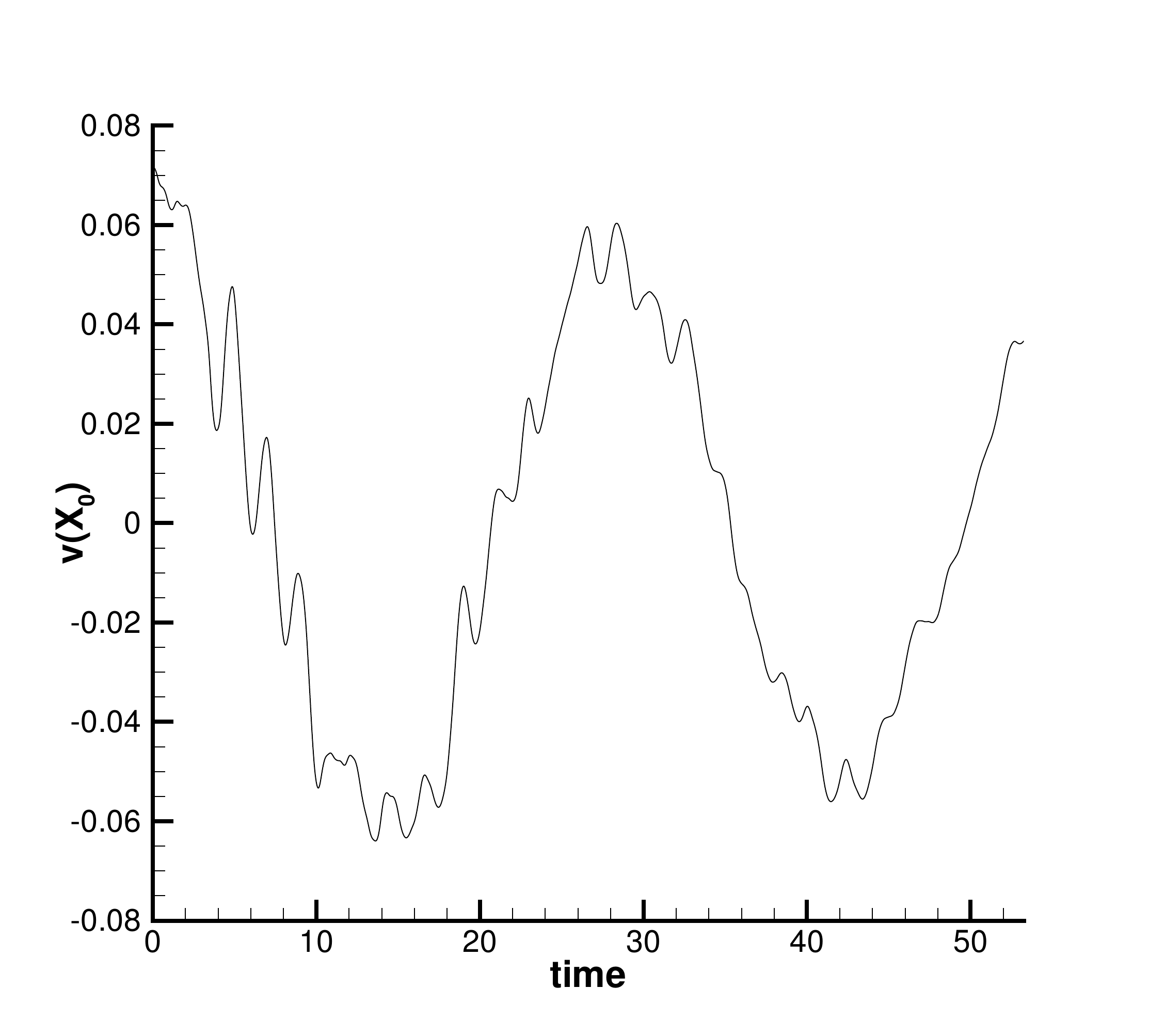} \\
	\caption{Beryllium bar test case: vertical velocity component of the point initially located at at $X_0=(0,0)$.} 
	\label{fig:Be_pick}
	\end{center}
\end{figure}

As noticed in \cite{Burton2015,Sambasivan_13} the observed oscillation period is of about $\omega_n=30$ instead of the theoretical one of $\omega_e=26.6266$, so our results are in agreement with what already obtained in literature. Finally, Figure \ref{fig:Be_pick} shows the time evolution of the vertical component of the velocity of the mesh point originally located at $X_0=(0,0)$, i.e. the barycenter of the bar. Again the results are in excellent agreement with the same plot reported in \cite{Burton2015}.

\subsection{Taylor bar impact} \label{sssec:taylor}
\begin{figure}[!htbp]
	\begin{center}
	\begin{tabular}{ccc} 
	\hspace{-1.5cm}
	\vspace{-0.15cm}
	\includegraphics[width=0.35\textwidth]{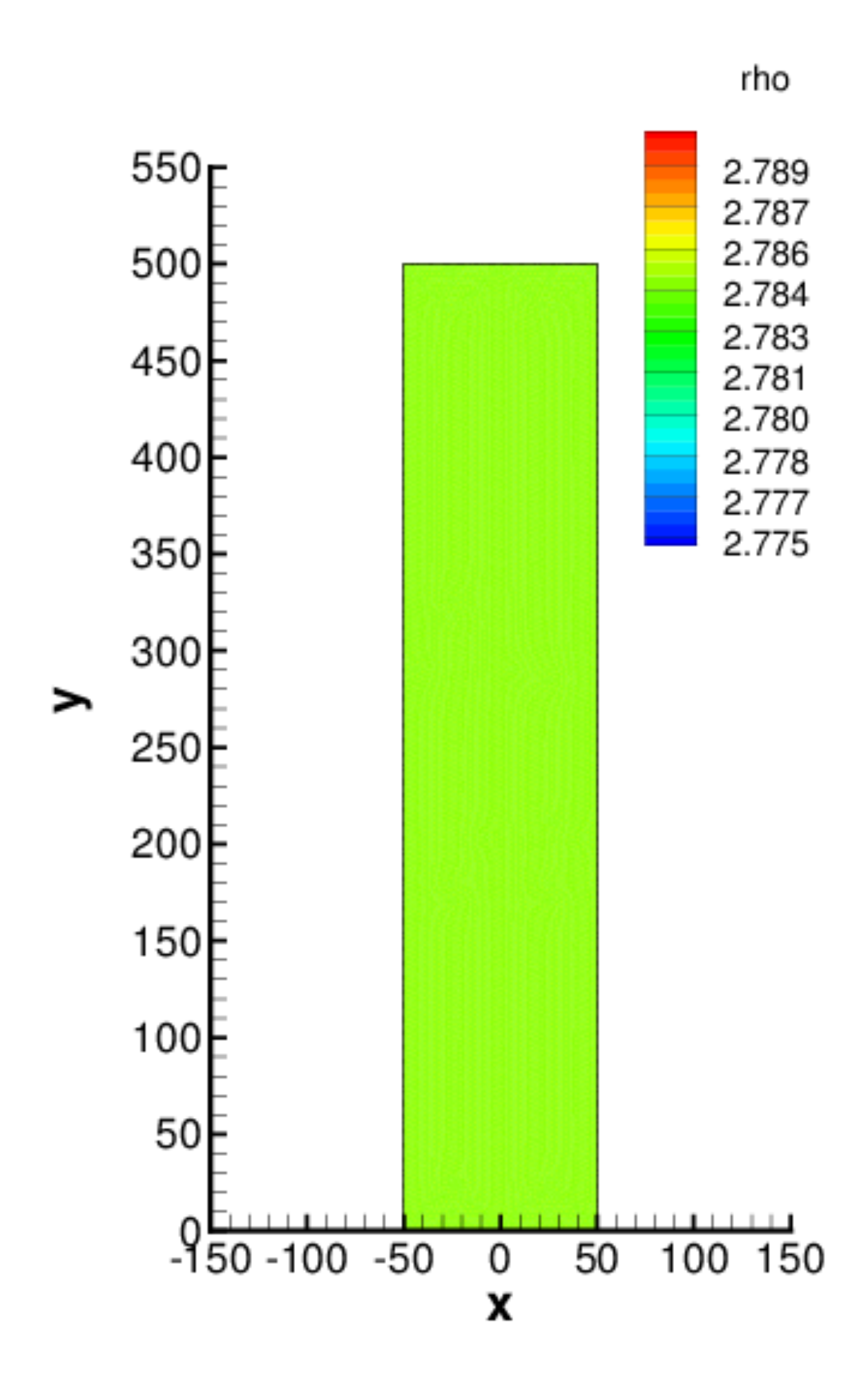} &
	\vspace{-0.15cm}
	\includegraphics[width=0.35\textwidth]{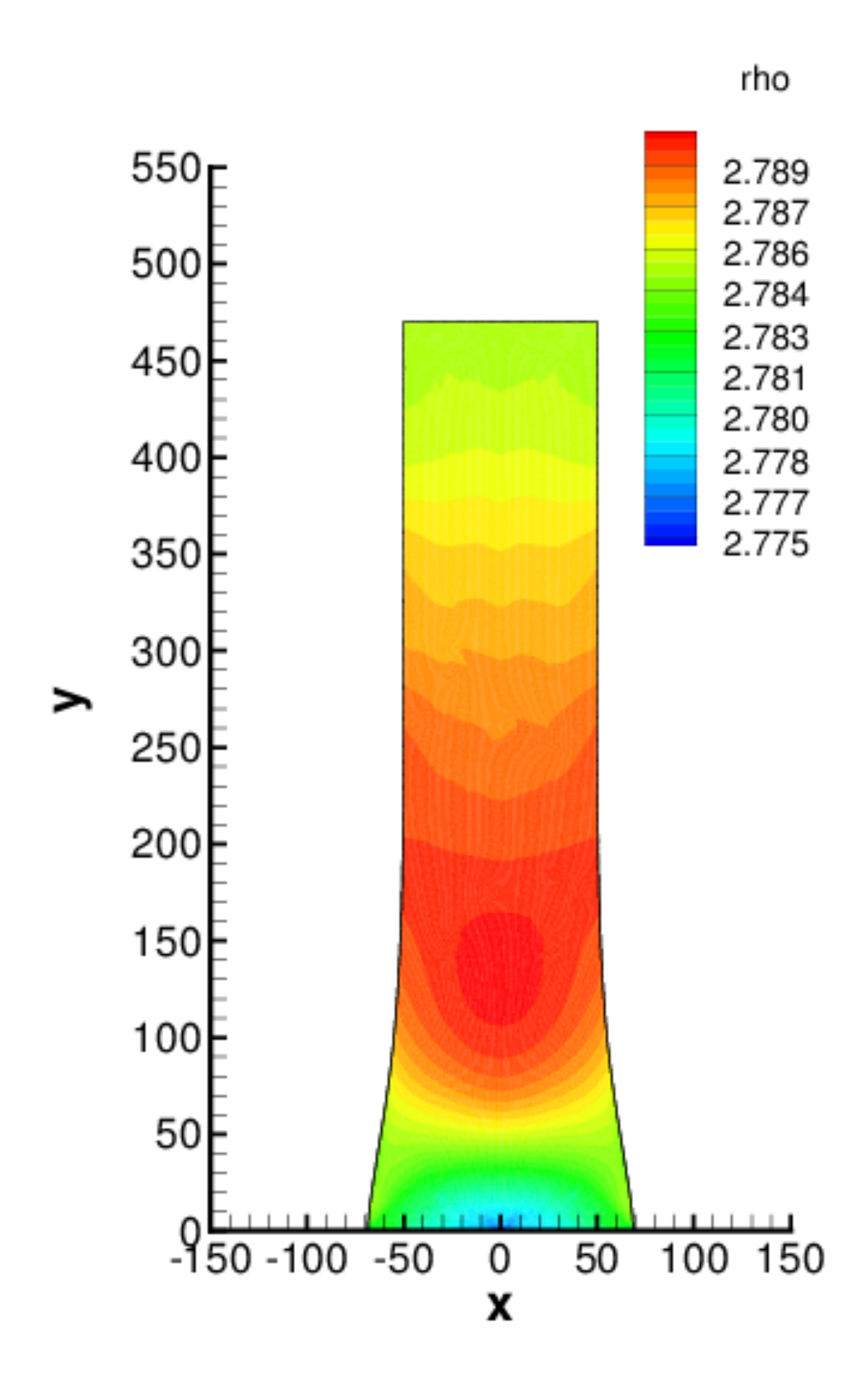} &
	\vspace{-0.15cm}
	\includegraphics[width=0.35\textwidth]{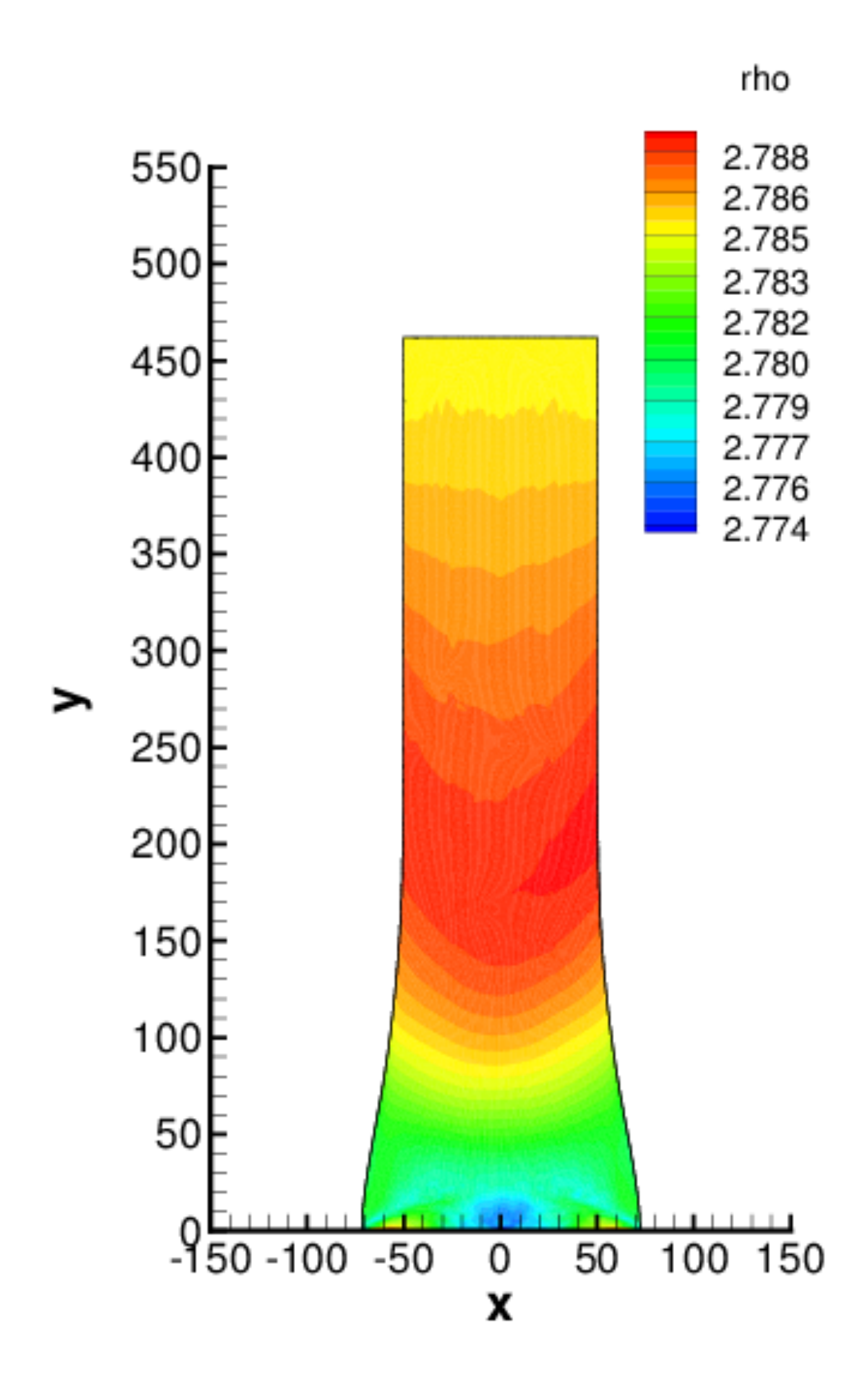} \\
	\hspace{-1.5cm}
	\includegraphics[width=0.35\textwidth]{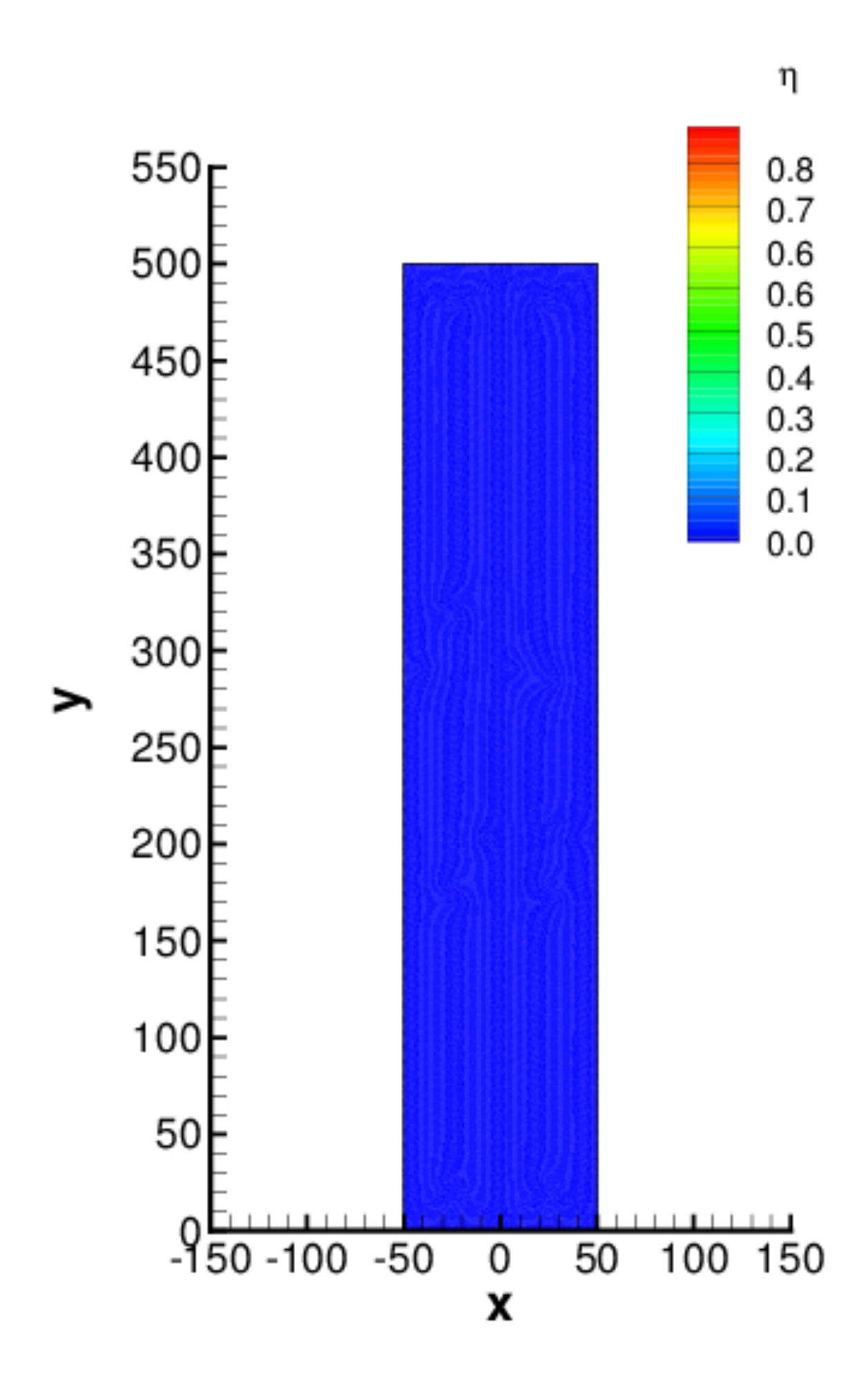} &
	\includegraphics[width=0.35\textwidth]{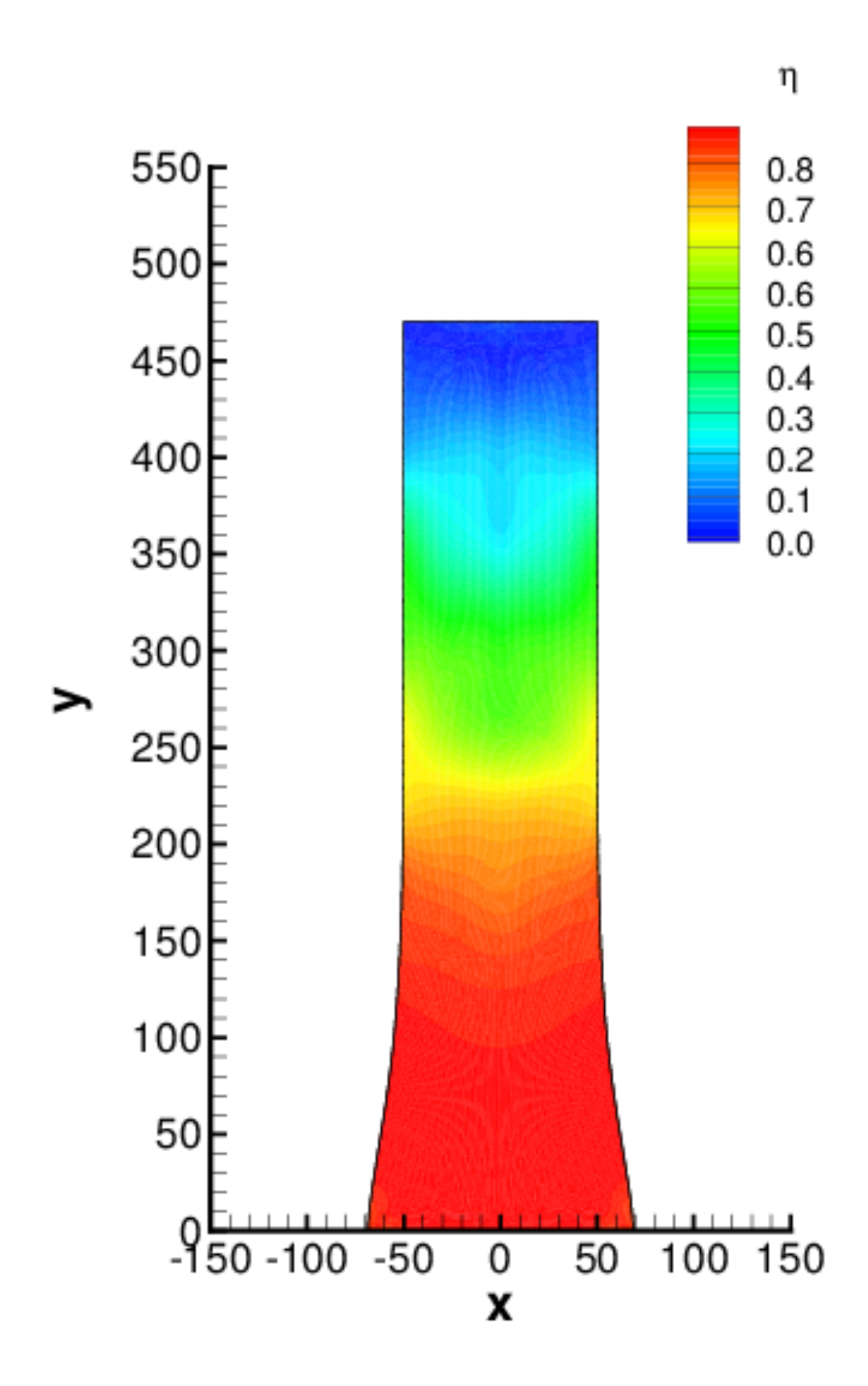} &
	\includegraphics[width=0.35\textwidth]{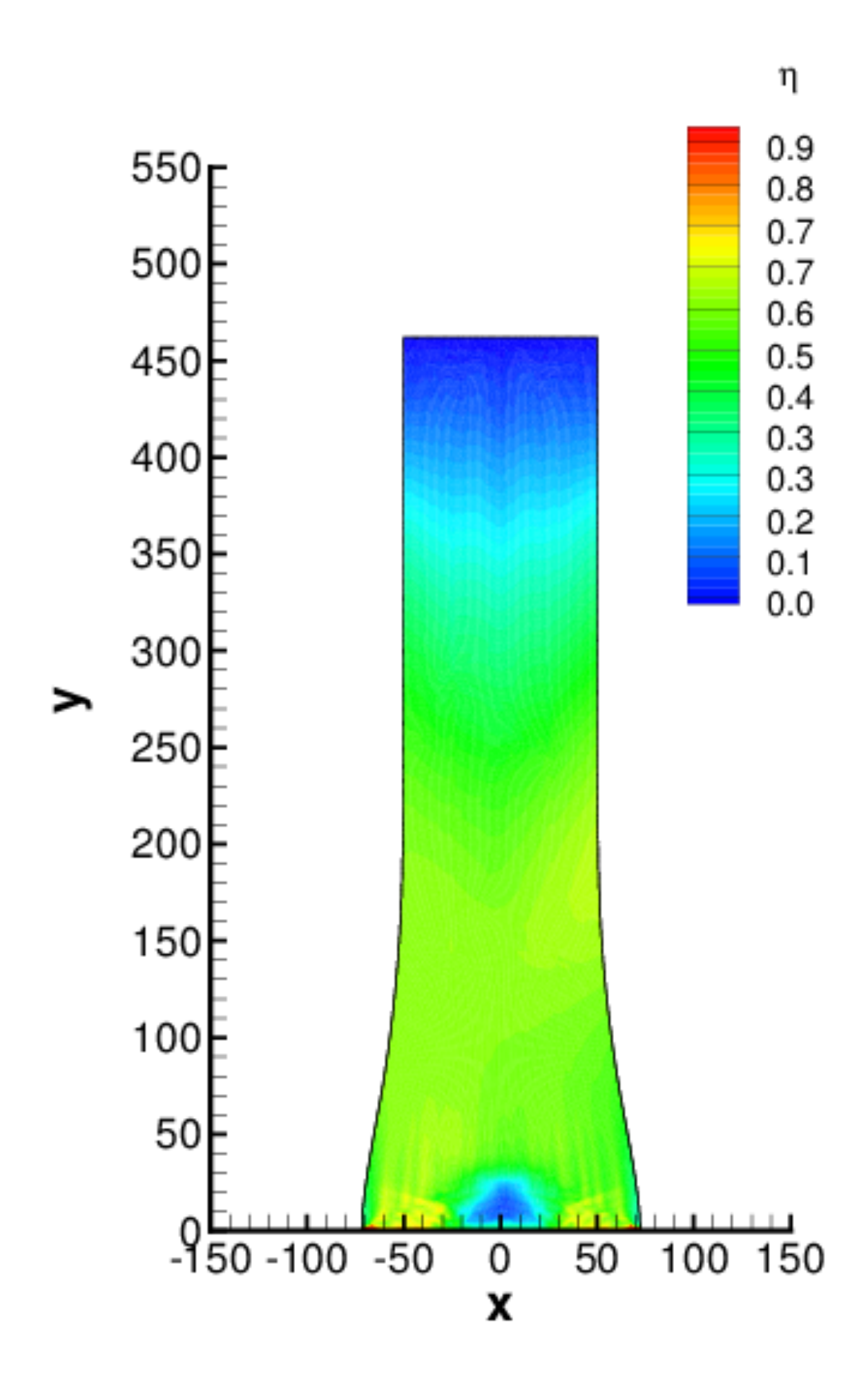} \\
	\end{tabular}
	\caption{Taylor bar impact problem: density distribution (top) and plastic deformation (bottom) at output times $t=0$, $t=0.0025$ and $t=0.005$.} 
	\label{fig:Taylor_bar}
	\end{center}
\end{figure}

\begin{figure}[!htbp]
	\begin{center}
	\begin{tabular}{cc} 
	\includegraphics[width=0.47\textwidth]{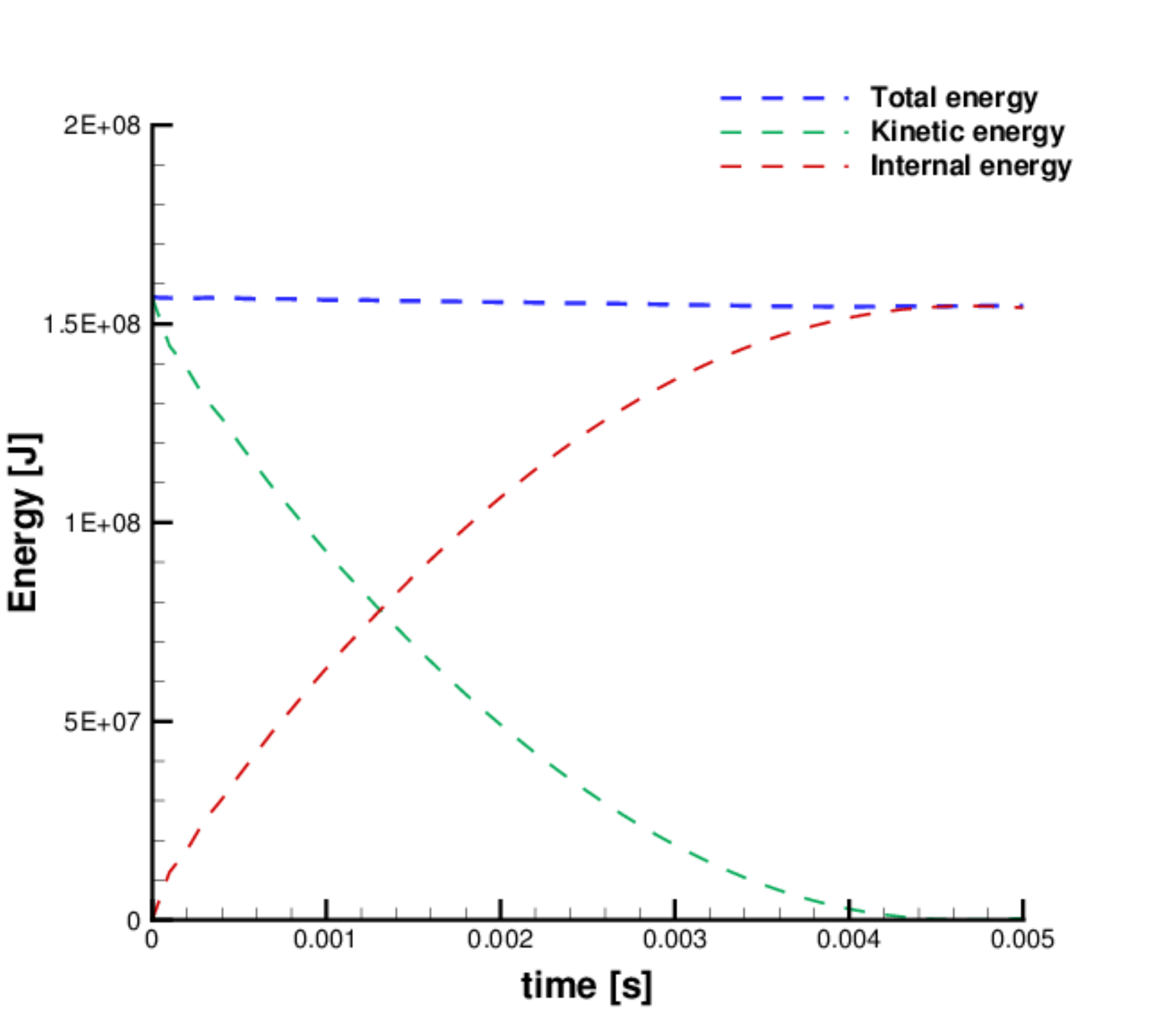} &
	\includegraphics[width=0.47\textwidth]{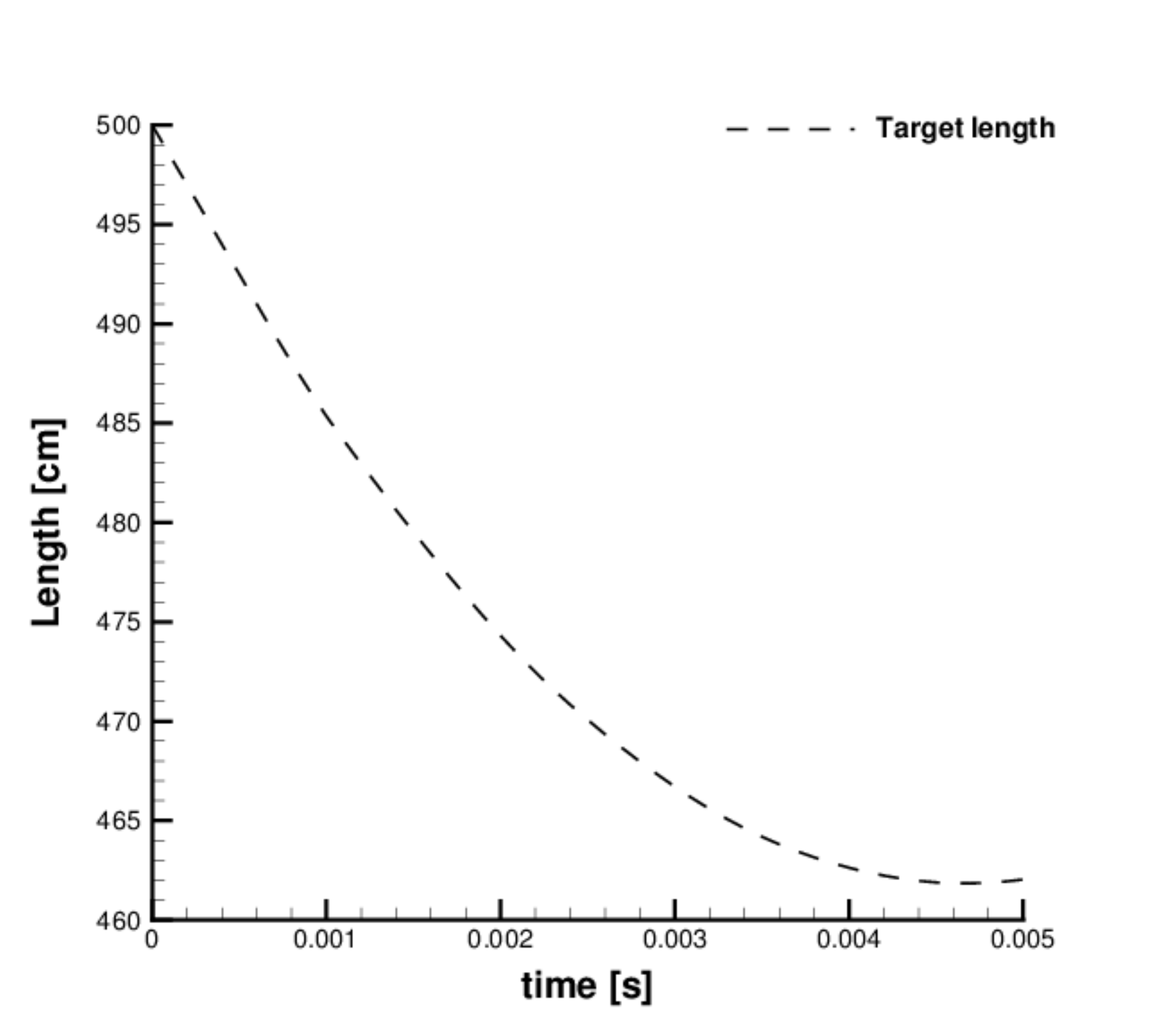} \\
	\end{tabular}
	\caption{Left: balance of kinetic, total and internal energy for the Taylor bar impact problem. Right: Length of the target versus time.} 
	\label{fig:Taylor_energy}
	\end{center}
\end{figure}

\begin{figure}[!htbp]
	\begin{center}
	\includegraphics[width=1.00\textwidth]{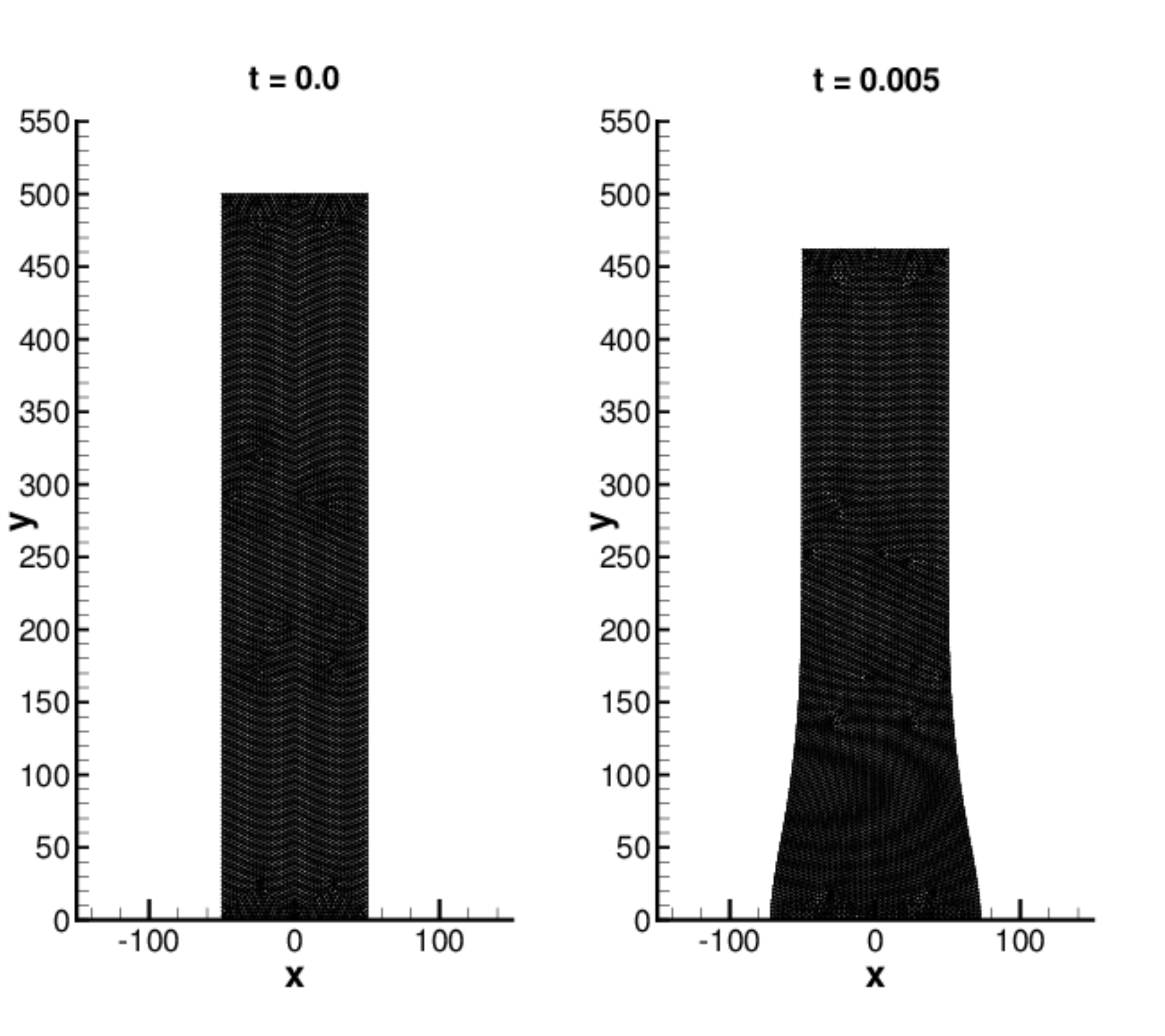} \\
	\caption{Taylor bar impact problem: initial (left) and final (right) mesh configuration.
	} 
	\label{fig:Taylor_mesh}
	\end{center}
\end{figure}

The Taylor bar impact is a classical test of an elasto-plastic target that impacts on a rigid solid wall \cite{Sambasivan_13,Maire_elasto_13,Sambasivan_presentation,Dobrev2014}. In this work we consider an aluminum bar with the initial length $L=500$ and thickness $H=100$. The parameters of the material are given in Table \ref{tab:dataEOS} and the target is initially moving with velocity $\v=(0,-0.015)$ towards a wall located at $y=0$. The initial condition is chosen as $\rho=\rho_0$, $p=p_0$, $\mathbf{A}=\mathbf{I}$ with the parameters $\tau_0=1$ and $n=20$ for the computation of the relaxation time \eqref{eqn.tauS}. We set free traction boundary conditions everywhere apart from the bottom boundary which is treated with a wall-type boundary condition. According to \cite{Maire_elasto_13,Dobrev2014} the final time of the simulation is $t=0.005$ and the computational domain is discretized with a total number of $N_E=12720$ triangles, corresponding to a mesh size of $h=3$. Here we adopt a classical source splitting for the treatment of the stiff sources that arise from the plastic deformation induced by the motion of the target. In Figure \ref{fig:Taylor_bar} we present the results computed with a third order accurate ALE ADER-WENO scheme with an Osher-type numerical flux \cite{LagrangeNC,Lagrange3D} which is less dissipative than the Rusanov flux \eqref{eqn.rusanov}: we plot the density distribution as well as the plastic rate $\eta=\frac{\sigma_I}{\sigma_0}$ at output times $t=0$, $t=0.0025$ and $t=0.005$. We note that the numerical solution is reasonably in agreement with what presented in \cite{Maire_elasto_13}, even though the models used are quite different. Furthermore during the impact the kinetic energy is totally dissipated into internal energy, as clearly shown in Figure \ref{fig:Taylor_energy}, and such a behavior has been observed also in \cite{Maire_elasto_13,Dobrev2014}. Finally, Figure \ref{fig:Taylor_mesh} depicts the initial and final mesh configurations, while the evolution of the target length is given in Figure \ref{fig:Taylor_energy} and we measure a final length of $L_f=462$ which perfectly fits the result achieved in \cite{Maire_elasto_13}.  

\section{Conclusion and Perspectives} \label{sec:conclusion}

The purpose of this paper was the numerical solution of the unified first order hyperbolic 
Peshkov \& Romenski \cite{PeshRom2014} (HPR) model of continuum mechanics, using a  
multi-dimensional ADER-WENO scheme on moving meshes in the direct ALE framework \cite{LagrangeNC,Lagrange2D,Lagrange3D}. 
The appealing property of the HPR model, which derives from the theory of nonlinear \textit{hyperelasticity} 
established by Godunov \& Romenski in \cite{GodunovRomenski72,Godunov:2003a}, is its 
ability to describe the behavior of inviscid and viscous compressible Newtonian and non-Newtonian 
fluids with heat conduction, and, at the same time, the behavior of elastic and elasto-plastic solids. 
In this paper we have shown that the family of high order ADER-WENO-ALE schemes can solve the complex 
governing PDE system of the HPR model in two limiting cases of the model, namely in the limit of 
inviscid and viscous Newtonian fluids, as well as in the limit of nonlinear hyperelasticity 
for elastic and elasto-plastic solids.  
In both cases the numerical results are comparable with results obtained from established
standard models, namely the Euler or Navier-Stokes equations for fluids, or the classical 
hypo-elastic model with plasticity \cite{Wilkins64}. 
More precisely, we have shown numerically that for pure fluid flow problems, the numerical scheme 
can achieve optimal order of accuracy for smooth flow, maintaining an essentially non oscillatory behavior 
in the presence of shock waves and steep fronts. On the other hand, we have also shown that classical 
elasto-plastic test cases can be simulated both in the reversible elastic regime (beryllium plate) 
or in situations where transition from elastic to plastic material behavior occurs (Taylor bar impact). 
Both limits of the HPR model can be nicely simulated by our high order one-step ADER-WENO-ALE schemes.  
Together with the computational results shown in \cite{Dumbser_HPR_16} the family of ADER-WENO schemes seems therefore 
to be a very promising candidate to simulate the full range of possible intermediate models embedded into  
the HPR formulation.
The moving mesh technique used in this paper is appealing when dealing with solid materials surrounded by fluids or gas, 
consequently we plan in the near future to test such situations, also adopting the idea of diffuse interface methods
as outlined in \cite{FavrGavr2012,GavrFavr2008,favrie2009solid}. 
We also plan to replace the WENO stabilization technique by the \textit{a posteriori} MOOD paradigm, see 
\cite{CDL3,ALE-MOOD2,ADER_MOOD_14} and its extension to the discontinuous Galerkin framework recently forwarded 
in \cite{Dumbser2014,Zanotti2015a}. 
Also the treatment of boundary conditions needs to be mathematically analyzed in more detail. 
Moreover we plan to explore even further the capability of the HPR model and compare with existing experimental data when possible.

\section*{Acknowledgments}
The authors would like to warmly acknowledge the help provided by I. Peshkov and E. Romenski  
for the design of proper boundary conditions and the computation of the relaxation time $\tau_1$ 
for the HPR model in the case of elasto-plastic solids. We also would like to 
thank S. Gavrilyuk for the inspiring discussions about the subject of hyperelasticity. 
Last but not least, the authors are grateful to S.K. Godunov for his great seminal ideas 
that are at the basis of the theoretical and numerical framework used in this paper. 

M.D. and W.B. have been financed by the European Research Council (ERC) under the
European Union's Seventh Framework Programme (FP7/2007-2013) with the
research project \textit{STiMulUs}, ERC Grant agreement no. 278267.
The authors acknowledge PRACE for awarding access to the SuperMUC 
supercomputer based at the LRZ in Munich, Germany. \\


\bibliography{biblio}
\bibliographystyle{plain}

\end{document}